\documentclass[preprint,11pt]{elsarticle}
\usepackage{amssymb}
\usepackage{mathrsfs}
\usepackage[OT1]{fontenc}
\usepackage{amsthm}
\usepackage{graphics}
\usepackage{epsfig}
\usepackage{booktabs}
\usepackage{tikz}
\usepackage{circuitikz}
\usepackage{tikz-cd}
\usetikzlibrary{shadings,patterns}
\usepackage{graphicx,graphics}
\usepackage{amssymb}
\usepackage{graphicx}
\usepackage{subfigure}
\usepackage{color}
\usepackage{comment}
\usepackage{tikz}
\usepackage{tikz,pgfplots}
\pgfplotsset{compat=1.12,axis lines=center}
\usepgfplotslibrary{fillbetween}
\usetikzlibrary{patterns}
\usepackage{hyperref}
\usepackage{a4wide}
\usepackage{amsmath}
\usepackage{amsfonts}
\usepackage{enumerate}
\usepackage{multicol}

\def\red#1{{\color{red}#1}}

\DeclareMathOperator{\ord}{ord}
\newcommand{\gen}[1]{\langle #1\rangle}

\makeatletter
\def\ps@pprintTitle{%
	\let\@oddhead\@empty
	\let\@evenhead\@empty
	\let\@oddfoot\@empty
	\let\@evenfoot\@oddfoot
}
\makeatother

\hypersetup{
	colorlinks=true,
	linkcolor=blue,
	filecolor=dviolet,      
	urlcolor=blue,
} 

\newtheorem*{theoremaux}{Theorem \theoremauxnum}
\gdef\theoremauxnum{1}

	\newtheorem{lemma}{\bf Lemma}[section]
	
	\newtheorem{theorem}{\bf Theorem}[section]
    \newtheorem{observation}{\bf Observation}[section]
	\newtheorem{proposition}[lemma]{\bf Proposition}
	\newtheorem{corollary}[lemma]{\bf Corollary}
	\newtheorem{definition}{\bf Definition}[section]
	\newtheorem{remark}{\bf Remark}[section]
	
	\newtheorem{problem}{\bf Problem}[section]
	
\journal{~}

\begin{document}

	\begin{frontmatter}
		
		
		
		\title{\bf On the difference of the non-generating graph and the intersection power graph of a finite group}

		
		%
		
 \author{Sudip Bera}
		\ead{sudip\_bera@dau.ac.in} 
		\address{Faculty of Mathematics,\\ 
		Dhirubhai Ambani University, Gandhinagar \\ India} 
  \cortext[cor1]{Corresponding author}
 \author{Sanjay Mukherjee}
       \ead{sanjaymukherjee@hri.res.in} 
		\address{Harish-Chandra Research Institute, \\ A CI of Homi Bhabha National Institute, \\ Chhatnag Road,
Jhunsi, Prayagraj - 211019, India\\}

		%

		\begin{abstract}
        For a finite group $G$, the non-generating graph ${\rm NGen}(G)$ is a simple graph with vertex set $G$, and two distinct vertices $x,y$ are adjacent in ${\rm NGen}(G)$ if the subgroup $\langle x,y \rangle$ is a proper subgroup of $G$. The intersection power graph $\mathcal{G}_I(G)$ is a simple graph with vertex set $G$, and two distinct vertices $x,y$ are adjacent in $\mathcal{G}_I(G)$ if $\langle x\rangle\cap \langle y\rangle \neq \{e\} $, and the identity element $e$ is adjacent to all other vertices of $G$. A finite group $G$ is called a $N_I$-group, if the graph $\mathcal{G}_I(G)$ is a spanning subgraph of ${\rm NGen}(G)$. For a $N_I$-group $G$, the difference graph $\mathcal{DNI}(G)$ is defined to be the finite simple graph with vertex set $G\setminus\{e\}$, and two distinct vertices $x,y$ are adjacent in $\mathcal{DNI}(G)$ if $x$ and $y$ are adjacent in ${\rm NGen}(G)$ but not adjacent in $\mathcal{G}_I(G)$.

        In this article, we study the difference graph $\mathcal{DNI}(G)$ of a finite group $G$. We first investigate the large class of $N_I$-groups, and classify $N_I$-groups for some classes of finite groups. Next, we classify finite groups $G$ for which $\mathcal{DNI}(G)$ is a complete graph or a null graph. In addition, we study the difference graph of the direct product of finite groups of mutually coprime orders, and show that the class of difference graphs of finite groups is a universal family of finite graphs, using the results obtained. Then we classify all finite non-$2$-generated groups and $2$-generated nilpotent groups $G$, such that $\mathcal{DNI}(G)$ is connected. We also explore the diameter of $\mathcal{DNI}(G)$ is the above cases. Further, we classify finite non-$2$-generated groups and $2$-generated nilpotent groups $G$, such that $\mathcal{DNI}(G)$ is Eulerian. In addition to that, we delve into various forbidden subgraph properties of the difference graph, including claw-free, chordal, and planarity conditions, and classify all finite non-$2$-generated groups and $2$-generated nilpotent groups $G$, such that $\mathcal{DNI}(G)$ satisfy these properties. We also study these properties for finite simple groups.
        \end{abstract}
		
		
	\end{frontmatter}

\section{Introduction}
The interplay between algebraic structures and graphs has emerged as a powerful paradigm for extracting and visualizing intrinsic properties of finite groups. Starting from the very idea of Cayley graph, several other notions-such as commuting graph, power graph, enhanced power graph, deep commuting graph etc.~encode subtle combinatorial and algebraic information through adjacency relations governed by cyclic subgroups, commutativity, Schur cover etc. The notion of commuting graph of a group was introduced by Brauer and Fowler \cite{braurflower} in $1955$, while examining the finiteness of the type of groups of even order whose centralizer follow certain conditions. Since then, several combinatorial properties of commuting graphs, as well as structural properties of the group via its commuting graph have been well explored \cite{arvind, PJC-Graphs-def-grp, eberhard, guralnick} in the literature. Add to this, Kelarev and Quinn \cite{combinatorialpropertyandpowergraphsofgroupskq1} introduced the notion of directed power graph of a group. The underlying undirected graph, simply termed the `power graph' of that group, was first considered by Chakrabarty et al.~in \cite{undpwrgraphofsemgmainsgc1}. For a detailed overview of several combinatorial studies on power graphs, we refer to the survey papers \cite{AKC, KSCC} and the references therein.

While each of these graphs has been studied independently and shown to reflect important structural features of the underlying group, a comparative investigation between them remains novel. To study that, the notion of enhanced power graph \cite{firstenhcedpwrstrctreaacns1} and deep commuting graph \cite{deepc} of finite groups were introduced such that they lie between the power graph and the commuting graph, and the four graphs form a one dimensional hierarchy of graphs defined over finite groups. In addition to that, difference between a pair of these graphs have also been considered in \cite{biswas, bonjak2024differencegraphsfiniteabelian, Kumar_panda}. 

Alongside these developments, the study of generating properties of finite groups has remained a central theme in finite group theory. A classical result asserts that every finite simple group is generated by two elements. Since finite simple groups serve as the building blocks of all finite groups, in much the same way that prime numbers are the building blocks of the integers, this theorem highlights the prominence of $2$-generated groups within finite group theory. More generally, the problem of understanding generating pairs in finite simple groups, and in finite groups as a whole, has attracted considerable attention and led to several landmark contributions \cite{Burness, GURALNICK2000743, Steinberg_1962}. Motivated by the early developments in this direction, the notion of the generating graph of a finite group originated in the work of Liebeck and Shalev \cite{Shalev}, and was subsequently formalized and studied in greater depth by Lucchini and Mar\'{o}ti in \cite{Lucchini_Maroti}. For a graph $\Gamma$ with vertex set $V(\Gamma)$ and edge set $E(\Gamma)$, two vertices $u,v\in V(\Gamma)$ are said to be \emph{adjacent} in $\Gamma$, if there is an edge between $u$ and $v$ in $\Gamma$. For a finite group $G$, the \textit{generating graph} of $G$, denoted by ${\rm Gen}(G)$, is a simple graph with vertex set $G$ and two distinct vertices $x$ and $y$ are adjacent in ${\rm Gen}(G)$ if $\langle x, y\rangle=G$. The complement of this graph, known as the \textit{non-generating graph}, was first introduced by Lucchini and Nemmi in \cite{Lucchini_nemmi}. In that work, the authors investigated several fundamental properties of the graph, including its connectivity and diameter. Since then, the non-generating graph of a finite group has attracted considerable attention \cite{Freedman_1, Freedman_2}, leading to a growing body of literature exploring its structural and combinatorial properties. 

Add to this, the notion of the intersection power graph of a group as a supergraph of the power graph was introduced by Bera in \cite{intersectionpwegraphb3}. For a group $G$, the \textit{intersection power graph} of $G$, denoted by $\mathcal{G}_I(G)$, is the graph with vertex set $G$ and two distinct vertices $x$ and $y$ are adjacent in $\mathcal{G}_I(G)$ if there exists a non-identity element $z\in G $ such that $x^m=z=y^n $, for some $m, n\in \mathbb{N}$, i.e.~$x\sim y$ in $\mathcal{G}_I(G)$ if $\langle x\rangle\cap \langle y\rangle \neq \{e\} $. Also, the identity element $e$ is adjacent to all other elements of $G$. Several algebraic and combinatorial properties of the intersection power graph have been extensively investigated in the literature, including its metric dimension, forbidden subgraph characterizations, and topological properties such as toroidality and projective planarity \cite{Ma_Li_Fu, Ma04052025, Ma_X}. Moreover, for a broad class of finite groups, the intersection power graph arises naturally as a spanning subgraph of the non-generating graph. For example, if $G$ is a finite group with trivial center, then any two elements whose cyclic subgroups intersect nontrivially necessarily generate a proper subgroup of $G$, since otherwise that nontrivial intersection is a central subgroup of $G$. Consequently, the intersection power graph of $G$ is a spanning subgraph of its non-generating graph. Motivated by this observation, in the present article we investigate the graph obtained by taking the difference between the non-generating graph and the intersection power graph of a finite group. We begin by introducing the following class of finite groups, for which the intersection power graph is a spanning subgraph of the non-generating graph. 

\begin{definition}\label{NI group}
A group $G$ is said to be a $N_I$-group if $\mathcal{G}_I(G)$ is a spanning subgraph of ${\rm NGen}(G)$.    
\end{definition}
Clearly, not all groups are $N_I$-group. For example, the group $C_{p}\times C_{p^2}$ is not a $N_I$-group, since the elements $(x,y)$ and $(e,y)$ generate $C_{p}\times C_{p^2}$, but $\langle (x,y)\rangle\cap\langle(e,y)\rangle\neq \{e\}$, where $x,y$ are generators of $C_p$ and $C_{p^2}$, respectively. Also, no nontrivial cyclic group is a $N_I$-group and every non-cyclic non-$2$-generated group is a $N_I$-group. Moreover, all simple groups are $N_I$-group.

It is easy to see that, for any $N_I$-group $G$ the identity element $e$ is a dominant vertex of both ${\rm NGen}(G)$ and $\mathcal{G}_I(G)$. Hence, $e$ is always an isolated vertex in the difference between ${\rm NGen}(G)$ and $\mathcal{G}_I(G)$. For this reason, we consider the vertex set of the difference graph to be $G\setminus \{e\}$, and the definition of the difference graph goes as follows.

\begin{definition}\label{Defn: Difference graph, DNI(G)}
For a $N_I$-group $G$, the difference graph of $G$, denoted by $\mathcal{DNI}(G)$ is a simple graph with vertex set $G\setminus\{e\}$, and two distinct vertices $a,b$ are adjacent in $\mathcal{DNI}(G)$ if $a\sim b$ in ${\rm NGen}(G)$ but $a\nsim b$ in $\mathcal{G}_I(G)$.
\end{definition}

It can be seen that for any finite group $G$, both ${\rm NGen}(G)$ and $\mathcal{G}_I(G)$ are invariant under any automorphism of $G$. As a consequence of this fact, we have the following simple observation.

\begin{observation}\label{Obs:DNI(G) is invariant under any automorphism}
$\mathcal{DNI}(G)$ is invariant under any automorphism of $G$.    
\end{observation}

\subsection{Organization of the article}
This article is a comprehensive study of the difference graph $\mathcal{DNI}$ of several classes of finite groups. The remainder of the article is organized as follows. In the next section, we present the necessary graph theoretic and group theoretic preliminaries. We also establish several auxiliary group theoretic results that will be used throughout the paper. Section \ref{Sec: The NI group} classifies certain classes of $2$-generated $N_I$-groups, such as finite metacyclic groups, and direct product of $2$-generated groups.

In Section \ref{Sec: Completeness and nullness}, we classify finite $N_I$-groups whose difference graph is either complete or null. Section \ref{Sec: Direct product of groups} is devoted to the study of the difference graph of direct products of finite $N_I$-groups of relatively prime orders. The results obtained therein are subsequently used in Section \ref{Sec: Universality}, where we investigate the universality of the difference graph of finite groups.

Sections \ref{Sec: Connectedness and diameter}, \ref{Sec: Eulerian}, \ref{Sec: Claw free}, \ref{Sec: Chordal}, and \ref{Sec: Planar} examine the connectedness, diameter, Eulerianity, claw-freeness, chordality, and planarity of the difference graph. In particular, we classify finite non-cyclic, non-$2$-generated groups and $2$-generated nilpotent $N_I$-groups for which the difference graph possesses these properties, and determine its diameter whenever it is connected.

In Section \ref{Sec: Difference graph of finite simple groups}, we investigate the aforementioned properties of the difference graph for several classes of finite simple groups. Finally, we conclude the article with a collection of open problems and directions for future research arising from the study of the difference graph $\mathcal{DNI}$ of finite $N_I$-groups.


\section{Preliminaries}
Throughout this paper, all graphs are finite, simple, and undirected, unless mentioned otherwise. Let $\Gamma$ be a graph with vertex set $V(\Gamma)$ and edge set $E(\Gamma)$. For $u,v\in V(\Gamma)$, we denote by $u \sim v$ the adjacency of $u$ and $v$ in $\Gamma$. The \textit{complement} of $\Gamma$, denoted by $\bar{\Gamma}$ is the graph with vertex set $V(\Gamma)$ and for two distinct vertices $u,v\in V(\Gamma)$, $u\sim v$ in $\bar{\Gamma}$ if and only if $u\nsim v$ in $\Gamma$. The \textit{neighborhood} of a vertex $v$ in $\Gamma$, denoted by $N_{\Gamma}(v)$, is the set of all vertices adjacent to $v$ in $\Gamma$. We define the \textit{degree} of $v$ in $\Gamma$ to be $\lvert N_{\Gamma}(v)\rvert$ and denote it by $d_{\Gamma}(v)$. By $N_{\Gamma}[v]$, we mean $N_{\Gamma}[v]=N_{\Gamma}(v)\cup \{v\}$. A vertex $v$ is called a \textit{dominant vertex} of $\Gamma$, if $N_{\Gamma}[v]=V(\Gamma)$. We denote the set of dominant vertices of $\Gamma$ by $_d\Gamma$. We say $\Gamma$ is complete if $_d\Gamma =V(\Gamma)$. $\Gamma$ is said to be \textit{connected}, if for any two distinct vertices $u,v\in V(\Gamma)$, there is a path between $u$ and $v$ in $\Gamma$. For a connected graph $\Gamma$, the \textit{distance} between vertices $u$ and $v$, denoted by ${\rm dist}(u, v)$, is the length of a shortest path from $u$ to $v$. The \textit{diameter} of $\Gamma$, denoted by ${\rm diam}(\Gamma$), is defined as ${\rm diam}(\Gamma)=\max\{{\rm dist}(u,v)|u,v\in V(\Gamma)\}.$ 

For a subset $S$ of $V(\Gamma)$, the subgraph induced in $\Gamma$ by $S$, or the \textit{induced subgraph} $\Gamma\vert_S$ is defined as the graph with vertex set $S$, and two vertices are adjacent in $\Gamma\vert_S$ if and only if they are adjacent in $\Gamma$. A graph is said to be \textit{chordal}, if it has no induced cycle of length more than $3$. We denote by $K_n$ the complete graph on $n$ vertices, and by $K_{a_1,a_2,\ldots,a_r}$ the complete $r$-partite graph with sizes of the parts being $a_1,a_2,\ldots,a_r$, respectively. The graph $K_{1,3}$ is called a \textit{claw}. A graph is said to be \textit{claw free}, if it does not contain any induced claw. In a graph $\Gamma$, an \textit{Eulerian trail} is a trail that contains every edge of $\Gamma$. An \textit{Eulerian circuit} or \textit{Eulerian cycle} in $\Gamma$ is an Eulerian trail that starts and ends on the same vertex of $\Gamma$. A finite graph is called \textit{Eulerian} if it contains an Eulerian circuit.

For a finite graph $\Gamma$, a \textit{planar embedding} of $\Gamma$ is a drawing of $\Gamma$ in the plane such that its edges intersect only at their endpoints. $\Gamma$ is said to \textit{planar} if $\Gamma$ admits such a planar embedding. It is easy to see that a finite graph is planar if and only if all its induced subgraphs are planar. For a connected planar graph $\Gamma$, let $F$ denote the number of regions into which a planar embedding of $\Gamma$ partitions the plane. Then the following relations are well known.
\begin{itemize}
    \item $\lvert V(\Gamma)\rvert-\lvert E(\Gamma)\rvert+F=2$.
    \item $\lvert E(\Gamma)\rvert\leq 3\lvert V(\Gamma)\rvert-6$.
\end{itemize}
\noindent Now, the next observation is a simple consequence of the above relations.
\begin{observation}\label{K_2222 not planar}
The graphs $K_5$, $K_{3,3}$ and $K_{2,2,2,2}$ are not planar.
\end{observation}

Let $\Gamma_1$ and $\Gamma_2$ be two finite graphs (not necessarily simple) with no parallel edges. Then the \textit{tensor product} of $\Gamma_1$ and $\Gamma_2$, denoted by $\Gamma_1\times \Gamma_2$, is a graph with the vertex set $V(\Gamma_1\times \Gamma_2)$ being the Cartesian product $V(\Gamma_1)\times V(\Gamma_2)$, and two distinct vertices $(u_1,v_1)$ and $(u_2,v_2)$ are adjacent in $\Gamma_1\times \Gamma_2$ if $u_1$ is adjacent to $u_2$ in $\Gamma_1$, and $v_1$ is adjacent to $v_2$ in $\Gamma_2$. The \textit{union} of $\Gamma_1$ and $\Gamma_2$, denoted by $\Gamma_1\cup \Gamma_2$, is simply the graph with the vertex set $V(\Gamma_1\cup \Gamma_2)=V(\Gamma_1)\cup V(\Gamma_2)$ and the edge set $E(\Gamma_1\cup \Gamma_2)=E(\Gamma_1)\cup E(\Gamma_2)$. Let $u,v\in V(\Gamma_1)\cup V(\Gamma_2)$. Then it is easy to see that $u\sim v$ in $\Gamma_1\cup \Gamma_2$ if and only if $u\sim v$ in $\Gamma_1$ or $u\sim v$ in $\Gamma_2$. Based on this simple fact, we have the following observation.
\begin{observation}\label{tensor product of graph union}
Let $\Gamma, \Gamma_1$ and $\Gamma_2$ be three graphs with no parallel edges. Then  $\Gamma\times(\Gamma_1\cup\Gamma_2)=(\Gamma\times\Gamma_1)\cup(\Gamma\times\Gamma_2).$  
\end{observation}

For two graphs $\Gamma_1$ and $\Gamma_2$ with disjoint vertex sets, the \textit{join graph} $\Gamma_1 \vee \Gamma_2$ of $\Gamma_1$ and $\Gamma_2$ is the graph obtained from the union of $\Gamma_1$ and $\Gamma_2$ by adding new edges from each vertex of $\Gamma_1$ to every vertex of $\Gamma_2$. Let $\gamma(e)$ be the graph with a single vertex $e$ and a loop on the vertex $e$. For a graph $\Gamma$ with $e\notin V(\Gamma)$, we denote by $\Gamma^e$ the join graph $\Gamma\vee \gamma(e)$. For other unexplained graph theoretic terminologies used in this paper, we refer to the book \cite{west}.

Throughout this article, we denote the cyclic group of order $n$ by $C_n$, unless stated otherwise. For $n\in\mathbb{N}$, let $[n]={1,2,\ldots,n}$. For positive integers $r$ and $k$, let $S_k(r):=1+r+\cdots+r^{k-1}$ and set $S_0(r)=0.$ We denote the identity element of any group by $e$, unless otherwise mentioned. Given a group $G$, we denote its center by $Z(G)$. For each non-identity element $a\in G,$ $G_a$ is the collection of all generators of the cyclic group $\langle a\rangle$. When $G$ is finite, we denote by $\pi(G)$ the number of distinct prime divisors of $\lvert G\rvert$. For any subset $H$ of $G$, we define the set $H^{-1}:=\{h^{-1}~\vert~h\in H\}$. For $x,y\in G$, we write $x^y=y^{-1}xy$ for the conjugate of $x$ by $y$. For a subset $S\subseteq G$ and an element $g\in G$, we define $S^g:= g^{-1}Sg$. Observe that $x\in Z(G)$ if and only if $x^y=x$ for every $y\in G$. Moreover, for all $x,y,z\in G$, $(xz)^y=x^yz^y$ and $(x^y)^z=x^{yz}$. For a subset $S\subseteq G$, the \textit{centralizer} of $S$ in $G$ is defined by $C_G(S)=\{g\in G\mid gs=sg \text{ for all } s\in S\}$. The \textit{normalizer} of $S$ in $G$, denoted by $\mathscr{N}_G(S)$ is defined by $\mathscr{N}_G(S)=\{g\in G\mid gS=Sg$. Note that for $a\in G$, we write $C_G(a)$ in place of $C_G(\{a\})$. It is well known that both $C_G(S)$ and $\mathscr{N}_G(S)$ are subgroups of $G$.

For $g,h\in G$, the \textit{commutator} of $g$ and $h$ is defined by $[g,h]=g^{-1}h^{-1}gh$. The subgroup generated by all commutators of $G$ is called the \emph{derived subgroup} (or \emph{commutator subgroup}) of $G$, and is denoted by $G'$ or $[G,G]$. Consequently, the quotient group $G/G'$ is abelian. For a prime number $p$, an \emph{elementary abelian $p$-group} is a finite abelian group with all the non identity elements having order $p$. The \emph{Frattini subgroup} of a group $G$, denoted by $\Phi(G)$, is the intersection of all maximal subgroups of $G$. If $G$ has no maximal subgroups (for example, when $G$ is trivial or a Prüfer group), then $\Phi(G)$ is defined to be $G$. $G$ is called a \textit{simple group} if $G$ is non-trivial and it has no non-trivial normal subgroup. Furthermore, $G$ is said to be \emph{nilpotent of class $n$} if its lower central series
$$
G=\gamma_1(G)\triangleright \gamma_2(G)\triangleright \cdots \triangleright \gamma_{n+1}(G)=\{e\}
$$
terminates, where
$$
\gamma_{i+1}(G)=[\gamma_i(G),G],~~ i\geq 1,
$$
and $\gamma_n(G)\neq {e}$. In particular, $\gamma_2(G)=G'=[G,G]$. A group is nilpotent of class $1$ if and only if it is abelian, while it is nilpotent of class $2$ if and only if $G'\leq Z(G)$. A group $G$ is called \emph{metacyclic} if it contains a cyclic normal subgroup $N$ such that the quotient group $G/N$ is also cyclic. 

For every positive integers $m, n, r$ with $m,n\geq 2,$ \(\gcd(r,m)=1\) and  $r^n\equiv1 \pmod m$, let
\[
G_{m,n,r}:=C_m\rtimes_r C_n
=\gen{a,b\mid a^m=b^n=1,\;b^{-1}ab=a^r}.
\] 
So, we have $D_{2m}=G_{m,2,-1}$. Also, the following groups are important in this article.
\begin{itemize}
    \item \emph{Generalized quaternion group:} $Q_{4n}:= \langle a, b ~\lvert~ a^{2n} = e, ~ b^2 = a^n, ~ b^{-1}ab = a^{-1} \rangle,~ n\geq 2.$
    \item $(C_4\times C_2)\rtimes C_2:=\langle a,b,c~\vert~a^2=b^4=c^2=[b,c]=e,~a^{-1}ba=bc,~a^{-1}ca=c\rangle$.
    \item $(C_4\times C_2)\rtimes C_4:=\langle a,b,c~\vert~a^4=b^4=c^2=[b,c]=e,~a^{-1}ba=bc,~a^{-1}ca=c\rangle$.
\end{itemize}

For a finite $2$-generated non-abelian group $G$, let $a,b\in G$ be a set of generators of $G$, such that $o(a)=o(b)=2$. Consider $c=ab$, and $o(c)=n$. Then, we have $G=\langle a,c\rangle$ with $c^n=a^2=e$, and $a^{-1}ca=c^{-1}$. Based on this simple fact, we have the following observation that is useful throughout this article.
\begin{observation}\label{Obs:generated by two elements of order 2}
Let $G$ be a finite non-abelian group generated by two elements of order $2$. Then $G$ is isomorphic to a dihedral group.   
\end{observation}

The rest of this section mostly deals with finite $p$-groups, for a prime $p$. The following results will play a crucial role in our subsequent analysis.
\begin{lemma}\cite[Chap 5, Corollary 1.2]{Gorenstein-group-book}\label{Lemma: Gorenstein_quotient by Frattini}
If $G/\Phi(G)$ is cyclic, then $G$ is cyclic.    
\end{lemma}
\begin{theorem}\cite[Chap 5, Theorem 1.3]{Gorenstein-group-book}\label{Theorem: Gorenstein_quotient by Frattini}
The Frattini factor group $G/\Phi(G)$ of a $p$-group $G$ is elementary abelian. Furthermore, $\Phi(G)=\{e\}$ if and only if $G$ is elementary abelian.  
\end{theorem}
\begin{theorem}\cite[Chap 5, Theorem 4.10(ii)]{Gorenstein-group-book}\label{Theorem: Gorenstein}
Let $G$ be a $p$-group with no non-cyclic abelian subgroup, then either $G$ is cyclic or $p=2$ and $G$ is isomorphic to $Q_{2^n}$, for $n\geq 3$.
\end{theorem}
\begin{theorem}\cite{miller_moreno}\label{Theorem:miller_moreno}
Let $G$ be a finite non-abelian non-$p$-group, such that every proper subgroup of $G$ is cyclic. Then $G\cong C_q\rtimes C_{p^m}$, where $q\equiv 1 \pmod p$. 
\end{theorem}

Based of the above results, we arrive at the following.
\begin{lemma}\label{Lemma: Cyclic abelianization}
 Let $G$ be a finite $p$-group with $G/G'$ cyclic. Then $G$ is cyclic.   
\end{lemma}
\begin{proof}
 Let $G$ be a finite $p$-group. Then from Theorem \ref{Theorem: Gorenstein_quotient by Frattini}, we have that $G/\Phi(G)$ is elementary abelian. Hence $G'\subseteq\Phi(G)$, i.e.~every maximal subgroup of $G$ contains $G'$. Now, since $G/G'$ is cyclic, $G$ has a unique maximal subgroup $M$ with $G/M$ cyclic. So $G/\Phi(G)$, and from Lemma \ref{Lemma: Gorenstein_quotient by Frattini} we have that $G$ is cyclic.
\end{proof}

\begin{lemma}\label{Lemma: one generator inside <commutator,center> implies proper subgroup}
Let $G$ be a finite $2$-generated non-abelian $p$-group, and $x\in \langle G',Z(G)\rangle$. Then $\langle x,y\rangle\subsetneq G$ for all $y\in G$.   
\end{lemma}
\begin{proof}
Let $G$ be a finite $2$-generated non-abelian $p$-group, and $x\in \langle G',Z(G)\rangle$. Assume that there exists $y\in G$ such that $\langle x,y\rangle=G$. Consider the natural projection map $\pi:G\rightarrow G/Z(G)$. Then $\pi(x)\in H'$ and $H=\langle \pi(x), \pi(y)\rangle$, where $H=G/Z(G)$. So, we have $H/H'$ is cyclic. Hence, Lemma \ref{Lemma: Cyclic abelianization} implies that $H$ is cyclic, i.e.~$G/Z(G)$ is cyclic. Now, it is easy to see that cyclic $G/Z(G)$ implies $G$ is abelian, which is contradictory. Hence, the result follows.     
\end{proof}

\begin{lemma}\label{Lemma:number of elements of order p in a p group}
Let $G$ be a finite $p$-group. Let $t$ be the number of elements of order $p$ in $G$. Then $t\equiv -1\pmod p$.   
\end{lemma}
\begin{proof}
Let $X:=\{x\in G~\vert~o(x)=p\}$. Then $\lvert X\rvert=t$. Let $X_1,X_2,\ldots,X_m$ be the distinct conjugacy classes of size at least $2$ in $X$, and $x_1,x_2,\ldots,x_m$ be elements from $X_1,X_2,\ldots,X_m$, respectively. Then it is easy to see that,  \begin{align*}
t&=\lvert Z(G)\cap X\rvert+\sum_{i=1}^m\lvert X_i\rvert\\
&=\lvert Z(G)\cap X\rvert+\sum_{i=1}^m\frac{G}{C_{G}(x_i)}.
\end{align*} 
Clearly $\frac{G}{C_{G}(x_i)}$ is divisible by $p$ for all $i\in[m]$, since $G$ is a $p$-group and $\lvert X_i\rvert\geq 2$. Also, $\lvert Z(G)\cap X\rvert$ is same as the number of elements of order $p$ is $Z(G)$. Since $Z(G)$ is abelian, the number of elements of order $p$ is $Z(G)$ is of the form $p^k-1$ for some $k\in\mathbb{N}$. So, $t\equiv -1\pmod p$.
\end{proof}

\begin{lemma}\label{Lemma: langle x^p,y_ranglesubsetneq G}
Let $G$ be a finite non-cyclic $p$-group such that $G=\langle x,y\rangle$, for some $x,y\in G$. Then $\langle x^p,y\rangle\subsetneq G$.       
\end{lemma}
\begin{proof}
$G$ is a finite $p$-group such that $G=\langle x,y\rangle$, for some $x,y\in G$. Assume $\langle x^p,y\rangle= G$. Then $x=x^{kp}y^ld$, for some $k,l\in\mathbb{N}$, and $d\in G'$. Now, take the natural projection map $f:G\rightarrow G/G'$ defined by $f(x)=\bar{x}$ and $f(y)=\bar{y}$. Then Lemma \ref{Lemma: Cyclic abelianization} implies that both $f(x)$ and $f(y)$ are nontrivial. So, $f(x)=f(x)^{kp}f(y)^l$, i.e.~$\bar{x}^{kp-1}=\bar{y}^{-l}$. Since ${\rm gcd}(kp-1,p)=1$, we have $\bar{x}\in \langle\bar{y}\rangle$, i.e.~$G/G'$ is cyclic, which is not possible by Lemma \ref{Lemma: Cyclic abelianization}. Thus, $\langle x^p,y\rangle\subsetneq G$.
\end{proof}

\section{The \(N_I\)-group}\label{Sec: The NI group}
In this section, we investigate $N_I$-groups for some important classes of $2$-generated groups. Recall that any cyclic group of order more than $2$ is not a $N_I$-group, and any group that is neither cyclic nor $2$-generated is a $N_I$-group. In the next subsection, we classify finite metacyclic groups that are $N_I$-groups.

\subsection{Metacyclic groups}
Before going into the main classification result, let us first prove some important lemmas. 
\begin{lemma}Let $t=\min\{k\ge1:S_k(r)\equiv 0\pmod m\}$ and
$d=\ord_m(r).$ Then $
S_k(r)\equiv0\pmod m$ if and only if $t\mid k,
$
and
$$
t=d\,\frac{m}{\gcd(m,S_d(r))}.
$$
So the condition $t=n$ is equivalent to
\[
\frac nd=\frac{m}{\gcd(m,S_d(r))}.
\]
\end{lemma}
\begin{proof}
		The identity
		\[
		(r-1)S_k(r)=r^k-1
		\]
		shows that \(S_k(r)\equiv0\pmod m\) implies \(d\mid k\).
		Since \(r^d\equiv1\pmod m\), for every \(u\ge1\) we have
		\[
		S_{du}(r)\equiv uS_d(r)\pmod m.
		\]
		Thus \(m\mid S_{du}(r)\) if and only if
		\(m/\gcd(m,S_d(r))\mid u\).
\end{proof}
\begin{lemma}\label{lem:arithmetic}Let $G=G_{m.n,r}$. 
With the notation above,
\[
|ab|=\operatorname{lcm}(n,t),\qquad
\gen a\cap\gen{ab}=\gen{a^{S_n(r)}},\qquad
\gen b\cap\gen{ab}=\gen{b^t}.
\]
\end{lemma}
\begin{proof}
The multiplication rule in \(G_{m,n,r}\) gives
\[
(ab)^k=a^{S_k(r)}b^k.
\]
This power is the identity precisely when both \(n\mid k\) and
\(t\mid k\), proving the formula for \(|ab|\).
A power of \(ab\) lies in \(A=\langle a\rangle\) precisely when its exponent is a
multiple of \(n\); hence
\[
\gen a\cap\gen{ab}=\gen{(ab)^n}=\gen{a^{S_n(r)}}.
\]
A power of \(ab\) lies in \(B=\langle b\rangle\) precisely when its exponent is a
multiple of \(t\); hence
\[
\gen b\cap\gen{ab}=\gen{(ab)^t}=\gen{b^t}.
\]
The corresponding intersection orders are
\(t/\gcd(t,n)\) and \(n/\gcd(t,n)\).
\end{proof}

\begin{corollary}\label{tnen}If $t\neq n,$ then $G_{m,n,r}$ is not a $N_I$-group.
\end{corollary}\begin{proof}If \(t\nmid n\), take \(x=a\) and \(y=ab\).
If \(n\nmid t\), take \(x=b\) and \(y=ab\).
\end{proof}

The converse requires a fact about split metacyclic \(p\)-groups.

\begin{lemma}\label{lem:local}
Let \(p\) be a prime and let
\[
P=\gen{c,e\mid c^{p^\alpha}=e^{p^\beta}=1,\;
ece^{-1}=c^s},
\]
where \(\alpha,\beta\ge1\), \(s\equiv1\pmod p\), and
\(s^{p^\beta}\equiv1\pmod{p^\alpha}\). Suppose that
\[
\min\{k\ge1:S_k(s)\equiv0\pmod{p^\alpha}\}=p^\beta.
\]
Then every generating pair \(u,v\) of \(P\) satisfies
\(\gen u\cap\gen v=\{1\}\).
\end{lemma}

\begin{proof}
Write \(u=c^i e^j\) and \(v=c^k e^\ell\).
The quotient obtained by imposing \(c^p=e^p=1\) is \(C_p^2\).
Since \(u,v\) generate \(P\), their coordinate vectors are linearly
independent modulo \(p\):
\[
i\ell-kj\not\equiv0\pmod p.
\]
It is enough to show that \(\gen u\) and \(\gen v\) have different
subgroups of order \(p\).

\medskip

Suppose that \(p\) is odd, or that \(p=2\) and \(s\equiv1\pmod4\).
The case \(\alpha=1\) is included by taking \(s=1\) as a
representative modulo \(p^\alpha\).
The elementary valuations of geometric sums give
\[
\min\{k\ge1:S_k(s)\equiv0\pmod{p^\alpha}\}=p^\alpha,
\]
so \(\beta=\alpha\). Expanding the geometric sum also gives
\[
S_{p^{\alpha-1}}(s^j)
 \equiv p^{\alpha-1}\pmod{p^\alpha}.
\]
Consequently,
\[
u^{p^{\alpha-1}}
 =c^{i p^{\alpha-1}}e^{j p^{\alpha-1}},
\qquad
v^{p^{\alpha-1}}
 =c^{k p^{\alpha-1}}e^{\ell p^{\alpha-1}}.
\]
These are nonidentity elements of order \(p\). Their cyclic
subgroups correspond to distinct lines in
\[
\gen{c^{p^{\alpha-1}},e^{p^{\alpha-1}}}\cong C_p^2,
\]
because the coordinate vectors are independent modulo \(p\).
Thus the cyclic subgroups of \(u\) and \(v\) intersect trivially.

\medskip

Suppose that \(p=2\), \(\alpha\ge2\), and \(s\equiv3\pmod4\).
Choose a positive representative for \(s\), and put
\[
\delta=\min\{\alpha,v_2(s+1)\}.
\]
The geometric-sum valuation yields
\[
\beta=\max\{1,\alpha-\delta+1\}.
\]

If \(\beta=1\), then \(s\equiv-1\pmod{2^\alpha}\), so \(P\) is
dihedral. A generating pair consists either of a rotation and a
reflection, or of two distinct reflections. Its cyclic subgroups
intersect trivially in either case.

Suppose now that \(\beta\ge2\). Then
\(\alpha=\beta+\delta-1\ge\beta+1\). Put
\[
U=c^{2^{\alpha-1}},\qquad V=e^{2^{\beta-1}}.
\]
For an element \(c^i e^j\) with nonzero coordinate vector modulo
\(2\), its unique involution is determined by its parity class:
\[
\begin{array}{c|ccc}
(i,j)\pmod2 &(1,0)&(0,1)&(1,1)\\ \hline
\text{unique involution}&U&V&UV.
\end{array}
\]
Indeed, for even \(j\), we have \(s^j\equiv1\pmod4\), and an
element with \(i\) odd has order \(2^\alpha\) and involution \(U\).
For odd \(j\), the valuation
\[
v_2\!\left(S_{2^h}(s^j)\right)=\delta+h-1
\qquad(h\ge1)
\]
shows that \(c^i e^j\) has order \(2^\beta\) and involution \(U^iV\).
The two generating vectors have distinct nonzero parity classes,
so their involutions are distinct. This proves the lemma.
\end{proof}

\begin{theorem}\label{metacyclic NI classification}
Let $G$ be a finite metacyclic group. Then $G$ is a $N_I$-group if and only if $G\cong G_{m,n,r}$ for a suitable choice of $(m, n, r)$ and
$
\min\{k\ge1:S_k(r)\equiv 0\pmod m\}=n.$
\end{theorem}

\begin{proof}
Assume that $G$ is a metacyclic $N_I$-group. Let $N\unlhd G$ with $N=\langle a\rangle$ and $G=N\langle b\rangle.$ Since $G=\langle a, b\rangle$, it must be $\langle a\rangle \cap \langle b\rangle=1,$ hence $G=\langle a\rangle \rtimes \langle b\rangle \cong G_{m,n,r}$ for a suitable choice of $(m,n,r)$. It remains to prove that $n=t.$

By Corollary \ref{tnen} this is a necessary condition.
For the converse, suppose that \(t=n\), and assume that a generating
pair \(x,y\) has nontrivial cyclic intersection. Every element of
\(\gen x\cap\gen y\) commutes with both generators, so
\[
\gen x\cap\gen y\le Z(G).
\]
Choose an element \(z\) of prime order \(p\) in this intersection.

Set \(g=\gcd(m,r-1)\) and \(d=\ord_m(r)\). Directly from the
presentation,
\[
Z(G)=\gen{a^{m/g}}\times\gen{b^d}.
\]
Also, \(a^{S_d(r)}\) is central because
\((r-1)S_d(r)=r^d-1\equiv0\pmod m\). Since \(t=n\),
\[
\frac nd=\frac{m}{\gcd(m,S_d(r))}
 =\left|a^{S_d(r)}\right|
\]
divides \(g\). Thus every prime divisor of \(|Z(G)|\) divides \(g\).
In particular,
\[
p\mid m,\qquad r\equiv1\pmod p.
\]
As \(S_n(r)\equiv0\pmod m\) and \(S_n(r)\equiv n\pmod p\), we also
have \(p\mid n\). Write
\[
\alpha=v_p(m)\ge1,\qquad \beta=v_p(n)\ge1.
\]

\medskip
\noindent\textbf{Case 1: \(v_p(d)=\beta\).}
The projection of the center to \(B\) is \(\gen{b^d}\), whose
order \(n/d\) is prime to \(p\). Therefore \(z\in A\).

Write \(x=a^i b^j\) and \(y=a^k b^\ell\), with nonnegative
exponents. Their images generate \(B\), so at least one of
\(j,\ell\) is prime to \(p\). Interchanging \(x,y\) if necessary,
assume \(p\nmid j\), and put
\[
L=\frac{n}{\gcd(j,n)}.
\]
Then \(L\) is the order of the image of \(x\) in \(B\), and
\[
\gen x\cap A=\gen{x^L},\qquad
x^L=a^{iS_L(r^j)}.
\]
Since \(n\mid jL\), Lemma~\ref{lem:arithmetic} and \(t=n\) give
\[
S_j(r)S_L(r^j)=S_{jL}(r)\equiv0\pmod{p^\alpha}.
\]
But \(S_j(r)\equiv j\pmod p\) is a unit modulo \(p^\alpha\).
Hence \(S_L(r^j)\equiv0\pmod{p^\alpha}\), and \(x^L\) has order
prime to \(p\). This contradicts \(z\in\gen x\cap A\).

\medskip
\noindent\textbf{Case 2: \(v_p(d)<\beta\).}
For each prime \(q\mid m\), define
\[
d_q=\ord_{q^{v_q(m)}}(r),\qquad
t_q=\min\{k\ge1:S_k(r)\equiv0\pmod{q^{v_q(m)}}\}.
\]
The Chinese remainder theorem and Lemma~\ref{lem:arithmetic} give
\[
t=\operatorname{lcm}_{q\mid m}t_q,\qquad
d=\operatorname{lcm}_{q\mid m}d_q,\qquad
\frac{t_q}{d_q}\text{ is a power of }q.
\]
For \(q\ne p\), it follows that
\[
v_p(t_q)=v_p(d_q)\le v_p(d)<\beta.
\]
Since \(v_p(t)=v_p(n)=\beta\), we must have \(v_p(t_p)=\beta\).
Moreover, \(r\equiv1\pmod p\) implies that both \(d_p\) and \(t_p\)
are powers of \(p\). Thus
\[
t_p=p^\beta.
\]

The order of \(r\) modulo \(p^\alpha\) is a power of \(p\) dividing
\(n\), hence it divides \(p^\beta\). Therefore the presentation
defines a surjection
\[
G\longrightarrow
P=\gen{c,e\mid c^{p^\alpha}=e^{p^\beta}=1,\;ece^{-1}=c^r},
\qquad a\longmapsto c,\quad b\longmapsto e.
\]
Its kernel has order \(mn/p^{\alpha+\beta}\), which is prime to
\(p\). The image of \(z\) is consequently nonidentity and belongs
to the cyclic subgroups of both images of \(x,y\).
Those images generate \(P\), contradicting
Lemma~\ref{lem:local}, because \(t_p=p^\beta\).
This completes the proof.
\end{proof}

\begin{corollary}\label{Cor: 2-generated p group}
The direct product \(C_m\times C_n\), with \(m,n\ge2\) is a $N_i$-group if and only
if \(m= n\).
\end{corollary}
\begin{proof}
For \(r=1\), we have \(S_k(1)=k\), so \(t=m\).
\end{proof}

\begin{corollary}\label{Cor: D_2n NI group}
The dihedral group $D_{2m}$, with \(m\ge2\), is a $N_I$-group.
\end{corollary}
\begin{proof}
We have \(S_1(-1)=1\) and \(S_2(-1)=0\), so \(t=2=n\).
\end{proof}

\begin{corollary}\label{faith}
If the action of \(B=\langle b\rangle\) on \(A=\langle a\rangle\) is faithful, then $G_{m,n,r}$ is a $N_I$-group if and only if
\[
S_n(r)\equiv0\pmod m.
\]
\end{corollary}
\begin{proof}
Faithfulness means \(d=n\), so \(n\mid t\). Therefore, \(t\ne n\)
is equivalent to \(t\nmid n\), which in turn is equivalent to
the displayed condition.
\end{proof}

The following examples illustrate the dependence on the action:
\[
\begin{array}{ccl}
\toprule
G & t & \text{Conclusion}\\
\midrule
C_8\rtimes_3 C_2 &4&
  \text{No: }(a,ab),\quad (ab)^2=a^4\ne1.\\
C_8\rtimes_7 C_2 &2&
  \text{Yes.}\\
C_4\rtimes_3 C_4 &2&
  \text{No: }(b,ab),\quad (ab)^2=b^2\ne1.\\
C_9\rtimes_4 C_9 &9&
  \text{Yes.}\\
\bottomrule
\end{array}
\]

\begin{corollary}\label{Cor: metacyclic p index not NI group}
	Suppose that a finite non-abelian group $G$ contains a normal cyclic subgroup of prime index $p$. Then $G$ is a $N_I$-group if and only if one of the following occur:
	\begin{enumerate}
		\item $p=2$ and $G\cong G_{m,2,-1}$ is a dihedral group.
		\item $p$ is odd, $G\cong G_{p^\alpha  u,p,r}$,  $\alpha\leq 1$ and $\gcd(u,p)=1.$ 
	\end{enumerate}
	\begin{proof}By Theorem \ref{metacyclic NI classification}, $G\cong G_{m,p,r}$ for suitable integers $m$ and $r$.	Since $G$ is non-abelian, the action of $B$ on $A$ is faithful, so, by Corollary \ref{faith}, $G$ is a $N_I$-groups if and only if $S_p(r)\equiv 0 \pmod m.$ If $p=2$, this is equivalent to say that $r=-1 \pmod m.$
		If $p$ is odd the conclusion follows from the fact that $$\frac{m}{\gcd(m,S_p(r))}=p^{\max\{\alpha-1,0\}}\gcd(u,r-1).$$ 
	\end{proof}
\end{corollary}

\subsection{Direct product}
Our next result deals with direct product of $N_I$-groups.
\begin{lemma}\label{2-generated NI}
Let $G_1$ and $G_2$ be two finite groups with ${\rm gcd}(\lvert G_1\rvert,\lvert G_2\rvert)=1$. Then the following are equivalent:
\begin{enumerate}
    \item Both $G_1$ and $G_2$ are $2$-generated $N_I$-groups,
    \item $G_1\times G_2$ is a $2$-generated $N_I$-group.
\end{enumerate}
\end{lemma}
\begin{proof}
Let $G_1=\langle a_1,a_2\rangle$ and $G_2=\langle b_1,b_2\rangle$. Then it can be seen that $\langle(a_1,b_1),(a_2,b_2)\rangle=G_1\times G_2$, since ${\rm gcd}(\lvert G_1\rvert,\lvert G_2\rvert)=1$. So, $G_1\times G_2$ is $2$-generated. Now, let $(c_1,d_1),(c_2,d_2)\in G_1\times G_2$ be two nontrivial elements of $G_1\times G_2$ such that $(c_1,d_1)\nsim(c_2,d_2)$ in ${\rm NGen}(G_1\times G_2)$. Then $\langle(c_1,d_1),(c_2,d_2)\rangle=G_1\times G_2$, hence $\langle c_1,c_2\rangle=G_1$ and $\langle d_1,d_2\rangle=G_2$. Since $G_1$ and $G_2$ are $N_I$-groups, we have $\langle c_1\rangle\cap\langle c_2\rangle=\{e\}=\langle d_1\rangle\cap\langle d_2\rangle$. So, $\langle (c_1,d_2)\rangle\cap\langle (c_2,d_2)\rangle=\{e\}$, i.e.~$(c_1,d_1)\nsim(c_2,d_2)$ in $\mathcal{G}_I(G_1\times G_2)$. Thus $G_1\times G_2$ is a $N_I$-group.

Conversely, let $G_1\times G_2=\langle(a_1,b_1),(a_2,b_2)\rangle$. Then clearly $G_1=\langle a_1,a_2\rangle$ and $G_2=\langle b_1,b_2\rangle$. So, both $G_1$ and $G_2$ are either cyclic or $2$-generated. Now, let $c_1,c_2\in G_1$ be two nontrivial elements of $G_1$ such that $c_1\nsim c_2$ in ${\rm NGen}(G_1)$. Then $\langle c_1,c_2\rangle=G_1$. Now, by an argument similar to the one used previously, we have $\langle (c_1,b_1),(c_2,b_2)\rangle=G_1\times G_2$. Since $G_1\times G_2$ is a $N_I$-group, we have $\langle (c_1,b_1)\rangle\cap\langle(c_2,b_2)\rangle=\{e\}$. Also, we have ${\rm gcd}(\lvert G_1\rvert,\lvert G_2\rvert)=1$. Hence, $\langle c_1\rangle\cap\langle c_2\rangle=\{e\}$, i.e.~$c_1\nsim c_2$ in $\mathcal{G}_I(G_1)$. So, $G_1$ is a $N_I$-group. Similarly, it can also be shown that $G_2$ is a $N_I$-group. As a consequence, both $G_1$ and $G_2$ are non-cyclic. Thus the result holds.    
\end{proof}

As a corollary of the above, we have the following.
\begin{corollary}\label{Cor:2-generated NI}
Let $G$ be a finite nilpotent group. Then $G$ is a $2$-generated $N_I$-group if and only of all its Sylow subgroups are $2$-generated $N_I$-groups.   
\end{corollary}

It is easy to see from the arguments used in the proof of Lemma \ref{2-generated NI} that for two finite $N_I$-groups $G_1$ and $G_2$ of relatively prime order, $G_1\times G_2$ is either generated by at least $3$ distinct elements, or $G_1\times G_2$ is a $2$-generated $N_I$-group. In either case, we have the following.

\begin{observation}\label{Obs: Direct product NI}
Let $G_1$ and $G_2$ be two finite $N_I$-groups of relatively prime order. Then $G_1\times G_2$ is a $N_I$-group.   
\end{observation}

Throughout the remainder of this article, our attention is primarily devoted to the difference graphs of finite $N_I$-groups. Henceforth, unless explicitly stated otherwise, all our groups are finite $N_I$-groups.

\section{Completeness and nullness}\label{Sec: Completeness and nullness}

In this section, we classify finite groups whose difference graphs are either a complete graph, or a null graph. We begin with the following lemma, which plays a crucial role in the proof of the main classification result.
\begin{lemma}\label{Lemma:Semidirect product null}
Let $H\cong C_q\rtimes_r C_{p^m}$, where $p,q$ are distinct primes, with $q\equiv 1~(\pmod p)$, and $r^p\equiv 1\pmod q$. Then $\mathcal{DNI}(H)$ is the null graph if and only if $m=1$.    
\end{lemma}
\begin{proof}
Let $H=\langle a,b \rangle $, with $a^q=b^{p^m}=1\text{ and }b^{-1}ab=a^r$. If $m\geq 2$, then $b^p\neq e$ and $b^{-p}ab^p=a^{r^p}=a$. Hence, $b^p\in Z(H)$. Since $H$ is non-abelian, we have $\langle a,b^p\rangle\subsetneq H$, and clearly $\langle a\rangle\cap\langle b^p\rangle=\{e\}$. So, $a\sim b^p$ in $\mathcal{DNI}(H)$ and $\mathcal{DNI}(H)$ is not null.

Now, consider $m=1$. Then $\lvert H\rvert=pq$, and every non-trivial proper subgroup of $H$ is cyclic of order $p$ or $q$. Let $x,y\in H$ be two distinct elements. Then the orders of $x$ and $y$ are either $p$ or $q$, since $H$ is not cyclic. If $o(x)\neq o(y)$, then we have $\langle x\rangle\cap\langle y\rangle=\{e\}$, i.e.~$x\nsim y$ in $\mathcal{G}_I(H)$. If $o(x)=o(y)$, then either $\langle x\rangle=\langle y\rangle$ or $\langle x\rangle\cap\langle y\rangle=\{e\}$. So, we have that $x\sim y$ in $\mathcal{G}_I(H)$ if and only if $\langle x\rangle=\langle y\rangle$. Thus, $H$ is a $N_I$-group. Now, consider two distinct elements $u,v\in H\setminus\{e\}$ such that $u\nsim v$ in $\mathcal{G}_I(H)$, i.e.~$\langle u\rangle\neq\langle v\rangle$. Then $\lvert\langle u,v\rangle\rvert>p,q$, i.e.~$\langle u,v\rangle=H$. Hence, we have $\bar{\mathcal{G}_I}(H)={\rm Gen}(H)$, i.e.~$\mathcal{DNI}(H)$ is the null graph.
\end{proof}

\begin{theorem}\label{Theorem: Completeness and nullness}
Let $G$ be a finite $N_I$-group. Then we have the following:
\begin{enumerate}
    \item $\mathcal{DNI}(G)$ is complete if and only if $G\cong C_2^k$, for some $k\in\mathbb{N}$ with $k\geq 3$,
    \item $\mathcal{DNI}(G)$ is the null graph if and only if either $G\cong C_p\times C_p$, or $G\cong C_q\rtimes_r C_p$, where $p,q$ are distinct primes with $q\equiv 1 \pmod p$.
\end{enumerate}
\end{theorem}
\begin{proof}
Let $G$ be a finite $N_I$-group.
\begin{enumerate}
    \item Let $\mathcal{DNI}(G)$ be complete. Then $\mathcal{G}_I(G)\vert_{G\setminus\{e\}}$ is the null graph. Suppose $G$ has an element $g$, with $o(g)\geq 3$. Then $g\neq g^{-1}$ and $g\sim g^{-1}$ in $\mathcal{G}_I(G)$, which contradicts our assumption. So, every non-trivial element of $G$ has order $2$. Hence, $G\cong C_2^k$ for some $k\in\mathbb{N}$. Since $G$ is a $N_I$-group, we have $k\geq 2$. Also, if $k=2$, then $G$ is two generated and any two non-trivial element of $G$ generates $G$, thus $\mathcal{DNI}(G)$ is null graph. So $G\cong C_2^k$, with $k\geq 3$.
    
    Conversely, let $G\cong C_2^k$ for some $k\in\mathbb{N}$, with $k\geq 3$. Then it is easy to see that $G$ is neither cyclic, nor $2$-generated, and every non-trivial element of $G$ is of order $2$. Thus, any two distinct nontrivial elements of $G$ are adjacent in ${\rm NGen}(G)$, but not adjacent in $\mathcal{G}_I(G)$. So $\mathcal{DNI}(G)$ is complete. 

    \item We consider the following two cases. 
    \begin{itemize}
        \item[\underline{{\bf $G$ is abelian:}}] If $G$ is an abelian non-$p$-group, then there exists $x,y\in G$ with ${\rm gcd}(o(x),o(y))=1$, i.e.~$x\nsim y$ in $\mathcal{G}_I(G)$. Also, $\langle x,y\rangle=\langle xy\rangle\neq G$, since $G$ is not cyclic. Thus, $x\sim y$ in $\mathcal{DNI}(G)$. 
        
        If $G$ is an abelian $p$-group, then either $G\cong C_p\times C_p$ or $G$ contains a proper subgroup $H\cong C_p\times C_p$. Without loss of generality, assume $H=G$ for $G\cong C_p\times C_p$. Consider two distinct non-trivial elements $x,y\in H$, such that $\langle x\rangle\cap\langle y\rangle=\{e\}$. So, when $H=G$ we have $\langle x,y\rangle=G$, and when $H\subsetneq G$ we have $\langle x,y\rangle\subsetneq G$. Hence, $x\sim y$ in $\mathcal{DNI}(G)$ if and only if $H\subsetneq G$. Since $x,y$ is arbitrary, we have that $\mathcal{DNI}(G)$ is null if and only if $G\cong C_p\times C_p$. 

        \item[\underline{{\bf $G$ is non-abelian:}}] If $G$ is a $p$-group, then Theorem \ref{Theorem: Gorenstein} implies that $G$ has a proper subgroup isomorphic to $C_p\times C_p$, and a set of arguments similar to the previous case imply that $\mathcal{DNI}(G)$ is not null.
        
        If $G$ is a nilpotent non-$p$-group, there exists $x,y\in G$, such that ${\rm gcd}(o(x),o(y))=1$ with $[x,y]=e$. Then $\langle x,y\rangle=\langle xy\rangle\neq G$, since $G$ is not cyclic. So $x\sim y$ in $\mathcal{DNI}(G)$, i.e.~$\mathcal{DNI}(G)$ is not null.  

        If $G\cong D_{2n}$ for some $n\in\mathbb{N}$ with $n\geq 3$, consider $x,y\in G$ such that $o(x)=n, o(y)=2$ and $y\notin Z(G)$. If $n$ is composite, let $m\geq 2$ be a proper divisor of $n$. Then $\langle x^m, y\rangle\subsetneq G$ with $\langle x^m\rangle\cap\langle y\rangle=\{e\}$, i.e.~$x^m\sim y$ in $\mathcal{DNI}(G)$. If $n$ is prime, then $D_{2n}\cong C_n\rtimes C_2$, and from Lemma \ref{Lemma:Semidirect product null} we have that $\mathcal{DNI}(C_n\rtimes C_2)$ is a null graph.

        Finally, let $G$ be neither nilpotent, nor dihedral. Clearly, $\pi(G)\geq 2$.  If $\lvert G\rvert$ is even, then there exists $x,y\in G$ such that $o(x)=2$ and $o(y)=p$, for some odd prime $p$. If $\langle x,y\rangle\neq G$, then $x\sim y$ in $\mathcal{DNI}(G)$, i.e.~$\mathcal{DNI}(G)$ is not null. If $\langle x,y\rangle=G$, let $x'=y^{-1}xy$. Then $o(x')=o(x)=2$ with $x\neq x'$, and since $G$ is neither abelian, nor dihedral, from Observation \ref{Obs:generated by two elements of order 2} we conclude that $\langle x,x'\rangle\subsetneq G$. Hence, $x\sim x'$ in $\mathcal{DNI}(G)$, i.e.~$\mathcal{DNI}(G)$ is not null. Now, assume $\lvert G\rvert$ is odd. Suppose $G$ has a non-cyclic proper subgroup $B$. If $\pi(B)\geq 2$, then there exists $x,y\in B$, such that ${\rm gcd}(o(x),o(y))=1$, hence $\langle x\rangle\cap\langle y\rangle=\{e\}$. If $B$ is a $p$-group, then in light of Theorem \ref{Theorem: Gorenstein} we can say that $B$ has a subgroup isomorphic to $C_p\times C_p$. So, some previous arguments yield that there exists $x,y\in B$, such that $\langle x\rangle\cap\langle y\rangle=\{e\}$. Hence, in both the cases we have $x\sim y$ in $\mathcal{DNI}(G)$, i.e.~$\mathcal{DNI}(G)$ is not null. Thus, we consider $G$ to be a non-abelian non-$p$-group with no non-cyclic proper subgroup. Then, Theorem \ref{Theorem:miller_moreno} implies that $G\cong C_q\rtimes_r C_{p^m}$, for distinct odd primes $p,q$, with $q\equiv 1~(\pmod p)$. Again, Lemma \ref{Lemma:Semidirect product null} says that $\mathcal{DNI}(C_q\rtimes_r C_{p^m})$ is a null graph if and only if $m=1$. Hence, $G\cong C_q\rtimes_r C_p$. 
    \end{itemize}
\end{enumerate}
This concludes the proof.
\end{proof}

\section{Direct product of groups}\label{Sec: Direct product of groups}

This section is devoted to the study of the difference graph of the direct product of finite groups of relatively prime orders. Let $\bar{\mathcal{G}_I}(G)$ denote the complement of the induced subgraph $\mathcal{G}_I\vert_{G\setminus\{e\}}$, for a finite group $G$. The following results describe the join graph $\mathcal{DNI}^e$ associated with the direct product of finite groups of relatively prime orders in terms of unions of tensor products of suitable graphs.
\begin{theorem}\label{Thm: direct product}
Let $G_1$ and $G_2$ be two groups with ${\rm gcd}(\lvert G_1\rvert,\lvert G_2\rvert)=1$. Then $\mathcal{DNI}(G_1\times G_2)^e=(\bar{\mathcal{G}_I}(G_1)^e\times\mathcal{DNI}(G_2)^e)\cup(\mathcal{DNI}(G_1)^e\times\bar{\mathcal{G}_I}(G_2)^e)$.   
\end{theorem}
\begin{proof}
Let $G_1$ and $G_2$ be two finite $N_I$-groups with ${\rm gcd}(\lvert G_1\rvert,\lvert G_2\rvert)=1$. Consider the graphs $\mathcal{DNI}(G_1\times G_2)^e$ and $(\bar{\mathcal{G}_I}(G_1)^e\times\mathcal{DNI}(G_2)^e)\cup(\mathcal{DNI}(G_1)^e\times\bar{\mathcal{G}_I}(G_2)^e)$. Then it is easy to see that $V(\mathcal{DNI}(G_1\times G_2)^e)=V((\bar{\mathcal{G}_I}(G_1)^e\times\mathcal{DNI}(G_2)^e)\cup(\mathcal{DNI}(G_1)^e\times\bar{\mathcal{G}_I}(G_2)^e))=G_1\times G_2$. Further, it can be seen that in both the graphs $\mathcal{DNI}(G_1\times G_2)^e$ and $(\bar{\mathcal{G}_I}(G_1)^e\times\mathcal{DNI}(G_2)^e)\cup(\mathcal{DNI}(G_1)^e\times\bar{\mathcal{G}_I}(G_2)^e)$, the identity element $e$ is a dominant vertex with a loop. Let $(a_1,b_1)$ and $(a_2,b_2)$ be two nontrivial elements of $G_1\times G_2$. If $(a_1,b_1)\sim(a_2,b_2)$ in $\mathcal{DNI}(G_1\times G_2)^e$, we have $\langle (a_1,b_1),(a_2,b_2)\rangle\neq G_1 \times G_2$ and $\langle(a_1,b_1)\rangle\cap\langle(a_2,b_2)\rangle=\{e\}$. Since ${\rm gcd}(\lvert G_1\rvert,\lvert G_2\rvert)=1$, we have $\langle a_1 \rangle\cap\langle a_2\rangle=\{e\}=\langle b_1 \rangle\cap\langle b_2\rangle$, and $\langle(a_1,a_2)\rangle\neq G_1$ or $\langle(b_1,b_2)\rangle\neq G_2$. Hence $(a_1,b_1)\sim(a_2,b_2)$ in $\bar{\mathcal{G}_I}(G_1)^e\times\mathcal{DNI}(G_2)^e$ or $(a_1,b_1)\sim(a_2,b_2)$ in $\mathcal{DNI}(G_1)^e\times\bar{\mathcal{G}_I}(G_2)^e$. 

Conversely, let $(a_1,b_1)\sim(a_2,b_2)$ in $(\bar{\mathcal{G}_I}(G_1)^e\times\mathcal{DNI}(G_2)^e)\cup(\mathcal{DNI}(G_1)^e\times\bar{\mathcal{G}_I}(G_2)^e)$. Then $(a_1,b_1)\sim(a_2,b_2)$ in $\bar{\mathcal{G}_I}(G_1)^e\times\mathcal{DNI}(G_2)^e$ or $(a_1,b_1)\sim(a_2,b_2)$ in $\mathcal{DNI}(G_1)^e\times\bar{\mathcal{G}_I}(G_2)^e$. If $(a_1,b_1)\sim(a_2,b_2)$ in $\bar{\mathcal{G}_I}(G_1)^e\times\mathcal{DNI}(G_2)^e$, we have $\langle a_1\rangle\cap \langle a_2\rangle=\{e\}=\langle b_1\rangle\cap \langle b_2\rangle$ and $\langle b_1,b_2\rangle\neq G_2$. So, $\langle (a_1,b_1)\rangle\cap \langle (a_2,b_2)\rangle=\{e\}$ and $\langle (a_1,b_1),(a_2,b_2)\rangle\neq G_1\times G_2$. Hence $(a_1,b_1)\sim(a_2,b_2)$ in $\mathcal{DNI}(G_1\times G_2)^e$. Similarly, it can be shown that if $(a_1,b_1)\sim(a_2,b_2)$ in $\mathcal{DNI}(G_1)^e\times\bar{\mathcal{G}_I}(G_2)^e$ then $(a_1,b_1)\sim(a_2,b_2)$ in $\mathcal{DNI}(G_1\times G_2)^e$. This concludes the proof.
\end{proof}

It is easy to see that if a group $G$ is at least $3$-generated, then $G$ is a $N_I$-group, and $\bar{\mathcal{G}_I}(G)=\mathcal{DNI}(G)$. So, as a simple corollary of Theorem \ref{Thm: direct product}, we have the following.

\begin{corollary}\label{Cor: direct product non 2-generated}
Let $G_1$ and $G_2$ be two non-cyclic non-$2$-generated groups such that ${\rm gcd}(\lvert G_1\rvert,\lvert G_2\rvert)=1$. Then $\mathcal{DNI}(G_1\times G_2)^e=\mathcal{DNI}(G_1)^e\times\mathcal{DNI}(G_2)^e$.    
\end{corollary}

Again, the following result is easy to observe.
\begin{observation}\label{Intersection power graph compliment observation}
For two non-cyclic groups $G_1$ and $G_2$ with ${\rm gcd}(\lvert G_1\rvert,\lvert G_2\rvert)=1$, we have $\bar{\mathcal{G}_I}(G_1\times G_2)^e=\bar{\mathcal{G}_I}(G_1)^e\times\bar{\mathcal{G}_I}(G_2)^e$.    
\end{observation}

As a consequence of the above results, we now arrive at the following corollary.
\begin{corollary}\label{direct product nilpotent}
Let $H$ be a finite $2$-generated nilpotent group with $H_1, H_2,\dots ,H_k$ be the Sylow subgroups of $H$. Then $\mathcal{DNI}(H)^e=\underset{i\in[k]}{\cup}(\bar{\mathcal{G}_I}(H_1)^e\times\cdots\times\mathcal{DNI}(H_i)^e\times\cdots\times\bar{\mathcal{G}_I}(H_k)^e)$. 
\end{corollary}
\begin{proof}
Let $H$ be a $2$-generated nilpotent group with $H_1, H_2,\dots ,H_k$ being the Sylow subgroups of $H$. Since $H$ is a $2$-generated $N_I$-group, corollary \ref{Cor:2-generated NI} implies that $H_i$ is also a $2$-generated $N_I$-group, for all $i\in [k]$. Now, from Theorem \ref{Thm: direct product} we have $\mathcal{DNI}(H)^e=(\bar{\mathcal{G}_I}(\prod_{i\in[k-1]}H_i)^e\times\mathcal{DNI}(H_k)^e)\cup(\mathcal{DNI}(\prod_{i\in[k-1]}H_i)^e\times\bar{\mathcal{G}_I}(H_k)^e)$. Then, from Observations \ref{tensor product of graph union} and \ref{Intersection power graph compliment observation} the rest follows. 
\end{proof}

\section{Universality}\label{Sec: Universality}
A class of graphs is said to be a \textit{universal family}, if any finite graph is isomorphic to a induced subgraph of some graph from that family. Our next result shows that the class of difference graphs is universal.

\begin{theorem}\label{Thm:DNI are universal}
Let $\Gamma$ be a finite, simple and undirected graph. Then there exists a finite abelian $N_I$-group $A$, such that $\Gamma$ is an induced subgraph of $\mathcal{DNI}(A)$.    
\end{theorem}
\begin{proof}
The proof proceeds via induction. Let $\Gamma_1$ be the trivial graph. Clearly, for any finite non-cyclic non-$2$-generated abelian group $A$, $\mathcal{DNI}(A)$ contains $\Gamma_1$ as an induced subgraph. Assume that the hypothesis is true for all graphs with at most $n-1$ vertices, for some positive integer $n\geq 2$. Let $\Gamma_{n}$ be a finite, simple and undirected graph with $n$ vertices. Let us denote the vertices of $\Gamma_{n}$ by $v_1,v_2,\dots,v_n$. Denote the induced subgraph $\Gamma_n\vert_{\{v_1,\dots,v_{n-1}\}}$ by $\Gamma_{n-1}$. Then $\Gamma_{n-1}$ is a finite, simple and undirected graph with $n-1$ vertices. So, by the inductive hypothesis, there exists a finite non-cyclic non-$2$-generated abelian group $H$, such that $\mathcal{DNI}(H)$ contains $\Gamma_{n-1}$ as an induced subgraph. Without loss of generality, we denote the vertices of that induced subgraph of $\mathcal{DNI}(H)$ by $v_1,v_2,\dots,v_{n-1}$ accordingly.   

Let $d_{\Gamma_n}(v_n)=n-1$. Consider a finite non-cyclic non-$2$-generated abelian group $K$, such that ${\rm gcd}(\lvert H\rvert,\lvert K\rvert))=1$. Clearly $H\times K$ is a finite non-cyclic non-$2$-generated abelian group. Let $k\in K$ be a nontrivial element. Then Corollary \ref{Cor: direct product non 2-generated} implies that $(v_i,e)\sim(e,k)$ in $\mathcal{DNI}(H\times K)$, for all $i\in[n-1]$, and $(v_i,e)\sim (v_j,e)$ in $\mathcal{DNI}(H\times K)$ if and only if $v_i\sim v_j$ in $\mathcal{DNI}(H)$, for all $i,j\in [n-1]$. This gives an embedding of $\Gamma_n$ in $\mathcal{DNI}(H\times K)$. Now, let $d_{\Gamma_n}(v_n)<n-1$. Without loss of generality, assume that $V(\Gamma_n)\setminus N_{\Gamma_n}[v_n]=\{v_1,v_2,\dots,v_m\}$, where $1\leq m\leq n-1$. Let $p_1,p_2,\ldots,p_m$ be $m$ distinct primes such that ${\rm gcd}(\lvert H\rvert,p_i)=1$, for all $i\in[m]$. Consider the group $K=\prod_{i\in [m]}C_{p_i}^3$. Then both $K$ and $H\times K$ are finite non-cyclic non-$2$-generated abelian groups. Let $x_i$ be an element of order $p_i$ in $K$, for $i\in [m]$. Consider the subset $S\subseteq (H\times K\setminus\{e\})$ defined as $S:=\{(v_1,x_1),(v_2,x_2),\ldots,(v_m,x_m),(v_{m+1},e),\ldots,(v_{n-1},e),(e,x_1x_2\cdots x_m)\}$. Clearly, $x_i\sim x_j$ and $x_i\nsim x_1x_2\cdots x_m$ in $\mathcal{DNI}(K)$, for $i,j\in[m]$ with $i\neq j$. Then Corollary \ref{Cor: direct product non 2-generated} implies that $(v_i,e)\sim (e,x_1x_2\cdots x_m)$ and $(v_j,x_j)\nsim (e,x_1x_2\cdots x_m)$ in $\mathcal{DNI}(H\times K)$, for all $i\in \{m+1,\ldots,n-1\}$, and $j\in [m]$. In addition, it can be seen that the adjacency among the vertices $(v_1,x_1),(v_2,x_2),\ldots,(v_m,x_m),(v_{m+1},e),\ldots,(v_{n-1},e)$ in $\mathcal{DNI}(H\times K)\vert_{S\setminus \{(e,x_1x_2\cdots x_m)\}}$ is completely determined by the adjacency among $v_1,v_2,\dots,v_{n-1}$ in $\mathcal{DNI}(H)\vert_{\{v_1,v_2,\dots,v_{n-1}\}}$. Thus, $\mathcal{DNI}(H\times K)\vert_S$ is an embedding of $\Gamma_n$ in $\mathcal{DNI}(H\times K)$. This concludes the proof.
\end{proof}

\section{Connectedness and diameter}\label{Sec: Connectedness and diameter}
In this section, we study the connectedness and diameter of the difference graph. We first study groups that are not $2$-generated and subsequently examine $2$-generated nilpotent groups. Recall that for each non-identity element $a\in G,$ $G_a$ is the collection of all generators of the cyclic group $\langle a\rangle$. The following results study adjacency among the vertices of the intersection power graph $\mathcal{G}_I(G)$, which will be useful later in this section.  
\begin{lemma}\label{Lemma: Ga has no edge in DNI(G)}
Let $G$ be a finite group. Then for each non-identity element $a\in G$,  $\mathcal{G}_I(G)\vert_{G_a}$ is a clique.
\end{lemma}
\begin{proof}
Let $G$ be a finite group, and $a\in G$ be a non-trivial element. Let $x, y$ be two arbitrary distinct elements in $G_a$. Since $\langle x\rangle=\langle y \rangle=\langle a \rangle$ and $a\neq e$, we have $x\sim y$ in $\mathcal{G}_I(G)$. Thus, $\mathcal{G}_I(G)\vert_{G_a}$ is a clique.
\end{proof}

\begin{lemma}\label{neighborhood intersection power graph}
Let $G$ be a finite $p$-group, and $x,y$ be two nontrivial elements of $G$. If $x\sim y$ in $\mathcal{G}_I(G)$, then $N_{\mathcal{G}_I(G)}[x]=N_{\mathcal{G}_I(G)}[y]$.    
\end{lemma}
\begin{proof}
Let $G$ be a finite $p$-group, and $x,y$ be two nontrivial elements of $G$, such that $x\sim y$ in $\mathcal{G}_I(G)$. Define $H=\langle x\rangle\cap \langle y\rangle$. Then $H\neq \{e\}$. Let $K$ be the unique subgroup of order $p$ in $H$. Then $K$ is the unique subgroup of order $p$ in both $\langle x\rangle$ and $\langle y\rangle$. Let $z\in N_{\mathcal{G}_I(G)}(x)\setminus\{y\}$ be an arbitrary nontrivial element. Then $\langle x\rangle\cap \langle z\rangle\neq \{e\}$. Since $K$ is the unique subgroup of order $p$ in $\langle x\rangle$, we have $K\subseteq\langle x\rangle\cap \langle z\rangle $, i.e.~$K\subseteq\langle z\rangle$. So $\langle y\rangle\cap \langle z\rangle\neq \{e\}$, i.e.~$z\in N_{\mathcal{G}_I(G)}(y)\setminus\{x\}$. Similarly, it can be shown that $N_{\mathcal{G}_I(G)}(y)\setminus\{x\}\subseteq N_{\mathcal{G}_I(G)}(x)\setminus\{y\}$. Now, since $e$ is inside both $N_{\mathcal{G}_I(G)}(x)\setminus\{y\}$ and $N_{\mathcal{G}_I(G)}(y)\setminus\{x\}$, we have $N_{\mathcal{G}_I(G)}[x]=N_{\mathcal{G}_I(G)}[y]$.  
\end{proof}

\begin{proposition}\label{Prop: IPG(G) is dominatable iff G is cyclic; G is any grp}
Let $G$ be a finite group. Then the intersection power graph $\mathcal{G}_I(G)$ has a dominant vertex distinct from the identity if and only if, for every prime divisor of $|G|$, the group $G$ possesses a unique subgroup of that prime order.
\end{proposition}
\begin{proof}
Let $p_1, \cdots, p_r$ be prime divisors of $|G|.$ Suppose the number of distinct cyclic subgroups of order $p_i$ in $G$ is not unique, for some $i\in [r].$ Let $H_1=\langle a_1\rangle, \cdots, H_k=\langle a_k\rangle, k\geq 3$ be the collection of all distinct cyclic subgroups of order $p_i$ of $G.$ Suppose $g\in G$ be a dominating vertex of $\mathcal{G}_I(G).$ Then $g\sim a_i$ in $\mathcal{G}_I(G)$ implies that $a_i\in\langle g\rangle,$ for all $i\in [k].$ This is not possible. Conversely, suppose that for every prime divisor $p$ of $|G|$, the group $G$ has a unique subgroup of order $p$. Let
$H_i=\langle a_i\rangle$
denote the unique subgroup of order $p_i$, where $p_1,p_2,\ldots,p_r$ are the distinct prime divisors of $|G|$. Since each $H_i$ is unique, it is normal in $G$. Moreover,
$H_i\cap H_j=\{e\}$
for all distinct $i,j\in [r]$. It follows that
$H_1H_2\cdots H_r \cong H_1\times H_2\times\cdots\times H_r,$
which is cyclic because the orders of the subgroups are pairwise coprime. Hence,
$H_1H_2\cdots H_r=\langle a_1a_2\cdots a_r\rangle.$
Since every element of prime order in $G$ belongs to one of the subgroups $H_i$, every nontrivial cyclic subgroup of $G$ intersects $\langle a_1a_2\cdots a_r\rangle$ nontrivially. Therefore, the vertex $a_1a_2\cdots a_r$ is adjacent to every other vertex of $\mathcal{G}_I(G)$, and hence it is a dominating vertex.
\end{proof}

\begin{lemma}\label{Lemma: x in Ga and y in Gb, if x sim y then any two elmnts in Ga and Gb edge connd in DNI(G)}
Let $a,b \in G$ be two distinct non-trivial elements of $G$, such that $G_a\neq G_b$. If an element of $G_a$ is adjacent to an element of $G_b$, then each element of $G_a$ is adjacent to every elements of $G_b$ in $\mathcal{DNI}(G)$.
\end{lemma}
\begin{proof}
Let $x \in G_a$ and $y \in G_b$ such that $x \sim y$ in $\mathcal{DNI}(G)$. Then $\langle x, y\rangle\neq G$ and $\langle x \rangle \cap \langle y \rangle=\{e\}$. Let $x_1 \in G_a$ and $y_1 \in G_b$ be two arbitrary elements. So $\langle x_1, y_1\rangle\neq G$, since $\langle x_1, y_1\rangle\subseteq\langle x, y\rangle\neq G$. In addition, $\langle x \rangle \cap \langle y \rangle=\{e\}$ implies $\langle x_1 \rangle \cap \langle y_1 \rangle=\{e\}$. Thus, we have $x_1\sim y_1$ in $\mathcal{DNI}(G)$. This concludes the proof.
\end{proof}

\subsection{Non 2-generated groups}
Throughout this subsection, we assume that $G$ is a finite non-$2$-generated group. Since every finite non-$2$-generated group is an $N_I$-group, it follows that
$\mathcal{DNI}(G)=\bar{\mathcal{G}_I}(G)\vert_{G\setminus\{e\}}$
where $\bar{\mathcal{G}_I}(G)\vert_{G\setminus\{e\}}$ denotes the subgraph of $\bar{\mathcal{G}_I}(G)$ induced by the vertex set $G\setminus{e}$. We denote the number of $p$-ordered subgroups of $G$ by $n_p,$ where $p$ is a prime number.

\begin{corollary}\label{Cor:Iso vertex in DNI(G); when G is not 2 gebnerated}
Let $G$ be a finite group that is not $2$-generated, and suppose that $n_p\geq 2$ for some prime divisor $p$ of $|G|$. Then the graph $\mathcal{DNI}(G)$ has no isolated vertices.
\end{corollary}
\begin{proof}
Since the identity element $e$ is excluded from the vertex set of $\mathcal{DNI}(G)$ by definition, Proposition~\ref{Prop: IPG(G) is dominatable iff G is cyclic; G is any grp} implies that $\mathcal{DNI}(G)$ has no isolated vertices.
\end{proof}

\begin{corollary}\label{Cor:np=1-isovertex}
Let $G$ be a finite group that is not $2$-generated. Assume that $n_p=1$ for every prime divisor $p$ of $|G|$. Then a vertex $g\in G\setminus \{e\}$ is isolated in $\mathcal{DNI}(G)$ if and only if the order of $g$ is divisible by every prime divisor of $|G|$.
\end{corollary}
\begin{proof}
If $\text{o}(g)$ is divisible by every prime divisor of $|G|$, then it follows from the proof of Proposition~\ref{Prop: IPG(G) is dominatable iff G is cyclic; G is any grp} that $g$ is an isolated vertex of $\mathcal{DNI}(G)$. On the other hand, if there exists a prime divisor $p$ of $|G|$ with $p\nmid \text{o}(g)$, then $g$ is adjacent in $\mathcal{DNI}(G)$ to every element of order $p$. Consequently, $g$ cannot be an isolated vertex.
\end{proof}

The following result is easy to observe.
\begin{lemma}\label{Lemma:Two distinct prime elemnts edge connd in GI(G) complement; G is no 2-generated}
Let $g, h\in G$ such that there is no prime $p$ which is a divisor of both $\text{o}(g)$ and $\text{o}(h).$ Then $g\sim h$ in $\mathcal{DNI}(G)$
\end{lemma}

\begin{lemma}\label{Lemma:p-grp-not-2-generated-structure-DNI-graph}
Let $G$ be a finite $p$-group that is neither cyclic nor generalized quaternion. Then the graph $\mathcal{DNI}(G)$ is a complete $r$-partite graph, where $r \ge 3$ is the number of distinct subgroups of $G$ of order $p$.
\end{lemma}

\begin{proof}
Since $G $ it is neither cyclic nor a generalized quaternion group. Then the number of distinct cyclic subgroups of order $p$ in $G$ is at least $3$. Let 
\[
\langle a_1 \rangle, \ldots, \langle a_r \rangle
\]
be the complete collection of distinct subgroups of $G$ of order $p$, where $r \ge 3$.

For each $i \in [r]$, define
\[
K_{a_i} = \{x \in G \setminus \{e\} : a_i \in \langle x \rangle\}.
\]
Then
\[
G \setminus \{e\} = \bigcup_{i=1}^r K_i,
\]
and moreover, $K_{a_i} \cap K_{a_j} = \emptyset$ whenever $i \neq j$.

Observe that each $K_{a_i}$ induces a complete subgraph in $\mathcal{G}_I(G)$, and there is no edge between any vertex of $K_{a_i}$ and any vertex of $K_{a_j}$ for $i \neq j$. Thus, the result follows.
\end{proof}
\begin{corollary}\label{Cor:connec-p-grp-n2gen}
 Let $G$ be a non $2$-generated group of prime power order. Then $\mathcal{DNI}(G)$ is connected with ${\rm diam}(\mathcal{DNI}(G))=2.$   
\end{corollary}
 \begin{proof}
 As $G$ is a non-$2$-generated group, $\mathcal{DNI}(G)=\bar{\mathcal{G}_I}(G)\vert_{G\setminus\{e\}}$. Hence, Lemma \ref{Lemma:p-grp-not-2-generated-structure-DNI-graph} implies that $\mathcal{DNI}(G)$ is connected with
${\rm diam}(\mathcal{DNI}(G))=2.$
\end{proof}   

\begin{theorem}\label{Thm:GI(G) complement is connd diam leq 3; when G is not 2-generated}
Let $G$ be a finite non-$2$-generated group with $\pi(G)\geq 2$. If $n_p\geq 2$ for some prime divisor $p$ of $|G|$, then $\mathcal{DNI}(G)$ is connected with ${\rm diam}(\mathcal{DNI}(G))\leq 3$.
\end{theorem}
\begin{proof}
Let 
\[
P=\{\, g\in G \mid \text{o}(g) \text{ is prime and } \langle g\rangle \neq \langle h\rangle \text{ for any distinct } g,h\in P \,\}.
\]
First, we prove that the vertices in $P$ form a clique in $\mathcal{DNI}(G)$. 

Let $g,h\in P$. Suppose $\text{o}(g)=p$ and $\text{o}(h)=q$, where $p$ and $q$ are two distinct primes. Then by Lemma~\ref{Lemma:Two distinct prime elemnts edge connd in GI(G) complement; G is no 2-generated}, we have $g\sim h$ in $\mathcal{DNI}(G)$.
Now suppose that $\mathrm{o}(g)=p=\mathrm{o}(h)$ and $\langle g \rangle \neq \langle h \rangle$. Then clearly $g \nsim h$ in $\mathcal{G}_I(G)$. Hence, $g \sim h$ in $\mathcal{DNI}(G)$.

\noindent
Let $x \in G \setminus (P\cup\{e\})$. Suppose there exists a prime $p$ dividing $|G|$ such that $p \nmid \mathrm{o}(x)$. Let $a \in G$ with $\mathrm{o}(a)=p$. Then $a \in P$, and by Lemma~\ref{Lemma:Two distinct prime elemnts edge connd in GI(G) complement; G is no 2-generated}, we have $x \sim a$ in $\mathcal{DNI}(G)$.

Now assume that every prime divisor of $|G|$ also divides $\mathrm{o}(x)$. Since $n_p\ge 2$ for a prime $p$ dividing $|G|,$ we can choose $y \notin \langle x \rangle$ such that $\mathrm{o}(y)=p$. Then $x \sim y$ in $\mathcal{DNI}(G)$. Hence, $\mathcal{DNI}(G)$ is connected. Moreover, from the above argument, it follows that the distance between any two vertices is at most $3$. Consequently,
\[
\mathrm{diam}(\mathcal{DNI}(G)) \leq 3.
\]
Again, it is easy to see that $\mathcal{DNI}(G)$ is not complete. Hence the theorem follows.
\end{proof}

Now, we consider a non-$2$-generated group $G$ with $|\pi(G)| \geq 2$ and $n_p = 1$ for every prime divisor $p$ of $|G|$. Let $I$ denote the set of all isolated vertices of the graph $\mathcal{DNI}(G)$. By adapting the proof of Theorem~\ref{Thm:GI(G) complement is connd diam leq 3; when G is not 2-generated}, we obtain the following result:
\begin{theorem}\label{Thm:connected-having-Isolated-N2Gen}
 The graph $\mathcal{DNI}(G)\vert_{G\setminus\{I\}}$ is connected with ${\rm diam}(\mathcal{DNI}(G)\vert_{G\setminus\{I\}})\leq 3$. 
\end{theorem}

\subsection{2-generated nilpotent groups}
We now examine the connectedness properties of the graph under the assumption that $G$ is a finite $2$-generated $N_I$-group, and $p$ is a prime number.
\begin{lemma}\label{Lemma: Diameter abelian p group}
Let $G$ be a $2$-generated abelian $p$-group. Then we have the following:
\begin{enumerate}
    \item If $G$ is elementary abelian, then $\mathcal{DNI}(G)$ is the null graph with vertex set $G\setminus \{e\}$.
    \item If $G$ is not elementary abelian, then $\mathcal{DNI}(G)$ is connected with ${\rm diam}(\mathcal{DNI}(G))=2$.  
\end{enumerate}
\end{lemma}
\begin{proof}
Let $G$ be a $2$-generated abelian $p$-group. Since $G$ is a $N_I$-group, Corollary \ref{Cor: 2-generated p group} implies that $G\cong C_{p^r}\times C_{p^r}$, for some $r\in \mathbb{N}$. 
\begin{enumerate}
    \item[$\mathbf{\underline{r=1:}}$] If $r=1$, then $G$ is an elementary abelian $p$-group. Let $x,y$ be two non-trivial elements of $G$. Then $o(x)=o(y)=p$. So, either $\langle x\rangle=\langle y\rangle$ or $\langle x\rangle\cap\langle y\rangle=\{e\}$. In first case, $x\sim y$ in $\mathcal{G}_I(G)$ and in the later case, $\langle x,y \rangle=G$, i.e.~$x\sim y$ in ${\rm Gen}(G)$. So, $x\nsim y$ in $\mathcal{DNI}(G)$. 
    
    \item[$\mathbf{\underline{r\geq 2:}}$] For $r\geq 2$, we have $\lvert G\rvert=p^{2r}$ let $x,y\in G$ be two non-trivial elements of $G$. Then $o(x),o(y)\leq p^r$. If $x\sim y$ in $\mathcal{G}_I(G)$, let $\langle z\rangle$ be the $p$-order subgroup of $\langle x\rangle\cap\langle y\rangle$. $G$ is non-cyclic, choose $z'\in G$ such that $o(z')=p$ and $\langle z\rangle\cap\langle z'\rangle=\{e\}$. Then $x\sim z'\sim y$ in $\bar{\mathcal{G}_I}(G)$. Now, $\lvert\langle x,z'\rangle\rvert=o(x)p<p^{2r}$ and $\lvert\langle y,z'\rangle\rvert=o(x)p<p^{2r}$, since $r\geq 2$. So, $\langle x,z'\rangle\neq G$ and $\langle y,z'\rangle\neq G$. Thus, $x\sim z'\sim y$ in $\mathcal{DNI}(G)$.

    Now, let $x\nsim y$ in $\mathcal{G}_I(G)$. Let $\tilde{x}$ and $\tilde{y}$ be two elements in $G$ such that $\langle\tilde{x}\rangle$ and $\langle\tilde{y}\rangle$ are the unique subgroups of order $p$ in $\langle x\rangle$ and $\langle y\rangle$, respectively. We have $\langle\tilde{x}\rangle\cap\langle\tilde{y}\rangle=\{e\}$. Since $[\tilde{x},\tilde{y}]=e$, we have $o(\tilde{x}\tilde{y})=p$ and $\langle \tilde{x}\rangle\cap\langle\tilde{x}\tilde{y}\rangle=\langle \tilde{y}\rangle\cap\langle\tilde{x}\tilde{y}\rangle=\{e\}$. Since $o(\tilde{x}\tilde{y})=p$, we have $\langle x\rangle\cap\langle\tilde{x}\tilde{y}\rangle=\langle y\rangle\cap\langle\tilde{x}\tilde{y}\rangle=\{e\}$. So, $x\sim \tilde{x}\tilde{y}\sim y$ in $\bar{\mathcal{G}_I}(G)$. Again, since $r\geq 2$ and $o(\tilde{x}\tilde{y})=p$, by similar arguments used before, we have Thus, $x\sim \tilde{x}\tilde{y}\sim y$ in $\mathcal{DNI}(G)$.
\end{enumerate}
\end{proof}

The following result is an easy consequence of Lemma \ref{Lemma:p-grp-not-2-generated-structure-DNI-graph}. 
\begin{lemma}\label{Lemma:G_I complement diameter 2}
Let $G$ be a $2$-generated $p$-group, and $x,y$ be two nontrivial elements of $G$. Then there exists a nontrivial element $z\in G$ such that $x\sim z\sim y$ is a path in $\bar{\mathcal{G}_I}(G)$.
\end{lemma}

\begin{lemma}\label{Lemma: non-cyclic <center, commutator> connectivity}
Let $G$ be a $2$-generated non-abelian $p$-group, such that $\langle Z(G),G'\rangle$ is neither cyclic, nor generalized quaternion. Then $\mathcal{DNI}(G)$ is connected with ${\rm diam}(\mathcal{DNI}(G))=2$.
\end{lemma}
\begin{proof}
Let $G$ be a $2$-generated non-abelian $p$-group, such that $\hat{G}:=\langle Z(G),G'\rangle$ is neither cyclic, nor generalized quaternion. Then Theorem \ref{Theorem: Gorenstein} implies that $\hat{G}$ has a subgroup isomorphic to $C_p\times C_p$. Let $x,y$ be two non-trivial elements of $G$. Then there exists $z\in \hat{G}$ such that $o(z)=p$ and $\langle x\rangle\cap\langle z\rangle=\langle y\rangle\cap\langle z\rangle=\{e\}$. In addition, Lemma \ref{Lemma: one generator inside <commutator,center> implies proper subgroup} implies that $\langle x,z\rangle\neq G$ and $\langle y,z\rangle\neq G$. So, $x\sim z\sim y$ is a path in $\mathcal{DNI}(G)$. Now, Theorem \ref{Theorem: Completeness and nullness} implies that $\mathcal{DNI}(G)$ is not complete. Thus, the result holds.
\end{proof}

\begin{proposition}\label{Prop: 2 group diameter 4 connectivity}
Let $G$ be a $2$-generated non-abelian $2$-group. Then $\mathcal{DNI}(G)$ is connected if and only if $G\ncong D_{2^n}$. Moreover, if $\mathcal{DNI}(G)$ is connected, then ${\rm diam}(\mathcal{DNI}(G))\leq 4$.    
\end{proposition}
\begin{proof}
Let $G$ be a $2$-generated non-abelian $2$-group. Clearly, $\lvert G\rvert\geq 8$. If $G\cong D_{2^n}$, then it is easy to see that all the elements   of order $2^{n-1}$ are the isolated vertices in $\mathcal{DNI}(G)$. So, $\mathcal{DNI}(G)$ is not connected. Otherwise, if $G\ncong D_{2^n}$, let $a,b\in G$ be two arbitrary distinct elements. Since $G$ is neither cyclic nor generalized quaternion, Theorem \ref{Theorem: Gorenstein} yields that $G$ has a subgroup isomorphic to $C_2\times C_2$. Then there exists $z\in G$, such that $o(z)=2$ and $\langle z\rangle\cap\langle a\rangle=\langle z\rangle\cap\langle b\rangle=\{e\}$. If $\langle a,z\rangle\neq G$, then we have $a\sim z$ in $\mathcal{DNI}(G)$. Otherwise, if $\langle a,z\rangle=G$, then we have $G=\langle a,az\rangle=\langle z,az\rangle$. Clearly $\langle a\rangle\cap\langle az\rangle=\langle z\rangle\cap\langle az\rangle=\{e\}$, since $G$ is a $N_I$-group. If $o(az)=2$, then $G=\langle z,az\rangle$ implies that $G\cong D_{2^n}$, which is contradictory. So, $o(az)\geq 4$. Then Lemma \ref{Lemma: langle x^p,y_ranglesubsetneq G} implies that $\langle a,(az)^2\rangle\subsetneq G$ and $\langle z,(az)^2\rangle\subsetneq G$, hence $a\sim (az)^2\sim z$ is a path in $\mathcal{DNI}(G)$. Similarly, it can be shown that $b$ is connected to $z$ via a path of length at most $2$ in $\mathcal{DNI}(G)$. Thus, $\mathcal{DNI}(G)$ is connected with ${\rm diam}(\mathcal{DNI}(G))\leq 4$. This concludes the proof.
\end{proof}

\begin{lemma}\label{Lemma: path implies normal subgroup}
Let $G$ be a $2$-generated non-abelian $p$-group, for an odd prime $p$, such that $G=\langle a,b\rangle$ and $\langle a\rangle\cap G'\neq \{e\}$. If $\langle a,G'\rangle$ is not cyclic, then there exist $z\in\langle a,G'\rangle$ and $a'\in \langle a\rangle\cap G'$ with $o(a')=p$, such that $a\sim z\sim a'\sim b$ is a path in $\mathcal{DNI}(G)$.    
\end{lemma}
\begin{proof}
Let $G$ be a $2$-generated non-abelian $p$-group, for an odd prime $p$, such that $G=\langle a,b\rangle$ and $\langle a\rangle\cap G'\neq \{e\}$. If $\langle a,G'\rangle$ is not cyclic, then by Theorem \ref{Theorem: Gorenstein}, $\langle a,G'\rangle$ has a subgroup isomorphic to $C_p\times C_p$. So, there exists $z\in \langle a,G'\rangle$ such that $\langle a\rangle\cap\langle z\rangle=\{e\}$. Now $\langle a,z\rangle \subseteq \langle a,G'\rangle\subsetneq G$, where the last proper containment follows from Lemma \ref{Lemma: Cyclic abelianization}. So, $a\sim z$ in $\mathcal{DNI}(G)$. Let $a'\in \langle a\rangle\cap G'$ such that $o(a')=p$. Then $\langle a'\rangle\cap\langle z\rangle=\{e\}$ and $\langle a',z\rangle \subseteq \langle a,G'\rangle\subsetneq G$. Then $a'\sim z$ in $\mathcal{DNI}(G)$. Now, since $G$ is a $N_I$-group, we have that $\langle a'\rangle\cap\langle b\rangle=\{e\}$, and from Lemma \ref{Lemma: Cyclic abelianization} we have that $\langle a', b\rangle\subsetneq G$. So, $a'\sim b$ in $\mathcal{DNI}(G)$. Thus, we have a path $a\sim z\sim a'\sim b$ in $\mathcal{DNI}(G)$. This concludes the proof.
\end{proof}

\begin{lemma}\label{Lemma: connected when non metacyclic}
Let $G$ be a finite $2$-generated non-abelian $p$-group, for an odd prime $p$, such that $G$ is not metacyclic. Then $\mathcal{DNI}(G)$ is connected with ${\rm diam}(\mathcal{DNI}(G))\leq 4$.    
\end{lemma}
\begin{proof}
Let $G$ be a finite $2$-generated non-abelian $p$-group, for an odd prime $p$, such that $G$ is not metacyclic. Let $a,b\in G$ be two arbitrary distinct non-trivial elements. If $G'\cap \langle a\rangle=G'\cap \langle b\rangle=\{e\}$, choose any non-trivial element $z\in G'$. Then $\langle a\rangle\cap\langle z\rangle=\langle b\rangle\cap\langle z\rangle=\{e\}$ and as a consequence of Lemma \ref{Lemma: Cyclic abelianization}, we have $\langle a,z\rangle\subsetneq G$ and $\langle b,z\rangle\subsetneq G$. So, $a\sim z\sim b$ is a path in $\mathcal{DNI}(G)$. Otherwise, if one of $G'\cap \langle a\rangle$ and $G'\cap \langle b\rangle$ is non-trivial, without loss of generality, we assume that $G'\cap \langle a\rangle\neq \{e\}$. 

Suppose $\langle a\rangle\cap \langle b\rangle=\{e\}$. If $\langle a,b\rangle\neq G$, then $a\sim b$ in $\mathcal{DNI}(G)$. Let $\langle a,b\rangle=G$. If $H:=\langle a,G'\rangle$ is cyclic, we have that $H$ is a cyclic normal subgroup of $G$, and $G/H=\langle a,b\rangle/H$ is cyclic, which is contradictory to our assumption that $G$ is not metacyclic. Thus $H$ is non-cyclic. Then Lemma \ref{Lemma: path implies normal subgroup} implies that there exists a path of length $3$ between $a$ and $b$.  

Now, suppose $\langle a\rangle\cap \langle b\rangle\neq\{e\}$. Then $G'\cap \langle b\rangle\neq \{e\}$. By Theorem \ref{Theorem: Gorenstein}, we have that $G$ has a subgroup isomorphic to $C_p\times C_p$. So, there exists $d\in G'$ such that $o(d)=p$ and $\langle a\rangle\cap\langle d\rangle=\langle b\rangle\cap\langle d\rangle=\{e\}$. If $\langle a,d\rangle\subsetneq G$ and $\langle b,d\rangle\subsetneq G$ then $a\sim d\sim b$ is a path in $\mathcal{DNI}(G)$. Otherwise, let one of $\langle a,d\rangle$ and $\langle b,d\rangle$ is $G$. Without loss of generality, we assume that $\langle a,d\rangle=G$. Define $H:=\langle a,G'\rangle$. If $H$ is cyclic, we have that $H$ is a cyclic normal subgroup of $G$, and $G/H=\langle a,d\rangle/H$ is cyclic, which is contradictory to our assumption that $G$ is not metacyclic. Thus $H$ is non-cyclic. Then by Lemma \ref{Lemma: path implies normal subgroup} we have a path $a\sim z_1\sim a'\sim d$ in $\mathcal{DNI}(G)$, where $z_1\in H$ and $a'\in G'\cap \langle a\rangle$ with $o(a')=p$. Similarly, it can be shown that $b\sim d$ or we have a path $b\sim z_2\sim b'\sim d$ in $\mathcal{DNI}(G)$, where $z_2\in \langle b,G'\rangle$ and $b'\in G'\cap \langle b\rangle$ with $o(b')=p$. Now, since $a'\in G'\cap \langle a\rangle$ with $o(a')=p$ and $b'\in G'\cap \langle b\rangle$ with $o(b')=p$, and $\langle a\rangle\cap \langle b\rangle\neq\{e\}$, we have $b'\in G_{a'}$. Hence, Lemma \ref{Lemma: x in Ga and y in Gb, if x sim y then any two elmnts in Ga and Gb edge connd in DNI(G)} implies that $a\sim z_1\sim b'\sim z_2\sim b$. This concludes the proof.
\end{proof}

\begin{proposition}\label{Prop: Metacyclic connectedness}
Let $G$ be a non-abelian metacyclic $p$-group, for an odd prime $p$. Then $\mathcal{DNI}(G)$ is connected with ${\rm diam}(\mathcal{DNI}(G))\leq 4$. 
\end{proposition}
\begin{proof}
Let $G$ be a non-abelian metacyclic $p$-group, for an odd prime $p$. Define the set $\mathcal{S}:=\{s\in G~\vert~\langle s\rangle \text{ is normal in } G \text{ and }G/\langle s\rangle\text{ is cyclic}\}$, and the set $\mathcal{T}:=\{t\in \mathcal{S}~\vert~o(t)\geq o(s)\text{ for all }s\in S\}$. Since $G$ is metacyclic, we have that $\mathcal{S}$ is non empty, so is $\mathcal{T}$. Let $a\in \mathcal{T}$ be an arbitrary element. Let $G/\langle a\rangle=\langle\tilde{b}\rangle$ and $b$ is a pre-image of $\tilde{b}$. Then $G=\langle a,b\rangle$, hence $\langle a\rangle\cap\langle b\rangle=\{e\}$.  If $G/\langle a\rangle\cong C_p$, then Corollary \ref{Cor: metacyclic p index not NI group} indicates that $G$ is not a $N_I$-group. So, we have $\lvert G/\langle a\rangle \rvert>p$. Define the set $\mathcal{Q}=\{g\in G\setminus\langle a\rangle~\vert~o(g)=p\}$. Then Theorem \ref{Theorem: Gorenstein} implies that $\mathcal{Q}$ is non empty. Choose an arbitrary element $c\in \mathcal{Q}$. Then we have $\langle a\rangle\cap\langle c\rangle=\{e\}$ and $\lvert \langle a,c\rangle\rvert=\lvert \langle a\rangle\rvert\times p<\lvert G\rvert$, i.e.~$\langle a,c\rangle\subsetneq G$. So, $a\sim c$ in $\mathcal{DNI}(G)$. Let $a_1\in \langle a\rangle$ with $o(a_1)=p$. Clearly $a_1\in G'$, since $G'\in \langle a\rangle$ and $G$ is non-abelian. So, Lemma \ref{Lemma: Cyclic abelianization} implies that $\langle a_1,c\rangle\neq G$. Hence, we have $a_1\sim c$ in $\mathcal{DNI}(G)$.
 
Recall the set $K_{a_1} = \{x \in G \setminus \{e\} : a_1 \in \langle x \rangle\}$. Let $x\in G\setminus K_{a_1}$. Then $\langle x\rangle\cap\langle a_1\rangle=\{e\}$ and Lemma \ref{Lemma: Cyclic abelianization} implies that $\langle x,a_1\rangle\neq G$. So, $x\sim a_1$ in $\mathcal{DNI}(G)$. Hence, $x\sim a_1\sim c\sim a$ is a path in $\mathcal{DNI}(G)$. Now, let $y\in K_{a_1}$ be an arbitrary element. Then $\langle y\rangle\cap\langle c\rangle=\{e\}$. If $y\in \langle a\rangle$, then Lemma \ref{Lemma: langle x^p,y_ranglesubsetneq G} implies that $\langle y,c\rangle\subseteq\langle a,c\rangle\subsetneq G$. Hence, $y\sim c$ in $\mathcal{DNI}(G)$. So, $y\sim c\sim a$. If $y\notin \langle a\rangle$, then either $G'\subseteq \langle y\rangle$ or $G'\nsubseteq\langle y\rangle$. If $G'\subseteq \langle y\rangle$, then $o(y)\leq o(a)$, since $\lvert \langle y,b\rangle\rvert=o(y)\times o(b)\leq \lvert G\rvert=o(a)\times o(b)$. So $\lvert \langle y,c\rangle\rvert\leq \lvert \langle a,c\rangle\rvert<\lvert G\rvert$, i.e.~$\langle y,c\rangle\subsetneq G$. So, $y\sim c\sim a$ in $\mathcal{DNI}(G)$. If $G'\nsubseteq\langle y\rangle$, then $\langle y, G'\rangle$ is not cyclic. If $\langle y,c\rangle\neq G$, then $y\sim c$ in $\mathcal{DNI}(G)$. So, $y\sim c\sim a$ is a path in $\mathcal{DNI}(G)$. Otherwise, if $\langle y,c\rangle=G$, then Lemma \ref{Lemma: path implies normal subgroup} implies that there exists $z\in\langle y, G'\rangle$ such that $y\sim z\sim a_1\sim c$ is a path in $\mathcal{DNI}(G)$. Hence, we have $y\sim z\sim a_1\sim c\sim a$ is a path in $\mathcal{DNI}(G)$. So, $\mathcal{DNI}(G)$ is connected.
 
Now, let $u,v$ be two arbitrary non-trivial elements in $G$. If both $u,v\in G\setminus K_{a_1}$, then we have $u\sim a_1\sim c\sim a$ and $v\sim a_1\sim c\sim a$ in $\mathcal{DNI}(G)$. Thus, $u\sim a_1\sim v$ in $\mathcal{DNI}(G)$. If $u\in G\setminus K_{a_1}$ and $v\in K_{a_1}$, then $u\sim a_1\sim c\sim a$ and $v\sim c\sim a$ or $v\sim z_v\sim a_1\sim c\sim a$ in $\mathcal{DNI}(G)$, for some $z_v\in G$. Then, either $u\sim a_1\sim c\sim v$ of $u\sim a_1\sim z_v\sim v$ in $\mathcal{DNI}(G)$. If both $u,v\in K_{a_1}$, then we have $u\sim c\sim a$ and $v\sim c\sim a$, or $u\sim c\sim a$ and $v\sim z_v\sim a_1\sim c\sim a$, or $v\sim c\sim a$ and $u\sim z_u\sim a_1\sim c\sim a$, or $u\sim z_u\sim a_1\sim c\sim a$ and $v\sim z_v\sim a_1\sim c\sim a$ in $\mathcal{DNI}(G)$, for some $z_u\in G$. So, we have $u\sim c\sim v$, or $u\sim c\sim a_1\sim z_v\sim v$, or $u\sim z_u\sim a_1\sim c\sim v$, or $u\sim z_u\sim a_1\sim z_v\sim v$ in $\mathcal{DNI}(G)$. Thus, we have ${\rm diam}(\mathcal{DNI}(G))\leq 4$. This concludes the proof.
\end{proof}

\begin{proposition}\label{Prop: diameter 2-generated nilpotent non-p $N_I$-group}
Let $G$ be a $2$-generated nilpotent group with $\pi(G)\geq 2$. Then $\mathcal{DNI}(G)$ is connected, and ${\rm diam}(\mathcal{DNI}(G))=2$.    
\end{proposition}
\begin{proof}
Let $G$ be a $2$-generated nilpotent group with $\pi(G)\geq 2$. Assume $G=G_1\times G_2\times\cdots\times G_n$, where $G_1,G_2,\ldots,G_n$ are Sylow $p$-subgroups of $G$, and $n\geq 2$. Then Lemma \ref{2-generated NI} implies that $G_i$ is a $2$-generated $N_I$-group, for all $i\in[n]$. Let $\mathfrak{g}_1\in G_1$ and $\mathfrak{g}_2\in G_2$ such that both $\mathfrak{g}_1$ and $\mathfrak{g}_2$ are of prime order. Then $o(\mathfrak{g}_1\mathfrak{g}_2)\geq 6$, hence $\mathfrak{g}_1\mathfrak{g}_2$ and $(\mathfrak{g}_1\mathfrak{g}_2)^2$ are two nontrivial elements of $G$ such that $\mathfrak{g}_1\mathfrak{g}_2\nsim(\mathfrak{g}_1\mathfrak{g}_2)^2$ in $\mathcal{DNI}(G)$. So, $\mathcal{DNI}(G)$ is not complete.
 
 Let $(a_1,b_1)$ and $(a_2,b_2)$ be two non-trivial elements of $G_1\times G_2$. Without loss of generality, we assume that $a_1\neq e$. Now, \cite[Theorem 4.1 and 4.3]{intersectionpwegraphb3} implies that for a finite $p$-group $H$, $\bar{\mathcal{G}_I}(H)$ has non identity isolated vertices if and only if $H$ is a cyclic group or a generalized quaternion group. Since both cyclic groups and generalized quaternion groups are non $N_I$-groups, we can see that both $\bar{\mathcal{G}_I}(G_1)$ and $\bar{\mathcal{G}_I}(G_2)$ have no isolated vertex other than the identity. Then Lemma \ref{Lemma:G_I complement diameter 2} implies that there exists a nontrivial element $a\in G_1$ such that $a_1\sim a\sim a_2$ in $\bar{\mathcal{G}_I}^e(G_1)$. So from Theorem \ref{Thm: direct product}, it can be seen that $(a_1,b_1)\sim (a,e)\sim(a_2,b_2)$ in $\mathcal{DNI}^e(G_1\times G_2)$. Again, since $(a_1,b_1), (a,e),(a_2,b_2)$ are all nontrivial elements of $G_1\times G_2$, we have $(a_1,b_1)\sim (a,e)\sim(a_2,b_2)$ in $\mathcal{DNI}(G_1\times G_2)$. This also proves the result for $n=2$.

 Now, let $3\leq m\leq n$ be a natural number such that for any two nontrivial elements $g_1,g_2\in G_1\times G_2\times\cdots\times G_{m-1}$, there exists a nontrivial element $g\in G_1\times G_2\times\cdots\times G_{m-1}$ such that $g_1\sim g\sim g_2$ in $\mathcal{DNI}(G_1\times G_2\times\cdots\times G_{m-1})$. Denote $\mathfrak{G}:=G_1\times G_2\times\cdots\times G_{m-1}$. Consider the group $G=\mathfrak{G}\times G_m$, where $G_m$ is a finite $2$-generated $N_I$-group of prime power order, and ${\rm gcd}(\lvert\mathfrak{G}\rvert,\lvert G_m\rvert)=1$. Let $(g_1,h_1),(g_2,h_2)$ be two non-trivial elements of $\mathfrak{G}\times G_m$. Without loss of generality, we assume that $g_1\neq e$. Then \cite[Theorem 4.1 and 4.3]{intersectionpwegraphb3} indicate that if $g_2=e$, then there exists a nontrivial $g\in \mathfrak{G}$, such that $g_1\sim g\sim g_2$ in $\bar{\mathcal{G}_I}^e(\mathfrak{G})$, and the inductive hypothesis indicate that if $g_2\neq e$, there exists a nontrivial $g\in \mathfrak{G}$, such that $g_1\sim g\sim g_2$ in $\mathcal{DNI}^e(\mathfrak{G})$. So, Theorem \ref{Thm: direct product} implies that $(g_1,h_1)\sim (g,e)\sim(g_2,h_2)$ in $\mathcal{DNI}^e(\mathfrak{G}\times G_m)$. Again, since $(g_1,h_1), (g,e),(g_2,h_2)$ are all nontrivial elements of $G$, we have $(g_1,h_1)\sim (g,e)\sim(g_2,h_2)$ in $\mathcal{DNI}(\mathfrak{G}\times G_m)$. This concludes the proof. 
\end{proof}

The following theorem summarizes all the results obtained in this section, and classifies all finite $2$-generated nilpotent $N_I$-groups for which the difference graph is connected.
\begin{theorem}\label{Theorem: connectedness 2-generated}
Let $G$ be a $2$-generated nilpotent group. Then we have the following:
\begin{enumerate}
    \item If $\pi(G)\geq 2$, then $\mathcal{DNI}(G)$ is connected with ${\rm diam}(\mathcal{DNI}(G))=2$. 
    \item If $G$ is an elementary abelian $p$-group, then $\mathcal{DNI}(G)$ is the null graph with vertex set $G\setminus \{e\}$.
    \item If $G$ is an abelian $p$-group with ${\rm exp}(G)>p$, then $\mathcal{DNI}(G)$ is connected, and we have ${\rm diam}(\mathcal{DNI}(G))=2$.
    \item If $G$ is a non-abelian $p$-group for an odd prime $p$, such that $G'$ is non-cyclic, then $\mathcal{DNI}(G)$ is connected with ${\rm diam}(\mathcal{DNI}(G))=2$.
    \item If $G$ is a non-abelian $p$-group for an odd prime $p$, such that $G'$ is cyclic, then $\mathcal{DNI}(G)$ is connected with ${\rm diam}(\mathcal{DNI}(G))\leq 4$.
    \item If $G$ is a non-abelian $2$-group, then $\mathcal{DNI}(G)$ is connected if and only if $G\ncong D_{2^n}$. Moreover, if $\mathcal{DNI}(G)$ is connected, then ${\rm diam}(\mathcal{DNI}(G))\leq 4$.
\end{enumerate}
\end{theorem}

\section{Eulerian}\label{Sec: Eulerian}
The objective of this section is to characterize the Eulerian nature of the difference graph. We start with a preliminary proposition that serves as a key ingredient in the proofs that follow.

\begin{proposition}\label{Prop:Deg(a) is even in DNI(G); when |G| is odd}
Let $G$ be a group of odd order. Then $\deg(a)$ is even for every non-isolated vertex $a\in V(\mathcal{DNI}(G)).$     
\end{proposition}
\begin{proof}
Let $a\sim x_1$ in $\mathcal{DNI}(G).$ Then by Lemma \ref{Lemma: x in Ga and y in Gb, if x sim y then any two elmnts in Ga and Gb edge connd in DNI(G)}, $a$ is edge connected to all the elements in $G_{x_1}.$  Let $x_1, \dots, x_m\in  V(\mathcal{DNI}(G))\setminus G_a$ be such that $a\sim x_i$ and $G_{x_i}\cap G_{x_j}=\emptyset,$ for all $i\neq j.$ Then the degree of the vertex $a$ is 
\begin{equation}\label{Eqn:form of Deg(a), for any a; G is not 2gented,Z(G)=e}
\text{deg}(a)=\phi(\text{o}(x_1))+\cdots+\phi(\text{o}(x_m)).
\end{equation}
It is given that $|G|$ is odd. So, the order of each non-identity element is odd and hence $\phi(\text{o}(x_i))$ is even, for all $i\in [m].$ Consequently, $\deg(a)$ is even. 
\end{proof}
\begin{proposition}\label{Prop:DNI(G) is not Eulerian, G is not 2-generated, Z(G)=e and |G| is even}
Let $G$ be a group of even order with $\pi(G)\geq 2$. Then $\mathcal{DNI}(G)$ is not Eulerian.
\end{proposition}
\begin{proof}
The order of the group is even. So, $a\in G$ such that $\text{o}(a)=2,$ and $a$ is not an isolated vertex in $\mathcal{DNI}(G)$. Let the number of distinct subgroups of order $2$ in $G$ be $r,$ where $r\geq 1.$ Let $a, b\in G$ be two distinct elements such that $\text{o}(a)=2=\text{o}(b).$ Then $\langle a\rangle\cap\langle b \rangle=\{e\}.$ As a result, $a\sim b$ in $\mathcal{DNI}(G)$. That is, $a$ is edge connected to all other $2$-ordered elements. So according to Equation \eqref{Eqn:form of Deg(a), for any a; G is not 2gented,Z(G)=e}, 
\begin{equation}\label{Eqn:Deg(a) is even iff, when Gis not 2gerted, Z(G)=e, ordr of a is 2}
\deg(a)=\sum\limits_{b\in G\setminus \{a\}:\text{o}(b)=2} \phi(2)+\sum\limits_{x\in G\setminus\{e\}:\text{o}(x)\neq 2 \text{ and }a\sim x}\phi(\text{o}(x)).    
\end{equation}
Now, the term $\left(\sum\limits_{x\in G\setminus\{e\}:\text{o}(x)\neq 2 \text{ and }a\sim x}\phi(\text{o}(x))\right)$ in Equation \eqref{Eqn:Deg(a) is even iff, when Gis not 2gerted, Z(G)=e, ordr of a is 2} is even (as $\phi(\text{o}(x))$ is even). Also, $\left(\sum\limits_{b\in G\setminus\{a\}:\text{o}(b)=2} \phi(2)\right)$ is even if and only if $r$ is odd. That is, $\deg(a)$ is even if and only if $r$ is odd. It is given that $|\pi(G)|\ge 2.$ Therefore, there exists an odd prime $p$ that divides $|G|.$ Let $c\in G$ be such that $\text{o}(c)=p.$ By Lemma \ref{Lemma:Two distinct prime elemnts edge connd in GI(G) complement; G is no 2-generated}, $c$ is edge connected with all the elements ordered $2$ of $G.$ As a result, 
\begin{equation}\label{Eqn:Deg(a) is even iff, when Gis not 2gerted, Z(G)=e, ordr of a is p, p is odd prime}
\deg(c)=r+\sum\limits_{x\in G\setminus\{e\}:\text{o}(x)\neq 2\text{ and }c\sim x}\phi(\text{o}(x)).    
\end{equation}
Now, $\left(\sum\limits_{x\in G\setminus\{e\}:\text{o}(x)\neq 2\text{ and } c\sim x}\phi(\text{o}(x))\right)$ is even and $r$ is odd implying that $\deg(c)$ is odd. That is, $\deg(a)$ is even if and only if $\deg(c)$ is odd. This completes the proof.
\end{proof}

\subsection{Non 2-generated groups}
Now, with the above results in hand, we have the following.
\begin{theorem}\label{Thm:DNI(G) is Eulerian iff, when G is not 2-generated}
Let $G$ be a non $2$-generated group with $n_p>1,$ for some prime $p$ dividing $|G|.$ Then $\mathcal{DNI}(G)$ is Eulerian if and only if either $|G|$ is odd, or $G$ is a $2$-group.    
\end{theorem}
\begin{proof}
Let $G$ be a non $2$-generated group. So it is an $N_I$-group. Since $n_p>1,$ by Corollary \ref{Cor:connec-p-grp-n2gen} and Theorem \ref{Thm:GI(G) complement is connd diam leq 3; when G is not 2-generated}, $\mathcal{DNI}(G)$ is connected. If $|G|$ is odd, then Proposition \ref{Prop:Deg(a) is even in DNI(G); when |G| is odd} implies that $\mathcal{DNI}(G)$ is Eulerian. Also, Proposition \ref{Prop:DNI(G) is not Eulerian, G is not 2-generated, Z(G)=e and |G| is even} implies that when $\pi(G)\geq 2$ and $|G|$ is even, then $\mathcal{DNI}(G)$ is not Eulerian. Now, if $G$ is a $2$-group, let $r$ be the number of elements of order $2$ in $G$. Then Lemma \ref{Lemma:number of elements of order p in a p group} implies that $r$ is odd. Let $g$ be a non-trivial element of $G$, and $x\in G$ is an element of order $2$. Then $g\sim x$ in $\mathcal{DNI}(G)$ if and only if $x\notin \langle g\rangle$. So, the number of elements of order $2$ in $G$ that are adjacent to $g$ is $r-1$, which is even. Now, from Equation (\ref{Eqn:Deg(a) is even iff, when Gis not 2gerted, Z(G)=e, ordr of a is 2}) we have $d_{\mathcal{DNI}(G)}(g)$ is even. Thus, $\mathcal{DNI}(G)$ is Eulerian. This concludes the proof.    
\end{proof}

\begin{theorem}
Let $G$ be a non $2$-generated group with $n_p=1,$ for every prime $p$ dividing $|G|.$ Then $\mathcal{DNI}(G)\vert_{G\setminus\{I\}}$ is Eulerian if and only if $|G|$ is odd with $\pi(G)\ge 2,$ where $I$ is the set of all isolated vertices in $\mathcal{DNI}(G).$ 
\end{theorem}
\begin{proof}
By Theorem~\ref{Thm:connected-having-Isolated-N2Gen}, the induced subgraph $\mathcal{DNI}(G)\vert_{G\setminus\{I\}}$
is connected. If $\pi(G)=1$, then $G$ is a $2$-generated group, but it is given that $G$ is a non $2$-generated group. Hence, we have that $\pi(G)\geq 2$. The remainder of the proof is an immediate consequence of Propositions~\ref{Prop:Deg(a) is even in DNI(G); when |G| is odd} and~\ref{Prop:DNI(G) is not Eulerian, G is not 2-generated, Z(G)=e and |G| is even}.
\end{proof}

\subsection{2-generated nilpotent groups}
\begin{proposition}\label{Prop:2-generated nilpotent non p group eulerian}
Let $G$ be a $2$-generated nilpotent group with $\pi(G)\geq 2$. Then $\mathcal{DNI}(G)$ is Eulerian if and only if $\lvert G\rvert$ is odd.
\end{proposition}
\begin{proof}
Let $G$ be a $2$-generated nilpotent $N_I$-group with $\pi(G)\geq 2$. So, Proposition \ref{Prop: diameter 2-generated nilpotent non-p $N_I$-group} implies that $\mathcal{DNI}(G)$ is connected. If $\lvert G\rvert$ is odd, then Proposition \ref{Prop:Deg(a) is even in DNI(G); when |G| is odd} shows that $\mathcal{DNI}(G)$ is Eulerian. If $\lvert G\rvert$ is even, let $G=H\times K$ with $H$ being the Sylow $2$-subgroup of $G$. Let $x\in H$ with $o(x)=2$, and $y\in K$ with $o(y)=p$ for some odd prime $p$ dividing $\lvert K\rvert$. Since, Lemma \ref{2-generated NI} implies that all Sylow $p$-subgroups of $G$ are $2$-generated, both $H$ and $K$ are $2$-generated. So, for any $z\in G$, we have $\langle x,z\rangle\subsetneq G$ and $\langle y,z\rangle\subsetneq G$. Hence, $N_{\mathcal{DNI}(G)}(x)=N_{\bar{\mathcal{G}_I}(G)}(x)$ and $N_{\mathcal{DNI}(G)}(y)=N_{\bar{\mathcal{G}_I}(G)}(y)$. Then using arguments similar to the proof of Proposition \ref{Prop:DNI(G) is not Eulerian, G is not 2-generated, Z(G)=e and |G| is even}, it can be shown that $\mathcal{DNI}(G)$ is not Eulerian. This concludes the proof.    
\end{proof}
Our next result is in the same theme of Theorem \ref{Thm:DNI(G) is Eulerian iff, when G is not 2-generated}. Since we have already classified $2$-generated nilpotent groups such that the difference graph is connected, the next result only considers the groups for which the difference graph is connected, and study their Eulerianess.
\begin{theorem}\label{Thm:DNI(G) is Eulerian iff, when G is 2-generated}
Let $G$ be a $2$-generated nilpotent group such that $\mathcal{DNI}(G)$ is connected. Then $\mathcal{DNI}(G)$ is Eulerian if and only if either $\lvert G\rvert$ is odd, or $G$ is a $2$-group.    
\end{theorem}
\begin{proof}
If $G$ is not a $p$-group, then Proposition \ref{Prop:2-generated nilpotent non p group eulerian} implies that $\mathcal{DNI}(G)$ is Eulerian if and only if $\lvert G\rvert$ is odd. Also, Proposition \ref{Prop:Deg(a) is even in DNI(G); when |G| is odd} implies that if $G$ is a $p$-group for some odd prime $p$ and $\mathcal{DNI}(G)$ is connected, then $\mathcal{DNI}(G)$ is Eulerian.

Now, if $G$ is a $2$-group, let $r$ denote the number of elements of order $2$ in $G$. Then Lemma \ref{Lemma:number of elements of order p in a p group} implies that $r$ is odd. If $G$ is abelian, then Lemma \ref{Lemma: Diameter abelian p group} says that $\mathcal{DNI}(G)$ is connected if and only if $G$ is not elementary abelian. Again, from Corollary \ref{Cor: 2-generated p group} and arguments used in the proof of Theorem \ref{Thm:DNI(G) is Eulerian iff, when G is 2-generated}, it can be seen that when $G$ is not elementary abelian, then $\mathcal{DNI}(G)$ is Eulerian.

Let $G$ be non-abelian $2$-group, and $g$ be a non-trivial element of $G$. If $\langle g,x\rangle\subsetneq G$ for every $x\in G$ with $o(x)=2$, then the number of elements of order $2$ in $G$ that are adjacent to $g$ in $\mathcal{DNI}(G)$ is $r-1$. Otherwise, if $\langle g,y\rangle=G$ for some $y\in G$ with $o(y)=2$, define the non-empty set $S_g:=\{y\in G~\vert~o(y)=2,~\langle g,y\rangle=G\}$. Then $y'\in S_g$, for any $y'\in G$ such that $y'$ is a conjugate of $y$ by some element of $\langle g\rangle$. Now for any $y\in S_g$, the size of the conjugacy class of $y$ by $\langle g\rangle$ is $\frac{o(g)}{\lvert\langle g\rangle\cap Z(G)\rvert}$, which is even since $g\notin Z(G)$. Hence, $\lvert S_g\rvert$ is even. Now, from Equation (\ref{Eqn:Deg(a) is even iff, when Gis not 2gerted, Z(G)=e, ordr of a is 2}) we have that
\begin{equation}\label{Eqn:Deg(g) is even iff, when G is 2gerted}
    d_{\mathcal{DNI}(G)}(g)=r-\lvert S_g\rvert-1+\sum\limits_{y\in G\setminus\{e\}:\text{o}(y)\neq 2 \text{ and }g\sim y}\phi(\text{o}(y)).  
\end{equation}
So, $d_{\mathcal{DNI}(G)}(g)$ is even. This concludes the proof. 
\end{proof}


\section{Claw free}\label{Sec: Claw free}
We now turn our attention to the characterization of claw-free difference graphs. Our goal is to determine those finite groups, whether $2$-generated or non-$2$-generated, whose associated difference graphs are claw-free.

\begin{proposition}\label{Prop:complete-r-partite-claw-free}
Every complete $r$-partite graph is claw-free if and only if each part has size at most $2$.
\end{proposition}

\begin{proof}
Let $G$ be a complete $r$-partite graph with partition 
\[
V(G)=V_1 \cup V_2 \cup \cdots \cup V_r,
\]
where each $V_i$ is an independent set, and every vertex in $V_i$ is adjacent to every vertex in $V_j$ for $i \ne j$.

Suppose, to the contrary, that $G$ contains an induced claw, that is, an induced subgraph isomorphic to $K_{1,3}$. Let $v$ be the center of the claw and let $x_1, x_2, x_3$ be its three leaves. Since $x_1, x_2, x_3$ are pairwise non-adjacent, they must all lie in the same part, say $V_i$, because vertices from distinct parts are adjacent in a complete $r$-partite graph.

However, $v$ must be adjacent to each of $x_1, x_2, x_3$. This is possible only if $v \notin V_i$, since vertices within the same part are non-adjacent. Hence $v \in V_j$ for some $j \ne i$.

Now observe that if $|V_i| \ge 3$, then any vertex in another part $V_j$ is adjacent to \emph{all} vertices of $V_i$. Thus, the induced subgraph on $\{v,x_1,x_2,x_3\}$ contains exactly the edges required for a $K_{1,3}$. Therefore, $G$ contains a claw whenever some part has at least three vertices.
Hence, a complete $r$-partite graph is claw-free if and only if each part has size at most $2$.
\end{proof}
\subsection{Non $2$-generated groups}

\begin{corollary}\label{Cor:DNI-claw-free-G-p-grp}
Let $G$ be a non-$2$-generated $p$-group. Then $\mathcal{DNI}(G)$ is claw-free if and only if the following two conditions hold:
\begin{enumerate}
    \item[\textup{(i)}] $p \in \{2,3\}$;
    \item[\textup{(ii)}] $o(g) \in \{1,2,3\}$ for every $g \in G$.
\end{enumerate}
\end{corollary}

\begin{proof}
The proof follows from Lemma \ref{Lemma:p-grp-not-2-generated-structure-DNI-graph}, and Proposition \ref{Prop:complete-r-partite-claw-free}.   
\end{proof}
\begin{theorem}\label{Thm:DNI(G) is claw free iff G iso to S3, G is not 2-generated}
Let $G$ be a non $2$-generated group with $\pi(G)\geq 2.$ Then $\mathcal{DNI}(G)$ is claw-free if and only if $\text{o}(g)\in \{2, 3\},$ for all $g\in G\setminus\{e\}.$    
\end{theorem}
\begin{proof}
Since $\pi(G) \ge 2$, there exist distinct primes $p$ and $q$ dividing $|G|$. Choose $a, b \in G$ with $\mathrm{o}(a) = p$ and $\mathrm{o}(b) = q$. Assume, without loss of generality, that $p \ge 5$. Then $|G_a| \ge 4$, so we may choose distinct elements $a_1, a_2, a_3 \in G_a$. By Lemma~\ref{Lemma: Ga has no edge in DNI(G)}, no two of $a_1, a_2, a_3$ are adjacent in $\mathcal{DNI}(G)$. On the other hand, Lemma~\ref{Lemma:Two distinct prime elemnts edge connd in GI(G) complement; G is no 2-generated} shows that $b$ is adjacent to each of $a_1, a_2, a_3$. Hence the induced subgraph on $\{b, a_1, a_2, a_3\}$ is a claw in $\mathcal{DNI}(G)$.

Now suppose that $|G| = 2^{\alpha}3^{\beta}$, where $\alpha, \beta \ge 1$, and assume that $G$ contains an element $g$ with $\varphi(\mathrm{o}(g)) \ge 3$. Since $G$ is non-$2$-generated, $G$ is in particular non-cyclic. Consequently, there exists a prime $p \in \{2,3\}$ dividing $|G|$ such that $G$ has $r$ distinct cyclic subgroups of order $p$, with $r \ge 3$. Choose $x \in G$ with $\mathrm{o}(x) = p$ and $x \notin \langle g \rangle$. Proceeding as before, the vertices $x, g_1, g_2, g_3$, where $g_1, g_2, g_3 \in G_g$, induce a subgraph isomorphic to $K_{1,3}$.

\noindent
We now consider the case $|G| = 2^{\alpha}3^{\beta}$, where $\varphi(\mathrm{o}(g)) \le 2$ for every $g \in G \setminus \{e\}$; equivalently, the order of each element of $G$ lies in $\{2,3,4,6\}$.

Suppose first that $G$ has an element $h$ of order $6$, and let $h, h' \in G_h$ be distinct. By Lemma~\ref{Lemma: Ga has no edge in DNI(G)}, $h \nsim h'$ in $\mathcal{DNI}(G)$. Since $G$ is non-cyclic, there exists $y \in G \setminus \{e\}$ with $\mathrm{o}(y) \in \{2,3\}$ and $y \notin \langle h \rangle$. Moreover, $G$ contains an element $y'$ such that $\mathrm{o}(y') = \mathrm{o}(y)$, $\langle y \rangle \ne \langle y' \rangle$, and $y' \in \langle h \rangle$. Then $y$ is adjacent to each of $h$, $h'$, and $y'$ in $\mathcal{DNI}(G)$, while $y'$ is adjacent to neither $h$ nor $h'$, since $y' \in \langle h \rangle$. Consequently, the induced subgraph on $\{y, h, h', y'\}$ is isomorphic to $K_{1,3}$, and hence $\mathcal{DNI}(G)$ contains a claw. A similar argument shows that if $G$ contains an element of order $4$, then $\mathcal{DNI}(G)$ again contains an induced $K_{1,3}$.

\noindent
Now we prove that if $\mathrm{o}(g) \in \{2,3\}$ for all $g \in G \setminus \{e\}$, then $\mathcal{DNI}(G)$ is $K_{1,3}$-free. Suppose that the induced subgraph on the vertices $x',x_1,x_2,x_3$ is isomorphic to $K_{1,3}$, where $\deg(x')=3$ and $\deg(x_i)=1$ for all $i \in \{1,2,3\}$. If $\mathrm{o}(x')=2$, then $\mathrm{o}(x_1)\in\{2,3\}$. We have the following cases:

\medskip
\noindent
\textbf{Case I:} Suppose $\mathrm{o}(x') = 2$ and $\mathrm{o}(x_1) = 2$. Then $\langle x' \rangle \ne \langle x_1 \rangle$, and $\mathrm{o}(x_2) \in \{2,3\}$. In either case, since neither $\langle x_1 \rangle \subseteq \langle x_2 \rangle$ nor $\langle x_2 \rangle \subseteq \langle x_1 \rangle$, we have $x_1 \sim x_2$ in $\mathcal{DNI}(G)$, contradicting the assumption that the induced subgraph on $\{x', x_1, x_2, x_3\}$ is $K_{1,3}$.

\medskip
\noindent
\textbf{Case II:} Suppose $\mathrm{o}(x') = 2$ and $\mathrm{o}(x_1) = 3$.
\begin{itemize}
    \item If $\mathrm{o}(x_2) = 2$, then by Lemma~\ref{Lemma:Two distinct prime elemnts edge connd in GI(G) complement; G is no 2-generated}, $x_1 \sim x_2$, again contradicting that the induced subgraph is $K_{1,3}$.
    \item If $\mathrm{o}(x_2) = 3$ and $\langle x_2 \rangle = \langle x_1 \rangle$, then $x_1 \nsim x_2$; however, $x_3$ is then adjacent to both $x_1$ and $x_2$, contradicting the structure of $K_{1,3}$.
    \item If $\mathrm{o}(x_2) = 3$ and $\langle x_2 \rangle \ne \langle x_1 \rangle$, then $x_1 \sim x_2$ in $\mathcal{DNI}(G)$, again contradicting that the induced subgraph on $\{x', x_1, x_2, x_3\}$ is $K_{1,3}$.
\end{itemize}

\medskip
\noindent
If $\mathrm{o}(x') = 3$, a similar argument shows that the induced subgraph of $\mathcal{DNI}(G)$ on $\{x', x_1, x_2, x_3\}$ cannot be isomorphic to $K_{1,3}$. Consequently, $\mathcal{DNI}(G)$ is $K_{1,3}$-free whenever $\mathrm{o}(g) \in \{2,3\}$ for all $g \in G \setminus \{e\}$.
\end{proof}
\subsection{2-generated nilpotent groups}
\begin{proposition}\label{Prop: Non p groups have claw}
 Let $G$ and $H$ be $2$-generated groups with ${\rm gcd}(\lvert G\rvert,\lvert H\rvert)=1$. Then $\mathcal{DNI}(G\times H)$ is not claw-free.   
\end{proposition}
\begin{proof}
$G$ and $H$ are two finite $2$-generated $N_I$-groups with ${\rm gcd}(\lvert G\rvert,\lvert H\rvert)=1$. Let $p$ and $q$ be the smallest prime divisors of $\lvert G\rvert$ and $\lvert H\rvert$, respectively. If one among $p$ and $q$ is greater than $3$, assume that $p>3$. Let $x\in G$ such that $o(x)=p$, and $x_1,x_2,x_3\in G_x$. Then $x_1,x_2,x_3$ forms an independent set in $\bar{\mathcal{G}_I}^e$, so in $\mathcal{DNI}^e(G)$. Let $y\in H$ be a non-trivial element of $H$. Then Theorem \ref{Thm: direct product} implies that the subgraph induced by the vertices $(x_1,e),(x_2,e),(x_3,e),(e,y)$ in $\mathcal{DNI}(G\times H)$, is isomorphic to $K_{1,3}$.

Again, if $p=2, q=3$, let $\langle x_1,x_2\rangle=G$ and $\langle y_1,y_2\rangle=H$. Then $x_1,x_2,y_1,y_2$ are non identity elements, with $\langle x_1,x_1x_2\rangle=\langle x_2,x_1x_2\rangle=G$. So $x_1,x_2,x_1x_2$ forms an independent set in $\mathcal{DNI}(G)$. Now, since $\lvert H\rvert$ is odd, Lemma \ref{Lemma: Ga has no edge in DNI(G)} implies that $y_1\nsim y_1^2$ in $\mathcal{DNI}(H)$. Also, since $\langle y_1\rangle=\langle y_1^2\rangle$, we have $\langle y_1,y_2\rangle=\langle y_1^2,y_2\rangle=H$. So, $y_1,y_1^2,y_2$ forms an independent set in $\mathcal{DNI}(H)$. Then from Theorem \ref{Thm: direct product}, we have that $(x_1,y_1),(x_2,y_2), (x_1x_2,y_1^2)$ forms an independent set in $\mathcal{DNI}(G\times H)$. Now, Lemma \ref{Lemma:G_I complement diameter 2} implies that there exists nontrivial $h\in H$, such that $y_1\sim h\sim y_2$ in $\bar{\mathcal{G}_I}(G)$. Clearly, $y_1^2\sim h$ in $\bar{\mathcal{G}_I}(G)$, since $\langle y_1\rangle=\langle y_1^2\rangle$. So, Theorem \ref{Thm: direct product} implies that the subgraph induced by the vertices $(x_1,y_1),(x_2,y_2),(x_1x_2,y_1^2),(e,h)$ in $\mathcal{DNI}(G\times H)$, is isomorphic to $K_{1,3}$. So, $\mathcal{DNI}(G\times H)$ is not claw-free.
\end{proof}
As a consequence of the above result, we have the following.
\begin{corollary}\label{Cor: Non p groups have claw}
Let $G$ be a $2$-generated nilpotent group with $\pi(G)\geq 2$. Then $\mathcal{DNI}(G)$ is not claw-free.   
\end{corollary}

\begin{theorem}
Let $G$ be a $2$-generated nilpotent group. Then $\mathcal{DNI}(G)$ is claw-free if and only if either $G$ is an elementary abelian $p$-group or $G$ is isomorphic to the dihedral group $D_8$.
\end{theorem}
\begin{proof}
Let $G$ be a finite $2$-generated nilpotent $N_I$-group, which is claw free. Then corollary \ref{Cor: Non p groups have claw} implies that $G$ is a $p$-group, for some prime $p$. If $G$ is elementary abelian, then Lemma \ref{Lemma: Diameter abelian p group} implies that $\mathcal{DNI}(G)$ is claw free. Now, if $G$ is abelian with $\lvert G\rvert>p^2$, then let $G=\langle x,y\rangle$. Clearly $\langle x\rangle\cap\langle y\rangle=\{e\}$. Now, Corollary \ref{Cor: 2-generated p group} implies that both $o(x)>p$ and $o(y)>p$. So, $x,x^2,x^3$ forms an independent set in $\mathcal{DNI}(G)$ and $x\sim y^p$, $x^2\sim y^p$, $x^3\sim y^p$ in $\mathcal{DNI}(G)$. So $\mathcal{DNI}(G)$ is not claw-free. 

Let $G$ be non-abelian. Since $G$ is a $N_I$-group, Theorem \ref{Theorem: Gorenstein} implies that $G$ has a subgroup isomorphic to $C_p\times C_p$.   
\begin{itemize}
    \item[$\mathbf{\underline{p\geq 5:}}$] If $p\geq 5$, let $x,y$ be generators of $C_p\times C_p$ inside $G$. Clearly $\langle x,y\rangle\subsetneq G$. Then $x,x^2,x^3,y$ forms a claw in $G$.

    \item[\underline{$\mathbf{p=2:}$}] If $G$ has an element of order more than $4$, let that element be $x$. Then theorem \ref{Theorem: Gorenstein} implies that there exists $y\in G$ such that $\langle x\rangle\cap\langle y\rangle=\{e\}$. So, Lemma \ref{Lemma: langle x^p,y_ranglesubsetneq G} implies that the set $x^2,x^4,x^6,y$ forms a claw in $\mathcal{DNI}(G)$. If every element of $G$ has order $2$, then $G$ is abelian. So, the only remaining case is when $2\leq o(g)\leq 4$ for all nontrivial element of $G$, with at least one element having order $4$.
    
    Let $G$ be a group with ${\rm exp}(G)=4$, and $x\in G$ such that $o(x)=4$. Then by previous arguments, there exists $y\in G$ such that $o(y)=2$ and $\langle y\rangle\cap \langle x\rangle=\{e\}$. If $\langle x,y\rangle\subsetneq G$, then $x,x^2,x^3,y$ forms a claw in $\mathcal{DNI}(G)$. Let $\langle x,y\rangle=G$. If $G'$ has an element of order $4$, let that element be $c$. Then $\langle c\rangle\cap\langle x\rangle=\{e\}$ or $\langle c\rangle\cap\langle y\rangle=\{e\}$. Without loss of generality, we assume that $\langle c\rangle\cap\langle y\rangle=\{e\}$. Then Lemma \ref{Lemma: Cyclic abelianization} implies $\langle c,y\rangle\subsetneq G$. So, $c,c^2,c^3,y$ forms a claw in $\mathcal{DNI}(G)$. Let every element of $G'$ be of order $2$. If $G'$ is not cyclic, then there exists $z\in G'$ with $\langle z\rangle\cap \langle x \rangle=\{e\}$. Again by Lemma \ref{Lemma: Cyclic abelianization}, $\langle x,z\rangle\subsetneq G$. So, $x,x^2,x^3,z$ forms a claw in $\mathcal{DNI}(G)$. Let $G'=\langle z\rangle$. Then $o(z)=2$. If $\langle z\rangle\cap \langle x\rangle=\{e\}$, then again $x,x^2,x^3,z$ forms a claw in $\mathcal{DNI}(G)$. Now, if $\langle z\rangle\subseteq \langle x\rangle$, then $\langle x\rangle$ is normal in $G$. So, $G$ is the non-abelian semidirect product of $\langle x\rangle$ and $\langle y\rangle$, which is isomorphic to $D_8$. Also, it is easy to see that $\mathcal{DNI}(D_8)$ is claw-free. 

    \item[\underline{$\mathbf{p=3:}$}] If $G$ has an element of order more than $9$, then by similar arguments used in the previous case, it can be seen that $\mathcal{DNI}(G)$ is not claw free. Let $G=\langle x,y\rangle$ be a non-abelian group such that every element of $G$ has order $3$. Then $o(x)=o(y)=3$. Let $z\in G'$, be a nontrivial element. Then Lemma \ref{Lemma: Cyclic abelianization} implies that $\langle z\rangle\cap\langle x\rangle=\langle z\rangle\cap\langle y\rangle=\{e\}$, hence $x,x^2,y,z$ forms a claw in $\mathcal{DNI}(G)$. So, the only remaining case is when $3\leq o(g)\leq 9$ for all nontrivial element of $G$, with at least one element having order $9$. 
    
    Let $G$ be such a group with ${\rm exp}(G)=9$, and $x\in G$ such that $o(x)=9$. Then by previous arguments, there exists $y\in G$ such that $o(y)=3$ and $\langle y\rangle\cap \langle x\rangle=\{e\}$. If $\langle x,y\rangle\subsetneq G$, then $x,x^2,x^3,y$ forms a claw in $\mathcal{DNI}(G)$. If $\langle x,y\rangle=G$, then proceeding as the previous case, it can be seen that $\mathcal{DNI}(G)$ is claw free only if $G$ is the unique non-abelian semidirect product $C_9\rtimes C_3$. Now, Corollary \ref{Cor: metacyclic p index not NI group} implies that $C_9\rtimes C_3$ is not a $N_I$-group. 
    \end{itemize}
    Thus, the result follows.
\end{proof}

\section{Chordal}\label{Sec: Chordal}
This section is devoted to the study of the chordality of difference graphs of finite groups. The principal results are contained in three theorems. Theorems~\ref{Thm: Chordal non 2-generated p group structure} characterize the chordality of difference graphs for non-$2$-generated groups, while Theorem~\ref{Theorem: DNI chordality 2-generated nilpotent} addresses the case of $2$-generated nilpotent groups.

\begin{lemma}\label{Lemma:Induced cycle of length 4 in DNI(G), Z(G)=e, not 2 genrated, phio(g), phi(o(h))geq 2, gcd =1}
Let $g, h \in G \setminus \{e\}$ such that $\varphi(\mathrm{o}(g)) \ge 2$ and $\varphi(\mathrm{o}(h)) \ge 2$. Suppose that $\gcd(\mathrm{o}(g), \mathrm{o}(h)) = 1$. Then $\mathcal{DNI}(G)$ contains an induced cycle of length $4$.
\end{lemma}

\begin{proof}
Clearly, $g \sim h$ in $\mathcal{DNI}(G)$ by Lemma~\ref{Lemma:Two distinct prime elemnts edge connd in GI(G) complement; G is no 2-generated}. Hence, by Lemma~\ref{Lemma: x in Ga and y in Gb, if x sim y then any two elmnts in Ga and Gb edge connd in DNI(G)}, each vertex of $G_g$ is adjacent to every vertex of $G_h$.

Let $g_1, g_2 \in G_g$ and $h_1, h_2 \in G_h$. By Lemma~\ref{Lemma: Ga has no edge in DNI(G)}, we have $g_1 \nsim g_2$ and $h_1 \nsim h_2$. Therefore, 
\[
g_1 \sim h_1 \sim g_2 \sim h_2 \sim g_1
\]
forms an induced cycle of length $4$ in $\mathcal{DNI}(G)$.
\end{proof}

\begin{lemma}\label{Lemma: Odd p groups GI not chordal}
Let $G$ be a group of prime power order, for some odd prime. Then $\mathcal{DNI}(G)$ contains an induced four cycle.  
\end{lemma}
\begin{proof}
Let $G$ be a $p$-group for some odd prime $p$. Since $G$ is a $N_I$-group, $G$ has a subgroup of the form $C_p\times C_p$. Let $x, y$ be a set of generators for  $C_p\times C_p\subseteq G$. Then $\langle x\rangle\cap\langle y\rangle=\{e\}$, so $x\sim y$ in $\mathcal{DNI}(G)$. Since $p\geq 3$, we have $x\sim y\sim x^2\sim y^2\sim x$ with $x\nsim x^2$ and $y\nsim y^2$ in $\mathcal{DNI}(G)$. This concludes the proof.
\end{proof}

\subsection{Non 2-generated groups}
\begin{theorem}\label{Thm: Chordal non 2-generated p group structure}
Let $G$ be a non $2$-generated $p$-group. Then we have the following:
\begin{itemize}
    \item If $G$ is abelian, then $\mathcal{DNI}(G)$ is chordal if and only if $G\cong C_{2^s}\times C_2^t$, where $s,t\in\mathbb{N}$ with $t\geq 2$. 
    \item If $G$ is non-abelian, then $\mathcal{DNI}(G)$ is chordal if and only if $p=2$, and $G$ has a unique involution having square roots. 
\end{itemize}  
\end{theorem}
\begin{proof}
Let $G$ be a non $2$-generated $p$-group. Clearly, $G$ is neither cyclic nor generalized quaternion, and by Theorem \ref{Theorem: Gorenstein}, we have that $G$ has a proper subgroup $H$, such that $H\cong C_p\times C_p$. Choose two nontrivial elements $x,y\in H$, such that $y\notin \langle x\rangle$. If $p\geq 3$, then we see that $x\sim y\sim x^2\sim y^2\sim x$ forms a induced $4$-cycle in $\mathcal{DNI}(G)$. Thus, for $p\geq 3$ $\mathcal{DNI}(G)$ is not chordal. So, throughout the rest of the proof, we assume $p=2$.
\begin{enumerate}
    \item If $G$ is abelian, let $G\cong C_{2^{k_1}}\times C_{2^{k_2}}\times \cdots\times C_{2^{k_l}}$, where $k_1\geq k_2\geq \cdots\geq k_l$ are positive integers for some $l\in \mathbb{N}$. If $k_2\geq 2$, then $G$ has a proper subgroup isomorphic to $C_4\times C_4$, and by the above argument we see that $\mathcal{DNI}(G)$ contains a induced $4$ cycle, hence not chordal. Otherwise, let $k_2=1$. If $k_1=1$, then $G\cong C_2^t$ with $t\geq 3$. So, by Theorem \ref{Theorem: Completeness and nullness} we have $\mathcal{DNI}(G)$ is complete, hence chordal. Otherwise, let $x$ be the element of order $2$ in $C_{2^{k_1}}$. Let $y\in G$ such that $o(y)\geq 4$. Then it is easy to see that $y^2\in C_{2^{k_1}}$, i.e.~$y\in K_x$. Also, since $G$ is not two generated, any two distinct elements of order $2$ are adjacent in $\mathcal{DNI}(G)$. So, if we denote by $S$ the set of all elements of order $2$ in $G$, and by $T$ the set of all elements of order more than $2$ in $G$, then we have $S=K_x\setminus \{x\}$. Hence, we have $\mathcal{DNI}(G)\vert _S$ is a complete graph and $\mathcal{DNI}(G)\vert _T$ is a null graph. Now, let $\mathcal{C}$ be an induced cycle in $\mathcal{DNI}(G)$. Then either $\mathcal{C}$ is of length $3$, or $\mathcal{C}$ contains at most two vertices from $S$, and those two vertices must be adjacent in $\mathcal{C}$. So, if $\mathcal{C}$ has length more than $3$, then $\mathcal{C}$ contains two adjacent vertices from $T$, which is not possible. Hence $\mathcal{C}$ is a three cycle. Thus, $\mathcal{DNI}(G)$ is chordal. Thus, the result holds.

    \item Let $G$ be non-abelian. Clearly, $G$ is neither cyclic nor generalized quaternion. Hence, $G$ has at least three distinct involutions. Let $x,y\in G$ be two distinct involutions having square roots in $G$. Then there exists $x_1,y_1\in G$ such that $x_1^2=x$ and $y_1^2=y$. So we have a induced $4$ cycle $x_1\sim y_1\sim x_1^3\sim y_1^3\sim x_1$ in  $\mathcal{DNI}(G)$. So,  $\mathcal{DNI}(G)$ is not chordal. Now, if $G$ has a unique involution having square root, let that involution be $z$. Let $S$ denote the set of all elements of order $2$ in $G$, and $T$ denote the set of all elements of order more than $2$ in $G$. Then it is easy to see that $S=K_z\setminus \{z\}$, and $\mathcal{DNI}(G)\vert _S$ is a complete graph and $\mathcal{DNI}(G)\vert _T$ is a null graph. Now, from arguments similar to the previous case, the rest follows.
\end{enumerate}
This concludes the proof.
\end{proof}
\noindent\textbf{\underline{Remark:}} Observe that the set of all finite groups satisfying the conditions given in Proposition \ref{Thm: Chordal non 2-generated p group structure} is non-empty. For example, consider the group $Q_8\times C_2$ which is $3$ generated and satisfies the conditions of the proposition.

\begin{proposition}\label{Prop:DNI(G)-chordal-Imply(G-non-2-generated)}
Let $G$ be a non-$2$-generated group with $\pi(G) \ge 2$. If $\mathcal{DNI}(G)$ is chordal, then $G$ satisfies the following conditions:
\begin{enumerate}
    \item[\textup{(i)}] $|G| = 2^{k}p^{r}$ for some prime $p$ and integers $k, r \ge 1$;
    \item[\textup{(ii)}] $G$ has no element of order $4$; and
    \item[\textup{(iii)}] $G$ has a unique subgroup of order $p$.
\end{enumerate}
\end{proposition}

\begin{proof}
First, suppose that $2 \nmid |G|$. Since $\pi(G) \geq 2$ and $2 \nmid |G|$, there exist two elements $g$ and $h$ satisfying the conditions of Lemma~\ref{Lemma:Induced cycle of length 4 in DNI(G), Z(G)=e, not 2 genrated, phio(g), phi(o(h))geq 2, gcd =1}. Hence, $\mathcal{DNI}(G)$ contains an induced cycle of length $4$, which implies that $\mathcal{DNI}(G)$ is not chordal. This is a contradiction. Therefore, $\pi(G)=2$, and hence $|G|=2^{k}p^{r}$, where $p$ is an odd prime.

Next, assume that $G$ contains an element of order $4$. Since $p>2$ and $\varphi(p)\geq 2$, the group $G$ also contains an element of order $p$ satisfying the conditions of Lemma~\ref{Lemma:Induced cycle of length 4 in DNI(G), Z(G)=e, not 2 genrated, phio(g), phi(o(h))geq 2, gcd =1}. Thus, $\mathcal{DNI}(G)$ contains an induced cycle of length $4$, contradicting the assumption that $\mathcal{DNI}(G)$ is chordal. Therefore, $G$ has no element of order $4$.

Finally, suppose that $G$ has at least two distinct subgroups of order $p$. Let $H_{1}=\langle a\rangle$ and $H_{2}=\langle b\rangle$ be two distinct subgroups of order $p$. Let $a_{1},a_{2}$ be distinct generators of $H_{1}$ and $b_{1},b_{2}$ be distinct generators of $H_{2}$. Then
\[
a_{1}\sim b_{1}\sim a_{2}\sim b_{2}\sim a_{1}
\]
forms an induced cycle of length $4$ in $\mathcal{DNI}(G)$, a contradiction. Hence, $G$ has a unique subgroup of order $p$.
\end{proof}

\begin{lemma}\label{Lemma: not having Induced cycle of 4;Z(G)=e and has unique p order element}
Let $G$ be a finite group satisfying all the conditions stated in Proposition \ref{Prop:DNI(G)-chordal-Imply(G-non-2-generated)}. Then $\mathcal{DNI}(G)$ does not have any induced cycle of length $4.$ 
\end{lemma}
\begin{proof}
 Let 
\[
\mathcal{C} : a_1 \sim a_2 \sim a_3 \sim a_4 \sim a_1
\]
be an induced cycle of length $4$ in $\mathcal{DNI}(G)$. 
First, assume that $\mathrm{o}(a_1)=2$. Then either $p \mid \mathrm{o}(a_2)$ or $p \nmid \mathrm{o}(a_2)$. 
We consider the following cases.

\medskip
\noindent
\textbf{Case I:} Suppose that $p$ does not divide $\text{o}(a_2).$ That is, $\mathrm{o}(a_2)=2$. 
If $\mathrm{o}(a_3)=2$ or $p$, then in either case we have $a_1 \sim a_3$ in $\mathcal{DNI}(G)$, which contradicts the fact that $\mathcal{C}$ is an induced cycle. 
Therefore, $\mathrm{o}(a_3)=2p^{i}$ for some $i\geq 1$. 
Hence, $\langle a_1\rangle \subseteq \langle a_3\rangle$ (otherwise, $a_1 \sim a_3$ in $\mathcal{DNI}(G)$). 

Since $G$ has a unique subgroup of order $p$ and $a_3 \sim a_4$ in $\mathcal{DNI}(G)$, it follows that $p \nmid \mathrm{o}(a_4)$. 
Thus, $\mathrm{o}(a_4)=2$. 
Consequently, $a_2 \sim a_4$, which contradicts the assumption that 
\[
\mathcal{C} : a_1 \sim a_2 \sim a_3 \sim a_4 \sim a_1
\]
is an induced cycle of length $4$.

\medskip
\noindent
\textbf{Case II:} Suppose that $p \mid \mathrm{o}(a_2)$. 
Since $G$ has a unique subgroup of order $p$, it follows that $\mathrm{o}(a_3)=2$. 
Otherwise, if $p \mid \mathrm{o}(a_3)$, then $a_2 \sim a_3$ in the intersection power graph $\mathcal{G}_I(G)$, which implies that $a_2$ and $a_3$ are not adjacent in $\mathcal{DNI}(G)$. 
In that case, we would have $a_1 \sim a_3$ in $\mathcal{DNI}(G)$, contradicting the assumption that
\[
\mathcal{C} : a_1 \sim a_2 \sim a_3 \sim a_4 \sim a_1
\]
is an induced cycle of length $4$ in $\mathcal{DNI}(G)$.

\medskip
\noindent
If $p \mid \mathrm{o}(a_1)$, then by a similar argument we conclude that $\mathcal{DNI}(G)$ does not contain any induced cycle of length $4$. 
This completes the proof.
\end{proof}
Now, from the proof of Lemma~\ref{Lemma: not having Induced cycle of 4;Z(G)=e and has unique p order element}, we obtain the following.

\begin{remark}\label{Rem:MaxLengthInducedPath-2-GIsNot2generated}
Let $G$ be a finite group satisfying all the conditions stated in Proposition~\ref{Prop:DNI(G)-chordal-Imply(G-non-2-generated)}. 
Then the maximum length of an induced path in $\mathcal{DNI}(G)$ is $2$. Consequently, $\mathcal{DNI}(G)$ does not contain any induced cycle of length greater than $4$.
\end{remark}

\begin{theorem}\label{thm:DNI-chordal-non-2-generated}
Let $G$ be a non-$2$-generated group with $\pi(G) \geq 2$. Then $\mathcal{DNI}(G)$
is chordal if and only if $G$ satisfies the following conditions:
\begin{enumerate}
    \item[\textup{(i)}]   $|G| = 2^{k}p^{r}$ for some primes $2,p$ and integers $k,r \geq 1$;
    \item[\textup{(ii)}]  $G$ has no element of order $4$;
    \item[\textup{(iii)}] $G$ has a unique subgroup of order $p$.
\end{enumerate}
\end{theorem}
\begin{proof}
The proof follows from Proposition \ref{Prop:DNI(G)-chordal-Imply(G-non-2-generated)}, Lemma \ref{Lemma: not having Induced cycle of 4;Z(G)=e and has unique p order element}, and Remark \ref{Rem:MaxLengthInducedPath-2-GIsNot2generated}.    
\end{proof}

\subsection{2-generated nilpotent groups}
\begin{lemma}\label{Lemma: DNI chordality for 2-generated non p groups}
Let $G$ be a $2$-generated nilpotent group with $\pi(G)\geq 2$. Then $\mathcal{DNI}(G)$ is not chordal. 
\end{lemma}
\begin{proof}
Let $G$ be a finite, $2$-generated, nilpotent, $N_I$-group with $|\pi(G)|\geq 2$. Then $G\cong H\times K$, where $K$ is a Sylow $p$-subgroup of $G$ for some odd prime $p$, and ${\rm gcd}(\lvert H\rvert,\lvert K\rvert)=1$. So, from Theorem \ref{Thm: direct product} we have $\mathcal{DNI}^e(G)=(\mathcal{DNI}^e(H)\times\bar{\mathcal{G}_I}^e(K))\cup (\bar{\mathcal{G}_I}^e(H)\times\mathcal{DNI}^e(K))$. Now, Lemma \ref{Lemma: Odd p groups GI not chordal} implies that $\bar{\mathcal{G}_I}^e(K)$ contains an induced four cycle $x\sim y\sim x^2\sim y^2\sim x$. Then $(e,x)\sim (e,y)\sim (e,x^2)\sim (e,y^2)\sim (e,x)$ forms an induced four cycle in $\mathcal{DNI}(H\times K)$. Thus $\mathcal{DNI}(G)$ is not chordal.       
\end{proof}
\begin{lemma}\label{Lemma: DNI chordality for 2-generated odd p groups}
Let $G$ be a finite, $2$-generated, $N_I$-group of prime power order, for some odd prime $p$. Then $\mathcal{DNI}(G)$ is chordal if and only if $G\cong C_p\times C_p$.   
\end{lemma}
\begin{proof}
$G$ is a finite, $2$-generated, $N_I$-group of prime power order, for an odd prime $p$. If $G\cong C_p\times C_p$, then Lemma \ref{Lemma: Diameter abelian p group} says that $\mathcal{DNI}(G)$ is a null graph, hence chordal. Otherwise, from Corollary \ref{Cor: 2-generated p group} and theorem \ref{Theorem: Gorenstein} it can be seen that $G$ contains a proper subgroup of the form $C_p\times C_p$, and Lemma \ref{Lemma: Odd p groups GI not chordal} implies that $\bar{\mathcal{G}_I}(C_p\times C_p)$ contains an induced four cycle. Thus, $\mathcal{DNI}(G)$ contains an induced four cycle. So, $\mathcal{DNI}(G)$ is not chordal.
\end{proof}

\begin{lemma}\label{Lemma: Abelian 2 group chordal}
Let $G$ be a $2$-generated abelian group of order $2^k$, for some $k\in \mathbb{N}$. Then $\mathcal{DNI}(G)$ is chordal if and only if either $G\cong C_2\times C_2$, or $G\cong C_4\times C_4$.   
\end{lemma}
\begin{proof}
Let $G$ be a finite $2$-generated abelian $N_I$-group of order $2^k$, for some $k\in \mathbb{N}$. Then Corollary \ref{Cor: 2-generated p group} implies that $G\cong C_{2^r}\times C_{2^r}$, where $r=\frac{n}{2}$. If $r\geq 3$ then $G$ has a proper subgroup $H$ isomorphic to $C_4\times C_4$. Let $H=\langle x,y\rangle$. Then we have $o(x)=o(y)=4$, and $\langle x\rangle \cap \langle y\rangle=\{e\}$. So, $x\sim y\sim x^3\sim y^3\sim x$ forms an induced four cycle in $\mathcal{DNI}(G)$. Thus, $\mathcal{DNI}(G)$ is not chordal. Now, if $r=1$ then Theorem \ref{Theorem: Completeness and nullness} implies that $\mathcal{DNI}(G)$ is null graph, hence chordal. Moreover, if $r=2$, let $H$ denote the set of all elements of order $2$ in $G$, and $K$ denote the set of all elements of order $4$ in $G$. Then $G=\{e\}\cup H\cup K$. Also, it is easy to see that no two distinct elements of $H$ generate $G$, and for any two distinct elements $k_1,k_2\in \langle K\rangle$ with $\langle k_1\rangle \cap \langle k_2\rangle={e}$, we have that $\langle k_1,k_2\rangle=G$. So, $\mathcal{DNI}(G)\vert_H$ is a complete graph, and $\mathcal{DNI}(G)\vert_K$ is a null graph. Now, let $\mathcal{C}$ be an induced cycle in $\mathcal{DNI}(G)$. Then either $\mathcal{C}$ is of length $3$, or $\mathcal{C}$ contains at most two vertices from $H$, and those two vertices must be adjacent in $\mathcal{C}$. So, if $\mathcal{C}$ has length more than $3$, then $\mathcal{C}$ contains two adjacent vertices from $K$, which is not possible. Hence $\mathcal{C}$ is a three cycle. Thus, $\mathcal{DNI}(G)$ is chordal. This concludes the proof.
\end{proof}

\begin{lemma}\label{Lemma: D_2^n is chordal if and only if $n=3$}
Let $D_{2^n}$ denote the dihedral group of order $2^n$. Then $\mathcal{DNI}(D_{2^n})$ is chordal if and only if $n=3$.    
\end{lemma}
\begin{proof}
Let $n\geq 4$. We have $D_{2^n}:=\langle a,b~\vert~a^{2^{n-1}}=b^2=e,~b^{-1}ab=a^{-1}\rangle$. Clearly, $\langle a\rangle\cap\langle b\rangle=\{e\}$. Consider the set $a^2,b,a^6,ab$. Here $\langle b,ab\rangle=D_{2^n}$, so $b\nsim ab$ and $a^2\nsim a^6$ in $\mathcal{DNI}(D_{2^n})$. Again, we have $\langle a^2\rangle=\langle a^6\rangle$. So, $\langle a^2,b\rangle=\langle a^6,b\rangle=\langle a^2,b~\vert~(a^2)^{2^{n-2}}=b^2=e,~b^{-1}a^2b=a^{-2}\rangle\cong D_{2^{n-1}}$. Hence, $a^2\sim b\sim a^6$ in $\mathcal{DNI}(D_{2^n})$. Similarly, it can be shown that $a^2\sim ab\sim a^6$ in $\mathcal{DNI}(D_{2^n})$. Thus, we have an induced four cycle $a^2\sim b\sim a^6\sim ab\sim a^2$ in $\mathcal{DNI}(D_{2^n})$. So, $\mathcal{DNI}(D_{2^n})$ is not chordal. 

Now, if $n=3$ then the graph $\mathcal{DNI}(D_{8})$ is the following:

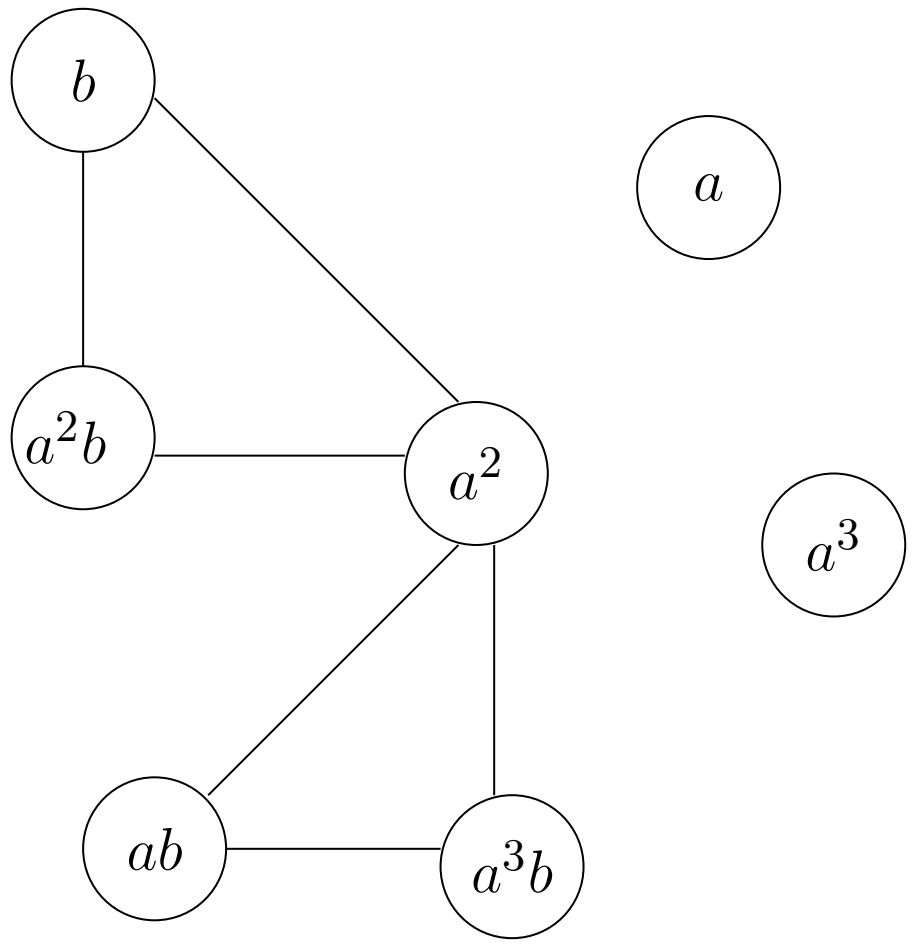
\begin{figure}[!ht]
\centering
\resizebox{0.35\textwidth}{!}{%
\begin{circuitikz}
\tikzstyle{every node}=[font=\fontsize{13.1pt}{17.0pt}\selectfont]
\draw  (3,10.5) circle (0.5cm);
\draw  (3,8) circle (0.5cm);
\draw  (5.75,7.75) circle (0.5cm);
\draw  (7.375,9.75) circle (0.5cm);
\draw  (8.25,7.25) circle (0.5cm);
\draw  (3.5,5.125) circle (0.5cm);
\draw  (6,5) circle (0.5cm);
\draw (3,10) to[short] (3,8.5);
\draw (3.5,10.375) to[short] (5.625,8.25);
\draw (3.5,7.875) to[short] (5.25,7.875);
\draw (3.875,5.5) to[short] (5.625,7.25);
\draw (5.875,7.25) to[short] (5.875,5.5);
\draw (4,5.125) to[short] (5.5,5.125);
\node [font=\fontsize{13.1pt}{17.0pt}\selectfont, inner xsep=0.080cm, inner ysep=0.085cm, rounded corners=0.020cm] at (2.875,8) {$a^2b$};
\node [font=\fontsize{13.1pt}{17.0pt}\selectfont, inner xsep=0.080cm, inner ysep=0.085cm, rounded corners=0.020cm] at (3,10.5) {$b$};
\node [font=\fontsize{13.1pt}{17.0pt}\selectfont, inner xsep=0.080cm, inner ysep=0.085cm, rounded corners=0.020cm] at (5.75,7.75) {$a^2$};
\node [font=\fontsize{13.1pt}{17.0pt}\selectfont, inner xsep=0.080cm, inner ysep=0.085cm, rounded corners=0.020cm] at (3.5,5.125) {$ab$};
\node [font=\fontsize{13.1pt}{17.0pt}\selectfont, inner xsep=0.080cm, inner ysep=0.085cm, rounded corners=0.020cm] at (6,5) {$a^3b$};
\node [font=\fontsize{13.1pt}{17.0pt}\selectfont, inner xsep=0.080cm, inner ysep=0.085cm, rounded corners=0.020cm] at (7.375,9.75) {$a$};
\node [font=\fontsize{13.1pt}{17.0pt}\selectfont, inner xsep=0.080cm, inner ysep=0.085cm, rounded corners=0.020cm] at (8.25,7.25) {$a^3$};
\end{circuitikz}
}%
\caption{The graph $\mathcal{DNI}(D_8)$.}
\label{fig:The graph {DNI}(D_8)}
\end{figure}
\noindent So, $\mathcal{DNI}(D_8)$ is chordal. This concludes the proof.
\end{proof}

\begin{lemma}\label{Lemma: Semidirect product C_2 chordal}
The group $(C_4\times C_2)\rtimes C_2$ is a $N_I$-group and the graph $\mathcal{DNI}((C_4\times C_2)\rtimes C_2)$ is chordal.   
\end{lemma}
\begin{proof}
Consider the group $(C_4\times C_2)\rtimes C_2=\langle a,b,c~\vert~a^2=b^4=c^2=[b,c]=e,~a^{-1}ba=bc,~a^{-1}ca=c\rangle$. It can be seen that the set $H:=\{ a,b^2,c,b^2c,ab^2,ac,ab^2c\}$ denotes the set of all elements of order $2$ in $(C_4\times C_2)\rtimes C_2$, and the set $K:=\{b, bc, b^{-1},b^{-1}c,ab,abc,ab^{-1},ab^{-1}c\}$ denotes the set of all elements of order $4$ in $(C_4\times C_2)\rtimes C_2$. Further, we have $\{e\}\cup H\cup K=(C_4\times C_2)\rtimes C_2$. Let $K_1:=\{b, bc\}$ and $K_2:=\{ab,abc\}$. Define $K_1^{-1}:=\{b^{-1},(bc)^{-1}\}$ and $K_2^{-1}=\{(ab)^{-1},(abc)^{-1}\}$. Then $K=K_1\cup K_1^{-1}\cup K_2\cup K_2^{-1}$. Again, we have $x^2=b^2$ for all $x\in K_1\cup K_1^{-1}$, and $y^2=b^2c$ for all $y\in K_2\cup K_2^{-1}$. In addition, it can be checked that the subgroups $\langle K_1,K_1^{-1}\rangle$ and $\langle K_2,K_2^{-1}\rangle$ are both abelian. Hence, we have $\langle K_1,K_1^{-1}\rangle\subsetneq (C_4\times C_2)\rtimes C_2$ and $\langle K_2,K_2^{-1}\rangle\subsetneq (C_4\times C_2)\rtimes C_2$. Since $(C_4\times C_2)\rtimes C_2$ is non-cyclic, we have shown that for any two distinct elements $x,y\in (C_4\times C_2)\rtimes C_2$ with $\langle x\rangle\cap \langle y\rangle\neq \{e\}$, we have $\langle x,y\rangle\subsetneq(C_4\times C_2)\rtimes C_2$. Thus $(C_4\times C_2)\rtimes C_2$ is a $N_I$-group.

Since, $(C_4\times C_2)\rtimes C_2$ is not isomorphic to the dihedral group, no two elements of order $2$ generate $(C_4\times C_2)\rtimes C_2$. So, the induced subgraph $\mathcal{DNI}((C_4\times C_2)\rtimes C_2)\vert_{H}\cong K_7$. Also, we have $\langle b,ab\rangle=\langle a,b\rangle=(C_4\times C_2)\rtimes C_2$, and $\langle b,abc\rangle=\langle b,abca^2\rangle=\langle b,ba\rangle=(C_4\times C_2)\rtimes C_2$. So, $b$ is not adjacent to any other vertex of $\mathcal{DNI}((C_4\times C_2)\rtimes C_2)\vert_K$. Again, since $bc=a^{-1}ba$, from Observation \ref{Obs:DNI(G) is invariant under any automorphism} we can conclude that $bc$ is also an isolated vertex of $\mathcal{DNI}((C_4\times C_2)\rtimes C_2)\vert_K$. Thus, we have $\mathcal{DNI}((C_4\times C_2)\rtimes C_2)\vert_K$ is the null graph. Now, using similar arguments as in the proof of Lemma \ref{Lemma: Abelian 2 group chordal}, it can be seen that $\mathcal{DNI}((C_4\times C_2)\rtimes C_2)$ is chordal. This concludes the proof. 
\end{proof}

\begin{lemma}\label{Lemma: Semidirect product C_4 chordal}
The group $(C_4\times C_2)\rtimes C_4$ is a $N_I$-group and the graph $\mathcal{DNI}((C_4\times C_2)\rtimes C_4)$ is chordal.   
\end{lemma}
\begin{proof}
Consider the group $(C_4\times C_2)\rtimes C_4=\langle a,b,c~\vert~a^4=b^4=c^2=[b,c]=e,~a^{-1}ba=bc,~a^{-1}ca=c\rangle$. Then $(C_4\times C_2)\rtimes C_4=\langle a,b\rangle$. Also, it can be seen that the set $H:=\{a^2,b^2, c,b^2c,a^2c,a^2b^2,a^2b^2c\}$ denotes the set of all elements of order $2$ in $(C_4\times C_2)\rtimes C_2$. Since, $(C_4\times C_2)\rtimes C_4$ is not isomorphic to the dihedral group, we have that the induced subgraph $\mathcal{DNI}((C_4\times C_2)\rtimes C_4)\vert_H$ is a complete graph. Define \[K:=\{a, ac, b, bc, ab, abc, a^2b, ab^2, a^2bc, ab^2c, a^3b, a^3bc \}.\] It can be checked that $K\cap K^{-1}=\emptyset$, and $K\cup K^{-1}$ is the set of all elements of order $4$ in $(C_4\times C_2)\rtimes C_4$. In addition, we have $\{e\}\cup H\cup K\cup K^{-1}=(C_4\times C_2)\rtimes C_4$. Now, define the sets $K_1=\{a, ab^2, ac, ab^2c\}$, $K_2=\{b, bc, a^2b, a^2bc\}$, and $K_3=\{ab, abc, a^3b, a^3bc\}$. Then we have $K=K_1\cup K_2\cup K_3$. In addition, the following statements are easy to check:
\begin{itemize}
    \item The set $K$ is closed under conjugation in $(C_4\times C_2)\rtimes C_4$.
    \item We have $x^2=a^2$ for all $x\in K_1$, $y^2=b^2$ for all $y\in K_2$, and $z^2=a^2b^2c$ for all $z\in K_3$.
    \item The subgroups $\langle K_1\rangle$, $\langle K_2\rangle$ and $\langle K_3\rangle$ are abelian.
\end{itemize}
Since $(C_4\times C_2)\rtimes C_4$ is non-cyclic, we have shown that for any two distinct elements $x,y\in (C_4\times C_2)\rtimes C_4$ with $\langle x\rangle\cap \langle y\rangle\neq \{e\}$, we have $\langle x,y\rangle\subsetneq(C_4\times C_2)\rtimes C_4$. Thus, $(C_4\times C_2)\rtimes C_4$ is a $N_I$-group. 

We will now show that the induced subgraph $\mathcal{DNI}((C_4\times C_2)\rtimes C_4)\vert_K$ is a null graph. To show that, we first observe that the induced subgraphs $\mathcal{DNI}((C_4\times C_2)\rtimes C_4)\vert_{K_1}$, $\mathcal{DNI}((C_4\times C_2)\rtimes C_4)\vert_{K_2}$, and $\mathcal{DNI}((C_4\times C_2)\rtimes C_4)\vert_{K_3}$ are null graphs. Also, we have the following:
\begin{enumerate}
\item \begin{itemize}
      \item $\langle a,ab \rangle=\langle a,a^2b \rangle=\langle a,a^3b \rangle=\langle a,b \rangle=(C_4\times C_2)\rtimes C_4$.
      \item $\langle a,abc \rangle=\langle a,a^2bc \rangle=\langle a,a^3bc \rangle=\langle a,bc \rangle=\langle a,a^{-1}bca \rangle=\langle a,b \rangle=(C_4\times C_2)\rtimes C_4$.
      \end{itemize} 
      So, $a$ is an isolated vertex in $\mathcal{DNI}((C_4\times C_2)\rtimes C_4)\vert_K$. Also, we have $ac=b^{-1}ab$. So, by Observation \ref{Obs:DNI(G) is invariant under any automorphism} we see that $ac$ is also an isolated vertex in $\mathcal{DNI}((C_4\times C_2)\rtimes C_4)\vert_K$.

\item \begin{itemize}
      \item $\langle ab^2,bc \rangle=\langle ab^2,(ab^2)^{-1}bc(ab^2)\rangle=\langle ab^2,b \rangle=\langle a,b \rangle=(C_4\times C_2)\rtimes C_4$.
      \item $\langle ab^2,abc \rangle=\langle ab^2,(ab^2)^{-1}abc(ab^2)\rangle=\langle ab^2,ab \rangle=\langle a,b \rangle=(C_4\times C_2)\rtimes C_4$.
      \item $\langle ab^2,a^2b \rangle=\langle ab^2,(ab^2)^{-1}a^2b(ab^2)\rangle=\langle ab^2,a^2bc \rangle=\langle ab^2,(ab^2)^2a^2bc \rangle=\langle ab^2,bc \rangle=(C_4\times C_2)\rtimes C_4$.
      \item $\langle ab^2,a^3b \rangle=\langle ab^2,(ab^2)^{-1}a^3b(ab^2)\rangle=\langle ab^2,a^3bc \rangle=\langle ab^2,(ab^2)^2a^3bc \rangle=\langle ab^2,abc \rangle=(C_4\times C_2)\rtimes C_4$.
      \end{itemize} 
      So, $ab^2$ is an isolated vertex in $\mathcal{DNI}((C_4\times C_2)\rtimes C_4)\vert_K$. Also, we have $ab^2c=b^{-1}ab^2b$. So, by Observation \ref{Obs:DNI(G) is invariant under any automorphism} we see that $ab^2c$ is also an isolated vertex in $\mathcal{DNI}((C_4\times C_2)\rtimes C_4)\vert_K$.

\item \begin{itemize}
      \item $\langle b, abc \rangle=\langle b, b^{-1}abcb \rangle=\langle b, ab \rangle=\langle b, a^3b \rangle=\langle b, a \rangle=(C_4\times C_2)\rtimes C_4$.
      \item $\langle b, a^3bc \rangle=\langle b, b^{-1}a^3bcb \rangle=\langle b, a^3b \rangle=(C_4\times C_2)\rtimes C_4$.
      \end{itemize}
      So, $b$ is an isolated vertex in $\mathcal{DNI}((C_4\times C_2)\rtimes C_4)\vert_K$. Also, we have $bc=a^{-1}ba$. So, by Observation \ref{Obs:DNI(G) is invariant under any automorphism} we see that $bc$ is also an isolated vertex in $\mathcal{DNI}((C_4\times C_2)\rtimes C_4)\vert_K$.
      
\item \begin{itemize}
      \item $\langle a^2b, abc \rangle=\langle a^2b, (a^2b)^{-1}abc(a^2b) \rangle=\langle a^2b, ab \rangle=\langle a,ab \rangle=(C_4\times C_2)\rtimes C_4$.
      \item $\langle a^2b, a^3bc \rangle=\langle a^2b, (a^2b)^{-1}a^3bc(a^2b) \rangle=\langle a^2b, a^3b \rangle=\langle a,a^2b \rangle=(C_4\times C_2)\rtimes C_4$.
      \end{itemize}
      So, $a^2b$ is an isolated vertex in $\mathcal{DNI}((C_4\times C_2)\rtimes C_4)\vert_K$. Also, we have $a^2bc=a^{-1}a^2ba$. So, by Observation \ref{Obs:DNI(G) is invariant under any automorphism} we see that $a^2bc$ is also an isolated vertex in $\mathcal{DNI}((C_4\times C_2)\rtimes C_4)\vert_K$.
\end{enumerate}
Thus, we have that the induced subgraph $\mathcal{DNI}((C_4\times C_2)\rtimes C_4)\vert_K$ is a null graph. Hence, from Lemma \ref{Lemma: x in Ga and y in Gb, if x sim y then any two elmnts in Ga and Gb edge connd in DNI(G)} we can conclude that $\mathcal{DNI}((C_4\times C_2)\rtimes C_4)\vert_{K\cup K^{-1}}$ is a null graph. Now, using similar arguments as in the proof of Lemma \ref{Lemma: Abelian 2 group chordal}, it can be seen that $\mathcal{DNI}((C_4\times C_2)\rtimes C_4)$ is chordal. This concludes the proof.
\end{proof}

\begin{theorem}\label{Theorem: DNI chordality 2-generated nilpotent}
Let $G$ be a $2$-generated nilpotent group. Then $\mathcal{DNI}(G)$ is chordal if and only if $G$ is isomorphic to one of the following groups:\\ $C_p\times C_p$ with $p$ being a prime, $C_4\times C_4$, $D_8$, $(C_4\times C_2)\rtimes C_2$, $(C_4\times C_2)\rtimes C_4$.   
\end{theorem}
\begin{proof}
Let $G$ be a finite, $2$-generated, nilpotent, $N_I$-group, such that $\mathcal{DNI}(G)$ is chordal. Then Lemma \ref{Lemma: DNI chordality for 2-generated non p groups} implies that $G$ is a $p$-group, for some prime $p$. Further, if $p$ is odd then Lemma \ref{Lemma: DNI chordality for 2-generated odd p groups} implies that $G\cong C_p\times C_p$. So the only situation we are left with is $p=2$.

If $G$ is an abelian $2$-group, then Lemma \ref{Lemma: Abelian 2 group chordal} implies that $G\cong C_2\times C_2$ or $G\cong C_4\times C_4$. In addition, if $G\cong D_{2^n}$, then Lemma \ref{Lemma: D_2^n is chordal if and only if $n=3$} implies that $G\cong D_8$. So, throughout the rest of the proof, we assume $G$ to be a non-abelian $2$ group such that $G\ncong D_{2^n}$. Let $G=[a,b]$. Then $\langle a\rangle\cap\langle b\rangle=\{e\}$. If $o(a)=o(b)=2$, then $G=\langle a,b\rangle\cong D_{2n}$, which is contradictory to our assumption. So, at least one of $a$ and $b$ has order more than $2$. Let $o(a)=2$ and $o(b)\geq 8$. If $o(ab)=2$, then we have $abab=e$, which implies $a^{-1}ba=b^{-1}$, i.e.~$G\cong D_{2n}$ which is contradictory. Hence, $o(ab)\geq 4$.  Now, $G=\langle a,b\rangle=\langle ab,b\rangle$. So, $\langle ab\rangle\cap\langle b\rangle=\{e\}$. Then Lemma \ref{Lemma: langle x^p,y_ranglesubsetneq G} implies $\langle b^2,ab\rangle\subsetneq G$. So, $b^2\sim ab\sim b^6\sim (ab)^2\sim b^2$ forms an induced four cycle in $\mathcal{DNI}(G)$. Thus, $\mathcal{DNI}(G)$ is not chordal. Also, if $o(a)\geq 4$ and $o(b)\geq 8$, then using the above arguments, it can be seen that $\mathcal{DNI}(G)$ is not chordal. If $o(a)=2$ and $o(b)=4$, then $o(ab)\geq 4$ since $G\ncong D_{2n}$. In addition, if $o(ab)\geq 8$ then the fact that $G=\langle b, ab\rangle$ and previous arguments yield $\mathcal{DNI}(G)$ is not chordal. So, $o(ab)=4$. Hence, the only scenario we are left with is $G$ is generated by two elements of order $4$.

Let $o(a)=o(b)=4$. We first show that $G$ is of nilpotency class $2$. Take $b'=a^{-1}ba$. Then $o(b')=o(b)=4$, and $b'=b[b,a]$, i.e. $\langle b,b'\rangle=\langle b,[b,a]\rangle$. Now, from Lemma \ref{Lemma: Cyclic abelianization} we have $\langle b,b'\rangle=\langle b,[b,a]\rangle\subsetneq G$. So, if $b^2\neq b'^2$, then $\langle b\rangle\cap\langle b'\rangle=\{e\}$, hence $b\sim b'\sim b^3\sim b'^3\sim b$ forms an induced $4$ cycle in $\mathcal{DNI}(G)$, which is contradictory. So, we have $b'^2=b^2$, i.e. $a^{-1}b^2a=b^2$. Hence $b^2\in Z(G)$. Similarly, it can be shown that $a^2\in Z(G)$. Now, since $\langle a\rangle\cap\langle b\rangle=\{e\}$, we have $\langle [a,b]\rangle\cap\langle a\rangle=\{e\}$ or $\langle [a,b]\rangle\cap\langle b\rangle=\{e\}$. Without loss of generality, we assume $\langle [a,b]\rangle\cap\langle a\rangle=\{e\}$. If $o([a,b])\geq 4$, then $a\sim[a,b]\sim a^3\sim [a,b]^3\sim a$ is an induced $4$ cycle in $\mathcal{DNI}(G)$, which is contradictory. So, we have $o([a,b])=2$. Then $a^{-1}[a,b]a=a^{-2}b^{-1}aba=b^{-1}a^{-1}ba=[b,a]=[a,b]$, and $b^{-1}[a,b]b=b^{-1}a^{-1}b^{-1}ab^2=b^{-1}a^{-1}ba=[b,a]=[a,b]$. So $[a,b]\in Z(G)$, hence $G$ is of nilpotency class $2$.

Define $c:=[a,b]$. Then $o(c)=2$, and it is easy to see that $b^{-1}ab=ac$ and $a^{-1}ba=bc$. If $c=[a,b]=b^2$, then $a^{-1}b^{-1}a=b$. So, $(ab)^2=a^2a^{-1}bab=a^2b^{-1}b=a^2$, i.e.~$\langle a\rangle\cap\langle ab\rangle\neq \{e\}$ with $G=\langle a,ab\rangle$, which is contradictory to the fact that $G$ is a $N_I$-group. So $c\neq b^2$, i.e.~$\langle b,c\rangle\cong C_4\times C_2$. Clearly $a\notin \langle b,c\rangle$. If $a^2\in \langle b,c\rangle$, then we have that $a^2\neq b^2$. Also, if $a^2=c$, then $b^{-1}ab=a^{-1}$, which implies $(ab)^2=ab^2b^{-1}ab=b^2$, i.e.~$\langle b\rangle\cap\langle ab\rangle\neq \{e\}$ with $G=\langle b,ab\rangle$, which is contradictory to the fact that $G$ is a $N_I$-group. So, we have $a^2=b^2c$. Then $(ab)^2=ab^2b^{-1}ab=b^2a^2c=b^4c^2=e$. Define $a_1=ab$. Then we have $G=\langle a_1,b\rangle$ and the only relations between $a_1$ and $b$ are, $a_1^2=b^4=c^2=[b,c]=e,~a_1^{-1}ba_1=bc$ and $a_1^{-1}ca_1=c$. So, $G\cong(C_4\times C_2)\rtimes C_2$. Otherwise, if $a^2\notin \langle b,c\rangle$, we have $G=\langle a,b\rangle$ and the only relations between $a$ and $b$ are, $a^4=b^4=c^2=[b,c]=e,~a^{-1}ba=bc$ and $a^{-1}ca=c$. So, $G\cong(C_4\times C_2)\rtimes C_4$. Now, from Lemmas \ref{Lemma: Semidirect product C_2 chordal} and \ref{Lemma: Semidirect product C_4 chordal} we have that $\mathcal{DNI}((C_4\times C_2)\rtimes C_2)$ and $\mathcal{DNI}((C_4\times C_2)\rtimes C_4)$ are chordal. Thus, the result holds.
\end{proof}

\section{Planar}\label{Sec: Planar}
In this section, we characterize all finite non-$2$-generated groups and finite nilpotent $2$-generated groups whose difference graphs are planar. The main results of this section are Theorems~\ref{Thm:DNI_G_Planar iff when G is not 2 genreted, Z(G)=e} and~\ref{Theorem: DNI planarity 2-generated nilpotent}, which provide complete characterizations of the planar difference graphs in these two classes of groups.

Throughout this section, for a prime number $p$, we denote by $n_p$ the number of distinct subgroups of order $p$ in $G$.

\subsection{Non 2-generated groups}
\begin{proposition}\label{Prop: Planar non 2-generated p group structure}
Let $G$ be a non $2$-generated $p$-group. Then we have the following:
\begin{itemize}
    \item If $G$ is abelian, then $\mathcal{DNI}(G)$ is not planar.
    \item If $G$ is non-abelian, then $\mathcal{DNI}(G)$ is planar if and only if $p=2$, and $G$ has exactly $3$ involutions with a unique involution having square root. 
\end{itemize}
\end{proposition}
\begin{proof}
    Let $G$ be a non $2$-generated $p$-group. So, it is an $\rm{NI}$-group. Clearly, $G$ is neither cyclic nor generalized quaternion. Then by Theorem \ref{Theorem: Gorenstein}, we have that $G$ has a proper subgroup $H$, such that $H\cong C_p\times C_p$. If $p=3$, we see that the induced subgraph $\mathcal{DNI}(G)\vert_H$ is isomorphic to $K_{2,2,2,2}$. Then from Observation \ref{K_2222 not planar} we see that $\mathcal{DNI}(G)\vert_H$ is not planar, hence $\mathcal{DNI}(G)$ is not planar. If $p\geq 3$, choose two nontrivial elements $x,y\in H$, such that $y\notin \langle x\rangle$. Then the subgraph induced by $\mathcal{DNI}(G)$ on the set $x,x^2,x^3,y,y^2,y^3$ is isomorphic to $K_{3,3}$, hence $\mathcal{DNI}(G)$ is not planar. Hence, throughout the rest of the proof, we assume $p=2$.

    For $p=2$, Lemma \ref{Lemma:number of elements of order p in a p group} implies that if $G$ has $r$ involutions, then $r$ is odd and $r\geq 3$. If $r\geq 5$, then it is easy to see that any five distinct involutions form a induced $K_5$ in $\mathcal{DNI}(G)$. So, throughout the rest of the proof, we assume that $r=3$.
\begin{enumerate}
    \item If $G$ is abelian, let $G\cong C_{2^{k_1}}\times C_{2^{k_2}}\times \cdots\times C_{2^{k_l}}$, where $k_1\geq k_2\geq \cdots\geq k_l$ are positive integers for some $l\in \mathbb{N}$. If $l\geq 3$, then $G$ has a subgroup isomorphic to $C_2\times C_2\times C_2$, which contradicts to our assumption $r=3$. So, we have $l=2$, i.e. $G$ is $2$-generated, which is contradictory. So we conclude that if $G$ is abelian and non $2$-generated, then $\mathcal{DNI}(G)$ is not planar. 

    \item Let $G$ be non-abelian. Let $a,b\in G$ be two distinct involutions having square roots in $G$. Then there exists $a_1,b_1\in G$ such that $a_1^2=a$ and $b_1^2=b$. So the subgraph induced by $\mathcal{DNI}(G)$ on the set $a,a^2,a^3,b,b^2,b^3$ is isomorphic to $K_{3,3}$, hence $\mathcal{DNI}(G)$ is not planar. Now, if $G$ has a unique involution having square root, let that involution be $z_1$. Let $z_2,z_3$ be the other two involutions in $G$ and let $K_{z_1}=\{w_1,w_2,\cdots, w_m, z_1\}$. Then it is easy to see that $K_{z_1}\cup \{z_2,z_3\}$ forms a partition of $G\setminus \{e\}$. So, Figure \ref{fig:The planar embedding} is a planar embedding of $\mathcal{DNI}(G)$.  
   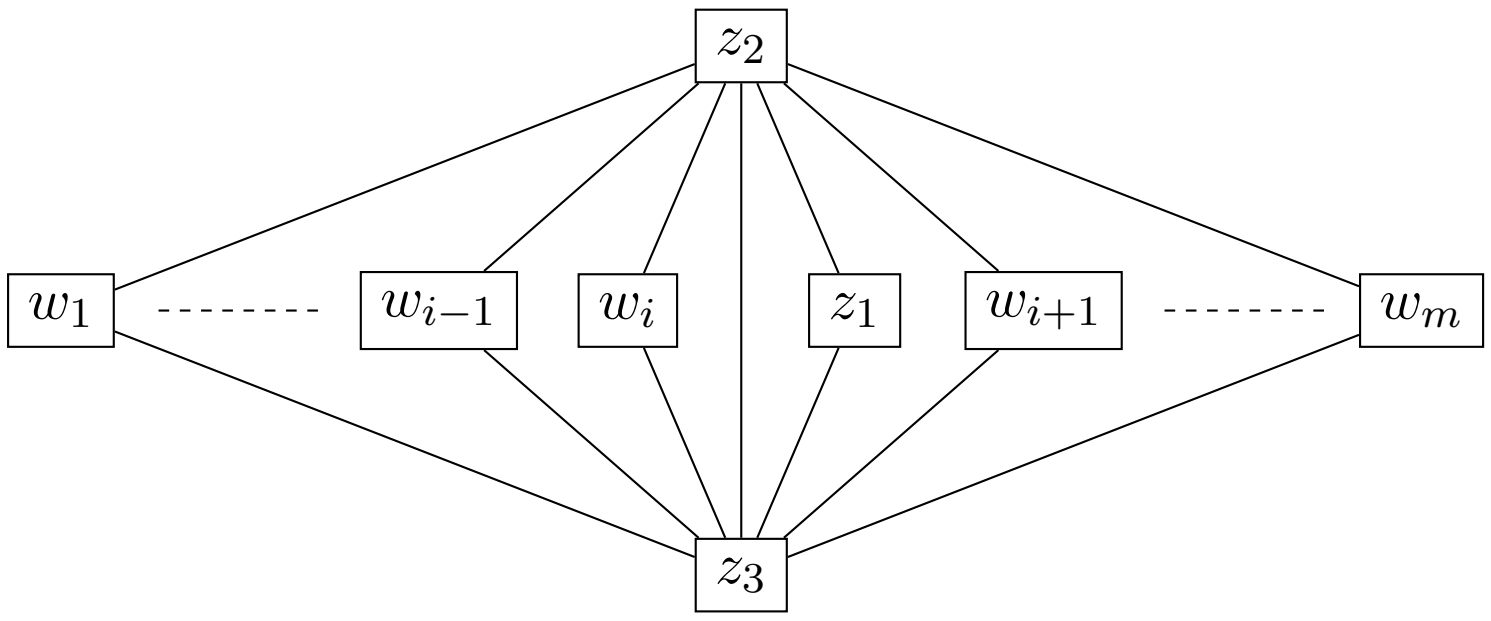
\begin{figure} 
    \begin{center}
\begin{tikzpicture}
  [scale=.5,auto=center,every node/.style={draw,fill=white}] 
    
  \node (a1) at (0,3.5) {$z_2$};  
  \node (a2) at (0,-3.5) {$z_3$};
  \node (a3) at (1.5,0) {$z_1$};  
  \node (a4) at (-1.5,0) {$w_i$};  
  \node (a5) at (4,0) {$w_{i+1}$};
  \node (a6) at (9,0) {$w_m$};  
  \node (a7) at (-4,0) {$w_{i-1}$};
  \node (a8) at (-9,0) {$w_1$};

  \draw (a1) -- (a2); 
  \draw (a1) -- (a3);
  \draw (a1) -- (a4); 
  \draw (a1) -- (a5);
  \draw (a1) -- (a6);  
  \draw (a1) -- (a7);
  \draw (a1) -- (a8);
  \draw (a2) -- (a3);
  \draw (a2) -- (a4); 
  \draw (a2) -- (a5);
  \draw (a2) -- (a6);  
  \draw (a2) -- (a7);
  \draw (a2) -- (a8); 
  \draw [dash pattern=on2pt off 2pt] (-5.6,0) -- (-7.8,0);
  \draw [dash pattern=on2pt off 2pt] (5.6,0) -- (7.8,0);
  \end{tikzpicture}
  \caption{Planar embedding of $\mathcal{DNI}(G)$.}
\label{fig:The planar embedding}
\end{center}  
\end{figure}
\end{enumerate}
This concludes the proof.
\end{proof}
\noindent\textbf{\underline{Remark:}} Observe that the set of all finite groups satisfying the conditions given in proposition \ref{Prop: Planar non 2-generated p group structure} is non-empty. For example, again consider the group $Q_8\times C_2$ which is $3$ generated and satisfies the conditions of the proposition.  

Here we first prove the following results to prove one of our main theorem.
\begin{proposition}\label{Prop:DNI_G_contains_K5_G is not 2 genertaed, Z(G)=e}
Let $G$ be a group such that a prime $p \geq 3$ divides $|G|$, and the number of subgroups of order $p$ in $G$ is greater than one. Then the graph $\mathcal{DNI}(G)$ contains an induced subgraph isomorphic to $K_5.$
\end{proposition}
\begin{proof}
It is known that the number of elements of order  $p$ in $G$, denoted $n_p$, satisfies the congruence $n_p \equiv -1 \pmod{p}$. Hence, there are at least four distinct subgroups of order $p$ in $G$ when $p \geq 3$. Let $a_1, a_2, a_3, a_4 \in G$ be elements of order $p$, such that $\langle a_i \rangle \ne \langle a_j \rangle$ for all $i \ne j.$  Since $|\pi(G)| \geq 2$, there exists a prime $q \ne p$  dividing $|G|,$ and hence an element $g \in G$ of order $q$, where $\gcd(p, q) = 1$. By Lemma~\ref{Lemma:Two distinct prime elemnts edge connd in GI(G) complement; G is no 2-generated}, the five elements $a_1, a_2, a_3, a_4, g$ are pairwise adjacent in $\mathcal{DNI}(G)$. Therefore, the induced subgraph on these vertices is isomorphic to $K_5$.
\end{proof}

\begin{corollary}
Let $G$ be a group such that the order of $G$ is not of the form $2^{\alpha}p^{\beta},$ where $p\geq 3,\alpha, \beta\geq 1.$ Then $\mathcal{DNI}(G)$ is not planar.    
\end{corollary}

\begin{lemma}\label{Lemma:DNI(G) has K33 iff, |G|=2mpn, unique subgrp of order p geq 3, exactly 3 subgrps of order 2, not 2-generated}
Let $|G| = 2^{\alpha} p^{\beta},$ where $p\geq 3$ is a prime, and suppose that $G$ has a unique subgroup of order $p$ and exactly $3$ subgroups of order $2$. Then $\mathcal{DNI}(G)$ has an induced subgraph isomorphic to $K_{3,3}$ if and only if the following conditions hold:
\begin{enumerate}
    \item $G$ has an element $a$ of order $2^k$ for some $k \geq 2$, and
    \item $G$ has an element $b$ of order $p^{\ell}$ for some $\ell \geq 2$.
\end{enumerate}
\end{lemma}

\begin{proof}
Suppose $G$ has elements $a$ and $b$ such that $\operatorname{o}(a) = 2^k$ with $k \geq 2$, and $\operatorname{o}(b) = p^{\ell}$ with $\ell \geq 2$. Then $|G_a|\geq 2$ and $|G_b|\geq 2.$ Let $x \in \langle a \rangle$ with $\operatorname{o}(x) = 2$ and $y \in \langle b \rangle$ with $\operatorname{o}(y) = p$. Let $A = \{a, a', x\}$ and $B = \{b, b', y\}$ where $a' \in G_a$ and $b' \in G_b.$ Then the induced subgraph on the vertex set $A \cup B$ is isomorphic to $K_{3,3}$, with partition sets $A$ and $B$.

Conversely, suppose that $G$ has no element of order $4$. So, $G$ has no element of order $d,$ where $4$ divides $d.$ Then every element has order of the form $2 \cdot p^s$ for some $s \geq 0$. Let $\langle h \rangle$ be the unique subgroup of order $p$. Suppose for contradiction that $\mathcal{DNI}(G)$ contains an induced subgraph isomorphic to $K_{3,3}$ with partition sets $X = \{x_1, x_2, x_3\}$ and $Y = \{y_1, y_2, y_3\}$.
Note that either all $x_i$ (for $i \in [3]$) have order divisible only by $2$, or all $y_j$ (for $j \in [3]$) do (in fact, $p$ divides both $\text{o}(x_i)$ and $\text{o}(y_j)$ for some $i, j\in [3]$ implies that $x_i\nsim y_j$ in $\mathcal{DNI}(G)$ since $h\in \langle x\rangle\cap\langle y\rangle$). Suppose the former. Then, since $G$ has no element of order greater than $2$, having $2$ is the only prime divisor, each $x_i$ must have order $2$. But then $\langle x_i \rangle \cap \langle x_j \rangle = \{e\}$ for $i \neq j$, which implies $x_i \sim x_j$ in $\mathcal{DNI}(G)$, contradicting that $X$ is an independent set in the bipartite graph $K_{3,3}$.
A similar contradiction arises if $G$ lacks an element of order $p^{\ell}$ for $\ell \geq 2$. Therefore, both conditions (1) and (2) are necessary for $\mathcal{DNI}(G)$ to contain an induced $K_{3,3}$.
\end{proof}
\begin{lemma}\label{Lemma:DNI(G)_no_K5 |G|=2mpn, exactly 3 2 order subgrps and unique p order subgrps G is not 2 gentreted, Z(G)=e}
Let $|G| = 2^{\alpha}p^{\beta}, p\geq 3$ be a prime, and suppose that $G$ has a unique subgroup of order $p$ and exactly $3$ subgroups of order $2$. Then $\mathcal{DNI}(G)$ has no induced subgraph isomorphic to $K_5$. 
\end{lemma}
\begin{proof}
Suppose, for contradiction, that there exists a subset of vertices $A = \{a_1, a_2, a_3, a_4, a_5\}$ in $\mathcal{DNI}(G)$ that induces a subgraph isomorphic to $K_5$. Assume there exist $a_i, a_j \in A$ such that both orders $\operatorname{o}(a_i)$ and $\operatorname{o}(a_j)$ are divisible by $p$. Since $G$ has a unique subgroup of order $p$, it follows that $\langle a_i \rangle \cap \langle a_j \rangle \neq \{e\}$, and hence $a_i \nsim a_j$ in $\mathcal{DNI}(G)$. This contradicts the assumption that all vertices in $A$ form a complete subgraph. Therefore, at most one element in $A$ can have order divisible by $p$. Without loss of generality, suppose $p \mid \operatorname{o}(a_i)$ for some $i \in [5]$, and that the remaining four elements in $A \setminus \{a_i\}$ all have orders that are powers of $2$. Since $G$ has only three subgroups of order $2$, and $|A \setminus \{a_i\}| = 4$, by the pigeonhole principle, there exist two elements $a_r, a_s \in A \setminus \{a_i\}$ such that $\langle a_r \rangle \cap \langle a_s \rangle \neq \{e\}$. Thus, $a_r \nsim a_s$ in $\mathcal{DNI}(G)$, which contradicts the assumption that $A$ induces a $K_5$. Hence, no such subset $A$ exists, and $\mathcal{DNI}(G)$ has no induced subgraph isomorphic to $K_5$.
\end{proof}
\begin{theorem}\label{Thm:DNI_G_Planar iff when G is not 2 genreted, Z(G)=e}
Let $G$ be a non $2$-generated group such that $|\pi(G)|\geq 2.$ Then $\mathcal{DNI}(G)$ is planar if and only if $|G| = 2^{\alpha} p^{\beta}$ for some $\alpha, \beta \geq 1$ with $n_2=3, n_p=1$, and one of the following holds:
\begin{enumerate}
    \item $G$ does not have any element of order $2^k$ for $k \geq 2$, or
    \item $G$ does not have any element of order $p^{\ell}$ for $\ell \geq 2$.
\end{enumerate}
\end{theorem}
\begin{proof}
We split the proof into the following cases:

\noindent
\textbf{Case 1:} First, suppose that $|\pi(G)|\geq 3.$ If there is a prime $p\geq 3$ such that $n_p>1,$ then by Proposition \ref{Prop:DNI_G_contains_K5_G is not 2 genertaed, Z(G)=e}, $\mathcal{DNI}(G)$ has an induced subgraph isomorphic to $K_{5}.$ If $n_2>1,$ since $|\pi(G)|\geq 3,$ then by Lemma \ref{Lemma:Two distinct prime elemnts edge connd in GI(G) complement; G is no 2-generated} it is easy to see that $\mathcal{DNI}(G)$ has an induced subgraph isomorphic to $K_5.$ 

\noindent
\textbf{Case 2:} Let $|\pi(G)|=2.$ Suppose $|G|=p^{\alpha}q^{\beta}$ and neither $p=2$ nor $q=2.$ Then by Proposition \ref{Prop:DNI_G_contains_K5_G is not 2 genertaed, Z(G)=e}, $\mathcal{DNI}(G)$ contains an induced subgraph isomorphic to $K_5.$ So, in this case $\mathcal{DNI}(G)$ is not planar. Now, let $|G|=2^{\alpha}p^{\beta},$ where $p\geq 3$ be a prime number. Let the number of subgroups of order $p$ in $G$ be $n_p>1.$ Then by Proposition \ref{Prop:DNI_G_contains_K5_G is not 2 genertaed, Z(G)=e}, $\mathcal{DNI}(G)$ has an induced subgraph isomorphic to $K_5.$ If $G$ has unique subgroup of order $p,$ then $n_2\geq 3.$ If $n_2>3,$ then $\mathcal{DNI}(G)$ has induced subgraph $K_5.$ That is, $\mathcal{DNI}(G)$ is not planar. So, finally, we assume that $n_2=3$ and $n_p=1.$ Then by Lemma \ref{Lemma:DNI(G)_no_K5 |G|=2mpn, exactly 3 2 order subgrps and unique p order subgrps G is not 2 gentreted, Z(G)=e},

$\mathcal{DNI}(G)$ does not have any induced subgraph isomorphic to $K_5.$ Moreover, by Lemma \ref{Lemma:DNI(G) has K33 iff, |G|=2mpn, unique subgrp of order p geq 3, exactly 3 subgrps of order 2, not 2-generated}, $\mathcal{DNI}(G)$ has an induced subgraph isomorphic to $K_{3,3}$ if and only if $G$ has elements $a, b$ such that $\text{o}(a)=2^k$ and $\text{o}(b)=p^{\ell},$ for some $k, \ell \geq 2.$ This completes the proof.
\end{proof}

\subsection{2-generated nilpotent groups}
\begin{lemma}\label{Lemma: DNI planarity for 2-generated odd p groups}
Let $G$ be a finite, $2$-generated, $N_I$-group of prime power order, for some odd prime $p$. Then $\mathcal{DNI}(G)$ is planar if and only if $G\cong C_p\times C_p$.   
\end{lemma}
\begin{proof}
Let $G$ be a finite, $2$-generated, $N_I$-group of prime power order, for some odd prime $p$. If $G\cong C_p\times C_p$, then Lemma \ref{Lemma: Diameter abelian p group} says that $\mathcal{DNI}(G)$ is a null graph, hence planar. Otherwise, from Corollary \ref{Cor: 2-generated p group} and Theorem \ref{Theorem: Gorenstein} it can be seen that $G$ contains a proper subgroup of the form $C_p\times C_p$. Now, if $p\geq 5$, then let $x,y$ be generators of $C_p\times C_p$. Then the subgraph induced by $\mathcal{DNI}(G)$ on the set $x,x^2,x^3,y,y^2,y^3$ is isomorphic to $K_{3,3}$. Hence $\mathcal{DNI}(G)$ is not planar. 

Again, if $p=3$, we see that the subgraph induced by $\mathcal{DNI}(G)$ on the proper subgroup isomorphic to $C_3\times C_3$ is isomorphic to $K_{2,2,2,2}$, which is not planar by Observation \ref{K_2222 not planar}. Thus, $\mathcal{DNI}(G)$ is not planar. 
\end{proof}

\begin{lemma}\label{Lemma: DNI planarity for 2-generated non p groups}
Let $G$ be a $2$-generated nilpotent group with $\pi(G)\geq 2$. Then $\mathcal{DNI}(G)$ is not planar. 
\end{lemma}
\begin{proof}
Let $G=H\times K$, where $K$ is a $p$-group for an odd prime, and ${\rm gcd}(\lvert H\rvert,\lvert K\rvert)=1$. Then from Theorem \ref{Thm: direct product}, it can be seen that $\bar{\mathcal{G}_I}(\{e\}\times K)$ is a subgraph of $\mathcal{DNI}(G)$. Also, $\bar{\mathcal{G}_I}(\{e\}\times K)$ is isomorphic to $\bar{\mathcal{G}_I}(K)$. Now, if $K\cong C_p\times C_p$, then the arguments used in the proof of Lemma \ref{Lemma: DNI planarity for 2-generated odd p groups} yield that $\bar{\mathcal{G}_I}(K)$ is not planar. Otherwise, Lemma \ref{Lemma: DNI planarity for 2-generated odd p groups} yields that $\mathcal{DNI}(K)$ is not planar, hence $\bar{\mathcal{G}_I}(K)$ is also not planar. Hence, $\bar{\mathcal{G}_I}(\{e\}\times K)$ is not planar, so is $\mathcal{DNI}(G)$. This concludes the proof.
\end{proof}

\begin{proposition}\label{proposition: DNI planarity for 2-generated 2 groups}
Let $G$ be a $2$-generated group with $\lvert G\rvert=2^k$, for some $k\in\mathbb{N}$. Then $\mathcal{DNI}(G)$ is planar if and only if $G$ is isomorphic to $C_2\times C_2$, $C_4\times C_4$, or $D_8$.    
\end{proposition}
\begin{proof}
$G$ is a finite, $2$-generated, $N_I$-group with $\lvert G\rvert=2^k$, for some $k\in\mathbb{N}$. If $G$ is abelian, then Corollary \ref{Cor: 2-generated p group} implies that $G\cong C_{2^{\frac{k}{2}}}\times C_{2^{\frac{k}{2}}}$. If $\frac{k}{2}\geq 3$, then $G$ has a proper subgroup isomorphic to $C_4\times C_4$. Hence, $\mathcal{DNI}(G)$ contains an induced subgraph isomorphic to $K_{3,3}$. So, $\mathcal{DNI}(G)$ is not planar. Now, if $\frac{k}{2}\leq 2$, then with Lemma \ref{Lemma: Diameter abelian p group} in hand, it can be seen that $\mathcal{DNI}(G)$ is planar.

Now, let $G$ be non-abelian. If $G\cong D_{2^k}$, then $G$ has a proper subgroup isomorphic to $D_{2^{k-1}}$. So, if $k\geq 4$, $D_{8}$ is a proper subgroup of $G$, and the five distinct involutions in $D_8$ form an induced $K_5$ in $\mathcal{DNI}(G)$. So, $\mathcal{DNI}(G)$ is not planar. In addition, if $k\leq 3$, i.e.~when $G\cong D_8$, it can be seen that $E(\mathcal{DNI}(G))$ consists only of the edges between the unique non-trivial central element of $D_8$ and all other non central involutions. So $\mathcal{DNI}(G)$ is planar.  

Let $G$ be a non-abelian $2$ group that is not isomorphic to the dihedral group. Then no pair of involutions can generate $G$, hence any two distinct involutions are adjacent in $\mathcal{DNI}(G)$. So, if $G$ has at least five distinct involutions, then $\mathcal{DNI}(G)$ contains an induced $K_5$, and hence $\mathcal{DNI}(G)$ is not planar. Otherwise, from Lemma \ref{Lemma:number of elements of order p in a p group} and Theorem \ref{Theorem: Gorenstein} it can be concluded that $G$ has exactly three distinct involutions $b_1, b_2, b_3$, and the set $\{e, b_1, b_2, b_3\}$ forms a subgroup isomorphic to $C_2\times C_2$ in $G$. Since $G$ is two generated, let $G=\langle x,y\rangle$. Then $\langle x\rangle\cap\langle y\rangle=\{e\}$. Clearly, both $x,y$ are not involutions. If one of $x,y$ is an involution, without loss of generality, we assume $o(y)=2$. Then $y\in \{b_1,b_2,b_3\}$. Without loss of generality, we assume that $y=b_1$ and $b_2\in \langle x\rangle$. Since, $G$ is non-abelian, we have $[x,y]\neq e$. Also, since $b_2\in \langle x\rangle$ and $\langle x\rangle\cap\langle y\rangle=\{e\}$, we have $x^{-1}yx\neq b_2$. Hence, $x^{-1}yx=b_3=b_2b_1$, i.e.~$yxy^{-1}=y^{-1}xy=xb_2\in \langle x\rangle$. So, $\langle x\rangle$ is normal in $G$, and $G\cong \langle x\rangle\rtimes\langle y\rangle$, with $b_2\in Z(G)$. Now, since $\langle x,y\rangle=G$, we have $\langle x,xy\rangle=G$. But $(xy)^2=xyxy=x^2b_2 \in\langle x\rangle$, as $yxy^{-1}=yxy=xb_2$. So $x^2b_2=e$, since $G$ is a $N_I$-group. Then $x^2=b_2$ and $G\cong D_8$, which contradicts our assumption. So, both $x,y$ have order at least $4$.

If at least one of $x,y$ has an order greater than $4$, let that element be $x$. Let $o(x^l)=4$, for some even integer $l$. Then Lemma \ref{Lemma: langle x^p,y_ranglesubsetneq G} implies that $\langle x^l,y\rangle\subsetneq G$. So, the subgraph induced by $\mathcal{DNI}(G)$ on the vertex set $x^l,x^{2l},x^{3l},y,y^2,y^3$ is isomorphic to $k_{3,3}$. Hence, $\mathcal{DNI}(G)$ is not planar. So, the only case we are left with is when $o(x)=o(y)=4$. If $o(x)=o(y)=4$, then we claim that $o(xy)=4$. Since $G=\langle x,y\rangle$, we have $G=\langle x,xy\rangle=\langle y,xy\rangle$. So by the above arguments, it follows that $o(xy)=4$, and  $\langle x\rangle\cap\langle y\rangle=\langle x\rangle\cap\langle xy\rangle=\langle y\rangle\cap\langle xy\rangle=\{e\}$. Without loss of generality, we assume that $b_1\in\langle x\rangle$, $b_2\in\langle y\rangle$, and $b_3\in\langle xy\rangle$. Then $b_1b_2=b_3$ implies that $x^2y^2=(xy)^2$, i.e. $xy=yx$, which contradicts our assumption of $G$ being non-abelian. Thus the result follows. 
\end{proof}

In light of Lemmas \ref{Lemma: DNI planarity for 2-generated odd p groups}, \ref{Lemma: DNI planarity for 2-generated non p groups} and Proposition \ref{proposition: DNI planarity for 2-generated 2 groups}, we arrive at the following classification result.
\begin{theorem}\label{Theorem: DNI planarity 2-generated nilpotent}
Let $G$ be a $2$-generated nilpotent group. Then $\mathcal{DNI}(G)$ is planar if and only either $G$ is an elementary abelian $p$-group, or $G$ is isomorphic to $C_4\times C_4$ or $D_8$.        
\end{theorem}

\section{Difference graph of finite simple groups}\label{Sec: Difference graph of finite simple groups}
\begin{lemma}\label{lem:prime}
Let $G$ be a non-abelian finite simple group, let $M$ be a maximal subgroup of $G$, and let $e\neq x\in M$. Then there exists an element $u\in M$ of prime order such that $\langle x\rangle\cap\langle u\rangle=\{e\}$.
\end{lemma}

\begin{proof}
Suppose, on the contrary, for every element $u\in M$ of prime order, $\langle x \rangle\cap\langle u \rangle\neq \{e\}$. Since $u$ has prime order, this forces $\langle u \rangle\subseteq\langle x \rangle$.

Fix a prime $p$ dividing $|M|$. Every element of order $p$ in $M$ lies in $\langle x \rangle$, which is cyclic and therefore has a unique subgroup $X_p$ of order $p$. Hence all elements of order $p$ in $M$ lie in $X_p$, and $X_p$ is exactly the subgroup generated by them. In particular $X_p$ is the unique subgroup of order $p$ of $M$. Thus, $X_p$ is normal in $M$, so $M\subseteq \mathscr{N}_G(X_p)$. By maximality of $M$, either $\mathscr{N}_G(X_p)=M$ or $\mathscr{N}_G(X_p)=G$. The latter would make $X_p$ a nontrivial proper normal subgroup of $G$, contradicting the simplicity of $G$. Thus $\mathscr{N}_G(X_p)=M$.

Now, let $g\in G\setminus M$. Then $g\notin \mathscr{N}_G(X_p)$, so $X_p^g\neq X_p$. Note that $X_p^g$ is the unique subgroup of order $p$ of $M^g$. If $M\cap M^g$ contained an element $y$ of order $p$, then $\langle y \rangle$ would be a subgroup of order $p$ of both $M$ and $M^g$, so $\langle y \rangle=X_p$ and $\langle y \rangle=X_p^g$, a contradiction. So $M\cap M^g$ has no element of prime order $p$. As $p$ was an arbitrary prime divisor of $|M|$, we get
\[
M\cap M^g=1\qquad\text{for all } g\in G\setminus M .
\]
Since $1\neq M\subseteq G$, this says that $G$ is a Frobenius group with Frobenius complement $M$. By Frobenius' theorem, $G$ then has a normal subgroup (the Frobenius kernel) of index $|M|$, which is nontrivial and proper. This contradicts the simplicity of $G$. Hence, the result holds.
\end{proof}

\begin{lemma}\label{lem:path}
Let $G$ be a non-abelian finite simple group, $M$ a maximal subgroup of $G$, and $x\neq y$ elements of $M\setminus\{e\}$. Then $x$ and $y$ are joined by a path in $\mathcal{DNI}(G)$.
\end{lemma}

\begin{proof}
By Lemma~\ref{lem:prime} there are $u,v\in M$ of prime order with $\langle x \rangle\cap\langle u \rangle=\{e\}$ and $\langle y \rangle\cap\langle v \rangle=\{e\}$. Since $\langle x,u \rangle\subseteq M\subsetneq G$ and $\langle y,v \rangle\subseteq M\subsetneq G$, the pairs $\{x,u\}$ and $\{y,v\}$ are edges of $\mathcal{DNI}(G).$

If $\langle u \rangle=\langle v \rangle$, then $x\sim u\sim y$ is a path (with $u=v$ or $u\ne v$ generating the same subgroup, so $\langle y \rangle\cap\langle u \rangle=\{e\}$ as well).

If $\langle u \rangle\neq\langle v \rangle$, then $u,v$ have prime order and generate distinct subgroups, so $\langle u \rangle\cap\langle v \rangle=\{e\}$; moreover $\langle u,v \rangle\subseteq M\subsetneq G$. Hence $u\sim v$ and $x\sim u\sim v\sim y$ is a path.
\end{proof}


\begin{theorem}
 Let $G$ be a finite non-Abelian simple group. Then $\mathcal{DNI}(G)$ is connected.   
\end{theorem}

\begin{proof}
Let $a\neq b$ be two non-identity elements of $G$. Now, \cite[Proposition 11]{Lucchini_nemmi} implies that  $\rm{NGen}(G)$ is connected. So, there is a path
\[
a=a_0\sim a_1\sim \cdots \sim a_t=b
\]
in $\rm{NGen}(G)$. Hence, for each $i$, $\langle a_i,a_{i+1} \rangle\subsetneq G$, so $a_i$ and $a_{i+1}$ lie in a common maximal subgroup of $G$. Thus, by Lemma~\ref{lem:path}, $a_i$ and $a_{i+1}$ are joined by a path in $\mathcal{DNI}(G)$. Concatenating these paths gives a path from $a$ to $b$ in $\mathcal{DNI}(G)$. Hence $\mathcal{DNI}(G)$ is connected.
\end{proof}

Next, we study finite non-abelian simple groups whose difference graphs are planar, chordal, and claw-free.

\begin{lemma}\label{even simple group non-planar}
Let $G$ be any finite  non-abelian simple group. Then $\mathcal{DNI}(G)$ is not planar.
\end{lemma}
\begin{proof}
 Let $G$ be any finite non-abelian simple group. Then by \cite[Bemerkung 8.11,  p.~207]{Endliche-gruppen-I-Huppert}, the Sylow
$2$-subgroup $S_2$ of a non-abelian simple group is neither cyclic nor generalized quaternion.
So, according to Theorem \ref{Theorem: Gorenstein}, $S_2$ has a subgroup $H$ isomorphic to $C_2\times C_2$. Since $G$ is simple, $H$ is not normal in $G$. So, there exists $g\in G$ such that $H^g\neq H$. So $\lvert H\cap H^g\rvert=1$ or $2$, and $\lvert H\rvert=\lvert H^g\rvert=4$. Then we have $\lvert H\cup H^g\rvert\geq 6$, i.e.~$H\cup H^g$ has at least $5$ involutions. Since $G$ is not isomorphic to a dihedral group, any two distinct involutions are adjacent in $\mathcal{DNI}(G)$. So any $5$ involution from $H\cup H^g$ forms a induced $K_5$ in $\mathcal{DNI}(G)$. Thus $\mathcal{DNI}(G)$ is not planar.
\end{proof}


\begin{lemma}\label{even simple group non-chordal}
Let $G$ be a finite $2$-generated $N_I$-group having a proper subgroup $H$ such that the order of $H$ is divisible by at least two distinct odd primes. Then $\mathcal{DNI}(G)$ is not chordal.
\end{lemma}
\begin{proof}
Let $H$ be a proper subgroup of $G$ and $p, q$ are two odd primes dividing the order of $H.$ Then both $\phi(p)$ and $\phi(q)$ are greater than equal to $2.$ Let $h_1, h_2, h_3, h_4\in H$ such that the order of $h_1, h_3$ is $p,$ and that of $h_2, h_4$ is $q.$ Since $H$ is proper subgroup and $G$ is $2$-generated, $h_i$ is edge connected with $h_j$ in $\mathcal{DNI}(G)$ if and only if $\langle h_i\rangle\cap\langle h_j\rangle=\{e\}.$ Therefore, $h_1\sim h_2\sim h_3\sim h_4\sim h_1$ is an induced $4$-length cycle in the difference graph $\mathcal{DNI}(G).$ Hence the result. 
\end{proof}

\begin{remark}\label{rem:DNI-not chordal-12-order-subgroup}
Proceeding exactly as in the proof of Lemma~\ref{even simple group non-chordal}, if a group has a proper subgroup of order \(12\), then the difference graph of that group is also non-chordal.  
\end{remark}

\begin{theorem}
Let $G$ be any finite non-abelian simple group. Then $\mathcal{DNI}(G)$ is not chordal.    
\end{theorem}
\begin{proof}
Any finite non-Abelian simple group $G$ has at least one proper subgroup of order either $12$ or $m,$ where $m$ is divisible by at least two distinct odd primes.
\end{proof}

\begin{lemma}\label{lemma:simple-grp-claw-free}
 Let $G$ be a finite $N_I$-group having a proper subgroup $H$ such that the order of $H$ is divided by at least two distinct odd primes. Then $\mathcal{DNI}(G)$ is not claw-free.   
\end{lemma}
\begin{proof}
Let $p, q$ be two distinct odd primes dividing the order of the proper subgroup $H$ of the group $G.$ Let $p\ge 3,$ and $q\ge 5.$ Suppose that $\text{o}(a)=p,$ and $b_1, b_2,$ and $b_3$ are elements of $H$ of order $q.$ Then clearly, $a$ is edge connected with $b_1, b_2, b_3$ and $\{b_1, b_2, b_3\}$ is an independent set. Hence they form a claw in $\mathcal{DNI}(G).$
\end{proof}
\begin{remark}\label{rem:simple-grp-not-claw-free}
From Lemma \ref{lemma:simple-grp-claw-free}, it is easy to see that, if a proper subgroup $H$ has elements of order $2$ and $m$ such that $\gcd(m, 2)=1$ and $\phi(m)\ge 3,$ then the graph $\mathcal{DNI}(G)$ is not claw free.   
\end{remark}

\begin{lemma}\label{lemma:dni-not-claw-free}
Let $G$ be a finite non-abelian simple group having an element of order $4.$ Then $\mathcal{DNI}(G)$ is not claw-free.    
\end{lemma}
\begin{proof}
 $G$ has Sylow $2$-subgroups.
Let $S_2$ be a Sylow $2$-subgroup of $G.$ By \cite[Bemerkung 8.11,  p.~207]{Endliche-gruppen-I-Huppert}, $S_2$ is neither cyclic nor generalized quaternion. Also, any simple group is divisible by either $3$ or $5$ implies that the Sylow $2$-subgroups are proper subgroups of $G.$
Now, Theorem \ref{Theorem: Gorenstein} implies that $S_2$ has a subgroup isomorphic to $C_2\times C_2.$ Also, it is given that $G$ has a cyclic subgroup $H=\{e, h_1, h_2, h_3\},$ where $\text{o}(h_1)=4=\text{o}(h_2),$ and $\text{o}(h_3)=2.$ Let $x\in C_2\times C_2$ such that $\text{o}(x)=2$ and $x\notin H.$ Since $H\subseteq S_2$ and $S_2$ is proper, clearly $x$ is edge connected with $h_1, h_2,$ and $h_3$ in $\mathcal{DNI}(G).$ Note that $\{h_1, h_2, h_3\}$ is an independent set. Hence the induced subgraph with the vertices $x, h_1, h_2, h_3$ is a claw in $\mathcal{DNI}(G).$
\end{proof}

\begin{lemma}\label{lemma:claw-free-proper-dihedral-grp}
Let $G$ be a finite $N_I$-group having a proper subgroup isomorphic to a dihedral group of order $2m\ge 8.$ Then $\mathcal{DNI}(G)$ is not claw-free.     
\end{lemma}
\begin{proof}
 $G$ has a proper subgroup $H$ isomorphic to a dihedral group of order $2m\ge 8.$ So, $H$ has an element $g$ of order $m,$ and  $m\ge 4$ implies that $\langle g\rangle$ has at least $3$ non-identity elements say $g, g_1$ and $g_2,$ where $\langle g\rangle=\langle g_1\rangle.$ Again, in $H$ the number of distinct $2$ ordered elements is greater than equal to $3.$ Let $h\in H$ with $\text{o}(h)=2$ and $h\notin \langle g\rangle.$ Clearly, $h$ is edge connected with $g, g_1$ and $g_2.$ Moreover, $\{g, g_1, g_2\}$ is an independent set. Hence the induced subgraph form by the vertices $\{g, g_1, g_2, h\}.$ This completes the proof. 
\end{proof}

\begin{lemma}\label{lemma:PSL-is-not-claw-free}
 Let \( G = \mathrm{PSL}(2, p^n) \), where \( p \) is a prime, and $n\ge 1.$ Then $\mathcal{DNI}(G)$ is not claw-free.
\end{lemma}
\begin{proof}

Let \( G = \mathrm{PSL}(2, p^n) \), where \( p \) is a prime $n\ge 2$, and let \( k = \gcd(p^n - 1, 2) \). Then by \cite[Satz 8.3, p.~193; Satz 8.4, p.~193]{Endliche-gruppen-I-Huppert} 
the group \( G \) contains a dihedral groups of order $2u$ and \( 2s \), where $u=\frac{p^n-1}{k}$ and \( s = \frac{p^n + 1}{k} \). For any prime $p$ and $n\ge 2,$ the order of at least one dihedral group is greater than equal to $8.$ Therefore, by Lemma \ref{lemma:claw-free-proper-dihedral-grp}, the graph $\mathcal{DNI}(G)$ is not claw free.
Now, let $n=1.$ If $p\geq 7,$ then $G$ has a dihedral group of order greater than equal to $8.$ So, according to Lemma \ref{lemma:claw-free-proper-dihedral-grp}, $\mathcal{DNI}(G)$ is not claw free. Since $\mathrm{PSL}(2, p^n)$ is simple except \(\mathrm{PSL}(2, 2)\) and \(\mathrm{PSL}(2, 3)\), the remaining case is the only group $\mathrm{PSL}(2, 5)\cong A_5.$ The alternating group $A_5$ has a dihedral group of order $10.$ Hence, by Lemma \ref{lemma:claw-free-proper-dihedral-grp}, the graph $\mathcal{DNI}(A_5)$ is not claw free. This competes the proof.

\end{proof}
Now, we consider all the finite non-abelian simple groups having Abelian Sylow $2$-subgroups. The next theorem gives the complete list of all such groups.
\begin{theorem}\cite[Theorem 2]{Gorenstein-CJM},\cite[Theorem I]{Walter-AOM}\label{Thm:cha-simple-grp-abelian-Sylow-2subgrp}
Let $G$ be finite non-Abelian simple group having Abelian Sylow $2$-subgroups. Then $G$ is one of the following:
\begin{enumerate}[(a)]
    \item \(\operatorname{PSL}(2,q)\) with \(q \equiv 3,5 \pmod{8}\);

    \item \(\operatorname{PSL}(2,2^n)\) with \(n \geq 2\);

    \item The Ree groups $R(q)$
   
   \item The sporadic Janko group \(J_1\).
\end{enumerate}

The Sylow \(2\)-subgroups are elementary abelian in all four cases, with orders
\[
2^2,\quad 2^n,\quad 2^3,\quad \text{and } 2^3,
\]
respectively.
\end{theorem}

\begin{lemma}\label{lemma:DNI-Ree-grp-not-claw-free}
 $\mathcal{DNI}(\mathrm{R}(q))$ is not claw free, where $q = 3^{2n+1}$ and $n \geq 1.$  
\end{lemma}

\begin{proof}
According to \cite[Satz 13.2, p. 292]{Finite-grp-III-Huppert-Blackburn}, the centralizer of an involution $j$ in $\mathrm{R}(q)$ is $\langle j\rangle\times \text{PSL}(2, q).$ Hence $\mathrm{R}(q)$ has a subgroup isomorphic to $\text{PSL}(2, q).$ Therefore, by Lemma \ref{lemma:PSL-is-not-claw-free}, $\mathcal{DNI}(G)$ is not claw-free.      
\end{proof}

\begin{theorem}\label{thm:simple-grp-diiference-not-planar}
The difference graph $\mathcal{DNI}(G)$ of any finite non-Abelian simple group is neither planar nor claw-free.
\end{theorem}
\begin{proof}
By Lemma \ref{even simple group non-planar}, the graph $\mathcal{DNI}(G)$ is not planar. 

We now prove that graph $\mathcal{DNI}(G)$ is not claw free for any non-Abelain simple group $G$. First, suppose that the Sylow $2$-subgroups of $G$ are not elementary abelian. Since $G$ is of even order, then $G$ must have an element of order $4.$ Hence, by Lemma \ref{lemma:dni-not-claw-free} is not claw free. Now, let the Sylow $2$-subgroups of $G$ are elementary Abelian. The complete list of such simple groups is depicted in Theorem \ref{Thm:cha-simple-grp-abelian-Sylow-2subgrp}. By Lemmas \ref{lemma:PSL-is-not-claw-free} and \ref{lemma:DNI-Ree-grp-not-claw-free}, the difference graphs of the groups described in (a), (b) and (c) are not claw free. Now, sporadic Janko group \(J_1\) has an involution whose centralizer is isomorphic to the direct product of the group of order two and the alternating group $A_5$
of order $60$. That is, $J_1$ has a proper subgroup $A_5$ such that the order of
$A_5$ is divided by at least two distinct odd primes. So, by Lemma \ref{lemma:simple-grp-claw-free}, $\mathcal{DNI}(J_1)$ is not claw-free. 
\end{proof}

\section*{Acknowledgements} 
The first-named gratefully acknowledges the financial support provided by the NBHM research project (Reference No.~02011/29/2025NBHM(RP)/R\&DII/11951). He sincerely thanks the National Board for Higher Mathematics (NBHM), India, for funding this research. He also acknowledges the excellent research environment and facilities provided by Dhirubhai Ambani University, Gandhinagar, Gujarat. The second-named author thanks Harish-Chandra Research Institute, India, for providing funding (in the form of post-doctoral fellowship) during this work. The authors thank Prof. Andrea Lucchini for sharing his ideas regarding the classification of metacyclic $N_I$-groups, and connectedness of the difference graphs of non-abelian simple groups. Moreover, the authors thank Prof. Peter Cameron for suggesting the problem that motivated this work.

\subsection*{Data Availability Statements}
Data sharing not applicable to this article as no datasets were generated or analysed during the current study.

\subsection*{Competing Interests} The authors have no competing interests to declare that are relevant to the content of this article.

\bibliographystyle{amsplain}
\bibliography{gen-inv-lcp.bib}

\end{document}